\documentclass[a4paper,12pt,amsart,frenchb]{article}

\usepackage{amsmath,amsbsy,amsfonts,amssymb}
\usepackage[french]{babel}
\newcommand{\ass}[2]{\vskip0.3cm\noindent
{\bf {#1}}. { \sl {#2}}\vskip0.3cm\noindent
}

\begin{document}

  \title{ Repr\'esentations  de r\'eduction unipotente pour $SO(2n+1)$, I: une involution}
\author{J.-L. Waldspurger}
\date{21 novembre 2016}
\maketitle

\bigskip

{\bf Introduction}

Soient $p$ un nombre premier et $F$ une extension finie de ${\mathbb Q}_{p}$. Soit $n\geq1$ un entier. On suppose que $p$ est grand relativement \`a $n$. On consid\`ere le groupe sp\'ecial orthogonal d'une forme quadratique non d\'eg\'en\'er\'ee  sur un espace de dimension $2n+1$ sur $F$. Plus pr\'ecis\'ement, on consid\`ere les deux formes possibles de ce groupe: la forme d\'eploy\'ee que l'on note $G_{iso}$ et la forme non quasi-d\'eploy\'ee $G_{an}$. Pour l'un ou l'autre  de ces indices $\sharp=iso$ ou $an$, notons $Irr_{unip,\sharp}$ l'ensemble des (classes d'isomorphismes de) repr\'esentations admissibles irr\'eductibles de $G_{\sharp}(F)$ qui sont de r\'eduction unipotente. Cette derni\`ere propri\'et\'e est d\'efinie de la fa\c{c}on suivante. Soit $\pi$ une repr\'esentation admissible irr\'eductible de $G_{\sharp}(F)$ dans un espace complexe $E$. Pour tout sous-groupe parahorique $K$ de $G_{\sharp}(F)$, notons $K^{u}$ son radical pro-$p$-unipotent et $E^{K^{u}}$ le sous-espace des \'el\'ements de $E$ fix\'es par $K^{u}$. De $\pi$ se d\'eduit une repr\'esentation de $K/K^{u}$ dans $E^{K^{u}}$. Le groupe $K/K^{u}$ s'identifie au groupe des points sur le corps r\'esiduel ${\mathbb F}_{q}$ de $F$ d'un groupe alg\'ebrique connexe d\'efini sur ${\mathbb F}_{q}$. En suivant Lusztig, on sait d\'efinir la notion de repr\'esentation unipotente d'un tel groupe. On dit que $\pi$ est de r\'eduction unipotente si et seulement s'il existe $K$ comme ci-dessus de sorte que $E^{K^{u}}$ soit non nul et que la repr\'esentation de $K/K^{u}$ dans $E^{K^{u}}$ soit unipotente. On note $Irr_{tunip,\sharp}$ le sous-ensemble des repr\'esentations dans  $Irr_{unip,\sharp}$ qui sont temp\'er\'ees. L'involution introduite par  Zelevinsky dans le cas des groupes $GL(n)$ a \'et\'e g\'en\'eralis\'ee aux autres groupes par Aubert et par Schneider et
Stuhler, cf. \cite{Au} et \cite{SS}. Notons-la $D$. Dans trois articles dont celui-ci est le premier, nous allons \'etudier les repr\'esentations appartenant \`a $Irr_{tunip,\sharp}$ ainsi que leurs images par l'involution $D$.  Notre but est de calculer le front d'ondes de $D(\pi)$ pour $\pi\in Irr_{tunip, \sharp}$. Cela sera fait dans le troisi\`eme article et nous y reviendrons le moment venu. D\'ecrivons plut\^ot le contenu des deux premiers.

On note $Irr_{tunip}$ la r\'eunion disjointe de $Irr_{tunip,iso}$ et $Irr_{tunip,an}$. La conjecture de Langlands, raffin\'ee par Deligne et Lusztig, param\'etrise $Irr_{tunip}$ de la fa\c{c}on suivante. Notons $W_{F}$ de groupe de Weyl de $F$ et $Sp(2n;{\mathbb C})$ le groupe symplectique complexe d'un espace de dimension $2n$. Pour un homomorphisme $\psi:W_{F}\times SL(2;{\mathbb C})\to Sp(2n;{\mathbb C})$, notons $S(\psi)$ le groupe des composantes connexes du centralisateur dans $Sp(2n;{\mathbb C})$ de l'image de $\psi$. C'est un produit fini de groupes ${\mathbb Z}/2{\mathbb Z}$ et on note $S(\psi)^{\wedge}$ son groupe de caract\`eres. Alors $Irr_{tunip}$ est conjecturalement param\'etr\'e par les classes de conjugaison par $Sp(2n;{\mathbb C})$ de couples $(\psi,\epsilon)$ o\`u

$\psi:W_{F}\times SL(2;{\mathbb C})\to Sp(2n;{\mathbb C})$ est un homomorphisme dont la restriction \`a $W_{F}$ est  non-ramifi\'ee, semi-simple, d'image born\'ee et dont la restriction \`a $SL(2;{\mathbb C})$ est alg\'ebrique;

$\epsilon\in S(\psi)^{\wedge}$.

Ce param\'etrage devient unique si on impose aux repr\'esentations ainsi param\'etr\'ees des conditions relatives \`a l'endoscopie et \`a l'endoscopie tordue. Nous y reviendrons ci-dessous. 

Le param\'etrage a \'et\'e obtenu par diff\'erents auteurs: Lusztig, cf.  \cite{L}; Moeglin, cf. \cite{M} th\'eor\`eme 5.2; Arthur, cf. \cite{Ar} th\'eor\`eme 2.2.1. Arthur a prouv\'e les propri\'et\'es relatives \`a l'endoscopie, du moins dans le cas du groupe $G_{iso}(F)$. Nous utilisons les constructions de Lusztig. Le but des deux premiers articles est de prouver que les repr\'esentations qu'il a construites  v\'erifient bel et bien
  les propri\'et\'es requises quant \`a l'endoscopie. 
  
Pour \'enoncer ces propri\'et\'es, il est commode de modifier le param\'etrage. Consid\'erons un  homomorphisme $\psi$ comme ci-dessus. Puisqu'il est non ramifi\'e, il est d\'etermin\'e par sa restriction $\rho$ \`a $SL(2;{\mathbb C})$ et par l'image $s$ d'un \'el\'ement de Frobenius de $W_{F}$.  On sait param\'etrer les classes de conjugaison d'homomorphismes alg\'ebriques $\rho:SL(2;{\mathbb C})\to Sp(2n;{\mathbb C})$ par les orbites unipotentes de $Sp(2n;{\mathbb C})$. Celles-ci sont elles-m\^emes param\'etr\'ees par l'ensemble ${\cal P}^{symp}(2n)$ des partitions symplectiques de $2n$, cf. 1.3.  Ainsi, pour tout $\lambda\in {\cal P}^{symp}(2n)$, fixons un homomorphisme $\rho_{\lambda}$ dans la classe param\'etr\'ee par $\lambda$. Notons $Z(\lambda)$ le commutant dans $Sp(2n;{\mathbb C})$ de l'image de $\rho_{\lambda}$. L'\'el\'ement $s$ doit \^etre un \'el\'ement semi-simple de $Z(\lambda)$. De plus, ses valeurs propres doivent \^etre de module $1$: c'est la condition "temp\'er\'ee" et on dira simplement que $s$ est "compact". On note $Z(\lambda,s)$ le commutant de $s$ dans $Z(\lambda)$,  ${\bf Z}(\lambda,s)$ le groupe des composantes connexes de $Z(\lambda,s)$ et ${\bf Z}(\lambda,s)^{\wedge}$ son groupe de caract\`eres.  On a les \'egalit\'es $S(\psi)={\bf Z}(\lambda,s)$ et $S(\psi)^{\wedge}={\bf Z}(\lambda,s)^{\wedge}$. Ainsi, l'ensemble $Irr_{tunip}$ est param\'etr\'e par l'ensemble des triplets $(\lambda,s,\epsilon)$, o\`u $\lambda\in {\cal P}^{symp}(2n)$, $s$ est un \'el\'ement semi-simple et compact de $Z(\lambda)$ et $\epsilon\in {\bf Z}(s,\lambda)^{\wedge}$ (plus exactement par l'ensemble des classes de conjugaison de tels triplets, en un sens facile \`a pr\'eciser). Nous notons $\pi(\lambda,s,\epsilon)$ l'\'el\'ement que Lusztig associe \`a un tel triplet.

{\bf Remarque.} En fait,   Lusztig  construit non pas les repr\'esentations temp\'er\'ees mais leurs  images par l'involution $D$;  notre $\pi(\lambda,s,\epsilon)$ est l'image par cette involution de la repr\'esentation construite par Lusztig. 

\bigskip

Pour tout ensemble $X$, notons ${\mathbb C}[X]$ l'espace vectoriel complexe de base $X$.  L'espace ${\mathbb C}[Irr_{tunip}]$ est somme directe de ${\mathbb C}[Irr_{tunip,iso}]$ et ${\mathbb C}[Irr_{tunip,an}]$. Pour $\pi\in {\mathbb C}[Irr_{tunip}]$, on note $\pi_{iso}$ et $\pi_{an}$ ses deux composantes. Introduisons l'ensemble $\mathfrak{Endo}_{tunip}$ des classes de conjugaison de triplets $(\lambda,s,h)$ o\`u $\lambda\in {\cal P}^{symp}(2n)$, $s$ et $h$ sont des \'el\'ements semi-simples de $Z(\lambda)$, $s$ et $h$ commutent entre eux, $s$ est compact et $h^2=1$. Pour un tel triplet, on a $h\in Z(\lambda,s)$ et, pour $\epsilon\in {\bf Z}(\lambda,s)^{\wedge}$, on peut \'evaluer $\epsilon$ en l'image de $h$ dans ${\bf Z}(\lambda,s)$. On note simplement $\epsilon(h)$ cette valeur. Posons
$$\Pi(\lambda,s,h)=\sum_{\epsilon\in {\bf Z}(\lambda,s)^{\wedge}}\pi(\lambda,s,\epsilon)\epsilon(h).$$
On note $\mathfrak{St}_{tunip}$ le sous-ensemble des $(\lambda,s,h)\in \mathfrak{Endo}_{tunip}$ tels que $h=1$.  On sait d\'efinir la notion de distribution stable pour le groupe $G_{iso}(F)$. On dit qu'un \'el\'ement de ${\mathbb C}[Irr_{tunip,iso}]$ est stable si sa distribution trace associ\'ee l'est. La premi\`ere propri\'et\'e caract\'erisant le param\'etrage est

(1) pour $(\lambda,s,1)\in \mathfrak{St}_{tunip}$, $\Pi_{iso}(\lambda,s,1)$ est stable.

La deuxi\`eme propri\'et\'e fait intervenir  l'endoscopie tordue qui relie $G_{iso}$ et le groupe $GL(2n)$, ou plus exactement un espace tordu sur ce groupe. Nous l'\'enoncerons pr\'ecis\'ement en 2.1 et la r\'esumons ici par l'assertion vague

(2) pour $(\lambda,s,1)\in \mathfrak{St}_{tunip}$, le transfert de $\Pi_{iso}(\lambda,s,1)$ \`a $GL(2n)$ par endoscopie tordue est une repr\'esentation bien d\'etermin\'ee de ce groupe.

Consid\'erons un couple $(n_{1},n_{2})\in {\mathbb N}^2$ tel que $n_{1}+n_{2}=n$. On note $G_{n_{1},iso}$ et $G_{n_{2},iso}$ les groupes similaires \`a $G_{iso}$ quand on remplace $n$ par $n_{1}$ ou $n_{2}$. De m\^eme, on affecte les objets introduits ci-dessus d'un indice $n_{1}$ ou $n_{2}$ quand ils sont relatifs \`a ces entiers. Consid\'erons $(\lambda_{1},s_{1},1)\in \mathfrak{St}_{tunip,n_{1}}$ et $(\lambda_{2},s_{2},1)\in \mathfrak{St}_{tunip,n_{2}}$. On en d\'eduit un triplet $(\lambda,s,h)\in \mathfrak{Endo}_{tunip}$ de la fa\c{c}on suivante. Le groupe $Sp(2n_{1};{\mathbb C})\times Sp(2n_{2};{\mathbb C})$ se plonge naturellement dans $Sp(2n;{\mathbb C})$. Par composition avec ce plongement, $\rho_{\lambda_{1}}\otimes \rho_{\lambda_{2}}$ devient un homomorphisme de $SL(2;{\mathbb C})$ dans $Sp(2n;{\mathbb C})$, qui est param\'etr\'e par la partition $\lambda=\lambda_{1}\cup \lambda_{2}$. L'\'el\'ement $(s_{1},s_{2})$ devient pareillement un \'el\'ement $s$ de $Sp(2n;{\mathbb C})$. Enfin, $h$ est l'image par le plongement du produit de l'identit\'e de $Sp(2n_{1};{\mathbb C})$ et de l'\'el\'ement $-1$ de $Sp(2n_{2};{\mathbb C})$. Le groupe $G_{n_{1},iso}\times G_{n_{2},iso}$ est le groupe endoscopique d'une donn\'ee endoscopique \'evidente de $G_{iso}$ ou de $G_{an}$. On sait d\'efinir le transfert \`a $G_{iso}(F)$ ou $G_{an}(F)$ d'une distribution stable sur $G_{n_{1},iso}(F)\times G_{n_{2},iso}(F)$. 

{\bf Remarque.} Pour travailler simultan\'ement avec les deux groupes $G_{iso}$ et $G_{an}$, on doit utiliser une variante un peu sophistiqu\'ee de l'endoscopie, cf. 2.1.  On peut parler d'endoscopie pour les formes int\'erieures pures.

\bigskip
 La troisi\`eme propri\'et\'e du param\'etrage est

(3)(a) $\Pi_{iso}(\lambda,s,h)$ est le transfert \`a $G_{iso}(F)$ de $\Pi_{iso}(\lambda_{1},s_{1},1)\otimes \Pi_{iso}(\lambda_{2},s_{2},1)$;

(3)(a) $-\Pi_{an}(\lambda,s,h)$ est le transfert \`a $G_{an}(F)$ de $\Pi_{iso}(\lambda_{1},s_{1},1)\otimes \Pi_{iso}(\lambda_{2},s_{2},1)$.

On souligne la pr\'esence du signe $-1$ dans (3)(b). Les trois propri\'et\'es ci-dessus d\'eterminent enti\`erement le param\'etrage. La propri\'et\'e (1) est connue: c'est le r\'esultat principal de \cite{MW}. La propri\'et\'e (2) l'est aussi, cf. \cite{W3}.  Nous \'enoncerons plus pr\'ecis\'ement les propri\'et\'es (3)(a) et (3)(b) en 2.1 ci-dessous et nous les d\'emontrerons
dans le deuxi\`eme article.  

La premi\`ere section du pr\'esent article introduit tout le mat\'eriel qui sera utilis\'e dans la suite. On reprend en particulier de nombreuses constructions faites dans \cite{MW}. La deuxi\`eme section est consacr\'ee \`a l'une de ces constructions, celle d'une certaine involution. Rappelons d'o\`u vient celle-ci. Introduisons le sous-ensemble $\mathfrak{Endo}_{unip-disc}$ des $(\lambda,s,h)\in \mathfrak{Endo}_{tunip}$ tels que $s^2=1$ et que le commutant commun $Z(\lambda,s)\cap Z(\lambda,h)$ 
de $s$ et $h$ dans $Z(\lambda)$ soit fini.  La famille des $\Pi(\lambda,s,h)$, quand $(\lambda,s,h)$ d\'ecrit $\mathfrak{Endo}_{unip,disc}$, est une base du sous-espace de ${\mathbb C}[Irr_{tunip}]$ engendr\'e par les repr\'esentations de r\'eduction unipotente qui sont elliptiques au sens d'Arthur, cf. \cite{Ar2}. L'espace ${\mathbb C}[\mathfrak{Endo}_{unip,disc}]$ s'identifie donc \`a ce sous-espace de ${\mathbb C}[Irr_{tunip}]$. Une compatibilit\'e \`a l'induction permet de ramener la preuve de (3)(a) et (3)(b) au cas o\`u $(\lambda,s,h)\in \mathfrak{Endo}_{unip,disc}$. Pour $\sharp=iso$ ou $an$, on doit calculer 
 autant que faire se peut le caract\`ere de $\Pi_{\sharp}(\lambda,s,h)$   
  sur les \'el\'ements fortement r\'eguliers de $G_{\sharp}(F)$. Gr\^ace \`a un r\'esultat d'Arthur, on peut se restreindre aux \'el\'ements elliptiques de $G_{\sharp}(F)$. Un tel \'el\'ement appartient \`a un sous-groupe compact et  la premi\`ere chose \`a faire est de restreindre nos repr\'esentations aux divers sous-groupes compacts maximaux de $G_{\sharp}(F)$. La construction de Lusztig est parfaitement adapt\'ee pour cela. Mais elle d\'ecrit ces restrictions en termes de repr\'esentations irr\'eductibles des groupes "r\'esiduels" (les groupes $K/K^{u}$ ci-dessus). Le calcul des caract\`eres de ces repr\'esentations n'est pas simple mais a \'et\'e effectu\'e par Lusztig. Pour cela, il a introduit les "faisceaux-caract\`eres" qui cr\'eent des fonctions-traces plus facilement calculables que les caract\`eres de repr\'esentations.  Les faisceaux-caract\`eres et les repr\'esentations irr\'eductibles d'un groupe r\'esiduel sont param\'etr\'es par un m\^eme ensemble combinatoire, appelons-le $X$. On a ainsi deux applications lin\'eaires $k$ et $Rep$ d\'efinies sur ${\mathbb C}[X]$: pour $x\in X$, $k(x)$ est la fonction-trace du faisceau-caract\`ere param\'etr\'e par $x$ et $Rep(x)$ est le caract\`ere de la repr\'esentation irr\'eductible  param\'etr\'ee par $x$. Lusztig a d\'efini une involution ${\cal F}^L$ de ${\mathbb C}[X]$ qui v\'erifie $k=Rep\circ {\cal F}^L$. Revenons \`a notre probl\`eme. Apr\`es avoir restreint nos repr\'esentations aux sous-groupes  compacts maximaux de $G_{\sharp}(F)$, on doit appliquer les involutions ${\cal F}^L$ relatives \`a ces sous-groupes. Le calcul para\^{\i}t bien compliqu\'e mais on l'a simplifi\'e dans \cite{MW} en introduisant une involution ${\cal F}$  de ${\mathbb C}[\mathfrak{Endo}_{unip,disc}]$ qui a la propri\'et\'e suivante. Soit $(\lambda,s,h)\in \mathfrak{Endo}_{unip,disc}$. On en d\'eduit une repr\'esentation du groupe r\'esiduel ci-dessus $K/K^{u}$. Soit $\varphi\in {\mathbb C}[X]$ tel que cette repr\'esentation (identifi\'ee \`a son caract\`ere)  soit $Rep(\varphi)$. Alors la repr\'esentation du m\^eme  groupe r\'esiduel associ\'ee \`a ${\cal F}(\Pi(\lambda,s,h))$ est \'egale \`a $Rep\circ {\cal F}^L(\varphi)$, c'est-\`a-dire \`a $k(\varphi)$.  Cette involution ${\cal F}$ rend possible la suite du calcul.

  {\bf Remarque.} En fait, on a d\'emontr\'e dans \cite{MW} une propri\'et\'e plus faible mais suffisante, \`a savoir que la projection cuspidale de la repr\'esentation du groupe r\'esiduel associ\'ee \`a ${\cal F}(\Pi(\lambda,s,h))$ est \'egale \`a la projection cuspidale de $Rep\circ {\cal F}^L(\varphi)$.
  
  \bigskip
  
Dans \cite{MW}, on a donn\'e une d\'efinition combinatoire assez compliqu\'ee de ${\cal F}$. Avec les notations ci-dessus, elle est tr\`es simple: elle consiste \`a \'echanger $s$ et $h$. En effet, pour $(\lambda,s,h)\in \mathfrak{Endo}_{unip,disc}$, on a aussi  $(\lambda,h,s)\in \mathfrak{Endo}_{unip,disc}$ et on montrera en 2.3 que ${\cal F}(\Pi(\lambda,s,h))=\Pi(\lambda,h,s)$.   Enfin, en 2.7,  nous l\`everons la restriction \'evoqu\'ee dans la remarque ci-dessus: on a bien \'egalit\'e entre les deux repr\'esentations de cette remarque et pas seulement de leurs projections cuspidales. Le r\'esultat de 2.7 est conditionnel: on admet les propri\'et\'es (3)(a) et (3)(b). Mais, comme on l'a dit, ces propri\'et\'es seront d\'emontr\'ees dans l'article suivant. 
 
{\bf Remarque.} A l'aide des constructions de Lusztig, on devrait pouvoir traiter non seulement le cas des param\`etres non ramifi\'es, mais celui des param\`etres mod\'er\'ement ramifi\'es.  J'ignore ce que devient notre involution ${\cal F}$ dans cette situation plus g\'en\'erale.

\bigskip

Un index des notations se trouve en fin de l'article.

 \bigskip 

\section{Les groupes et leurs repr\'esentations}

\subsection{Groupes orthogonaux}
Soit $F$ un corps local non-archim\'edien de caract\'eristique nulle. On note $\mathfrak{o}$ son anneau d'entiers
  et on fixe une uniformisante $\varpi$. On note ${\mathbb F}_{q}=\mathfrak{o}/\varpi\mathfrak{o}$ le corps r\'esiduel, $q$ \'etant son nombre d'\'el\'ements;  $p$ la caract\'eristique de ${\mathbb F}_{q}$; $\vert .\vert _{F}$ la valeur absolue de $F$; $val_{F}$ la valuation. On a $\vert \varpi\vert _{F}=q^{-1}$ et $val_{F}(\varpi)=1$.  On suppose $p\not=2$. On note $F^{\times}$ le groupe multiplicatif $F-\{0\}$;  $\mathfrak{o}^{\times}$ le groupe des unit\'es;  $F^{\times2} $,  $\mathfrak{o}^{\times2}$ et ${\mathbb F}_{q}^{\times2} $ les sous-groupes des carr\'es dans $F^{\times}$, $\mathfrak{o}^{\times}$ et ${\mathbb F}_{q}^{\times}$.

Soit $d\geq1$ un entier et soit $V$ un espace vectoriel sur $F$ de dimension $d$, muni d'une forme bilin\'eaire sym\'etrique et non d\'eg\'en\'er\'ee $Q$. Le  d\'eterminant $det(Q)$ est bien d\'efini dans $F^{\times}/F^{\times 2}$. On pose $\eta(Q)=(-1)^{[d/2]}det(Q)$, o\`u, pour $x\in {\mathbb R}$, $[x]$ est la partie enti\`ere de $x$. La valuation $val_{F}(\eta(Q))$ est bien d\'efinie dans ${\mathbb Z}/2{\mathbb Z}$.  On note $O(Q)$ le groupe orthogonal de $(V,Q)$, $SO(Q)$ ou $O^+(Q)$ le groupe sp\'ecial orthogonal et $O^-(Q)$ la composante connexe non neutre de $O(Q)$. 

 Les m\^emes d\'efinitions s'appliquent si $V$ est un espace sur ${\mathbb F}_{q}$. La terme $\eta(Q)$ est alors un \'el\'ement de ${\mathbb F}^{\times}/{\mathbb F}^{\times 2}$. On identifie ce groupe \`a $\mathfrak{o}^{\times}/\mathfrak{o}^{\times 2}$.  Pour mieux distinguer les corps de base, on notera par lettres grasses ${\bf SO}(Q)$ etc.. les groupes relatifs aux espaces d\'efinis sur ${\mathbb F}_{q}$.

\bigskip
Revenons \`a un espace quadratique $(V,Q)$ d\'efini sur $F$. Pour un $\mathfrak{o}$-r\'eseau $L\subset V$, notons $L^*=\{v\in V;\forall v'\in L, Q(v,v')\in \mathfrak{o}\}$. Il existe des r\'eseaux $L$ presque autoduaux, c'est-\`a-dire tels que $\varpi L^*\subset L\subset L^*$. Pour un tel r\'eseau, posons $l'=L/\varpi L^*$, $l''=L^*/L$. On note $d'$, resp. $d''$, la dimension sur ${\mathbb F}_{q}$ de $l'$, resp. $l''$. La forme $Q$ se r\'eduit en une forme quadratique non d\'eg\'en\'er\'ee $Q'$ sur $l'$ et la forme $\varpi Q$ se r\'eduit en une forme quadratique non d\'eg\'en\'er\'ee $Q''$ sur $l''$. Les termes $\eta(Q')$ et $\eta(Q'')$ sont ind\'ependants  du choix du r\'eseau $L$. On les note $\eta'(Q)$ et $\eta''(Q)$. 
On v\'erifie les relations
$$d'\equiv d+val_{F}(\eta(Q))\,\,mod\,\, 2{\mathbb Z},\,\,d''\equiv val_{F}(\eta(Q))\,\,mod\,\, 2{\mathbb Z},$$
$$(-1)^{d'd''}\eta'(Q)\eta''(Q)\varpi^{val_{F}(\eta(Q))}=\eta(Q).$$

Consid\'erons maintenant un entier $d\geq1$ et un \'el\'ement $\eta\in F^{\times}/F^{\times 2}$. A l'exception des cas $d=1$ et $(d,\eta)=(2,1)$, on sait qu'il y a deux classes d'isomorphie d'espaces quadratiques $(V,Q)$ comme ci-dessus tels que $\eta(Q)=\eta$. On les note $(V_{iso},Q_{iso})$ et $(V_{an},Q_{an})$, les indices $iso$ et $an$ \'etant d\'etermin\'es par les relations suivantes

si $d+val_{F}(\eta)$ est pair, $\eta''(Q_{iso})\in {\mathbb F}^{\times 2}$ et $\eta''(Q_{an})\not\in {\mathbb F}_{q}^{\times 2}$;

si $d+val_{F}(\eta)$ est impair, $\eta'(Q_{iso})\in  {\mathbb F}^{\times 2}$ et $\eta'(Q_{an})\not\in {\mathbb F}_{q}^{\times 2}$.

 Dans les cas particuliers $d=1$ ou $(d,\eta)=(2,1)$, la d\'efinition ci-dessus conduit \`a consid\'erer l'unique espace quadratique $(V,Q)$ comme \'etant $(V_{iso},Q_{iso})$. 
 
Le groupe $SO(Q_{iso})$ est quasi-d\'eploy\'e. Le groupe $SO(Q_{an})$ en est une forme int\'erieure.  
 Il n'est pas quasi-d\'eploy\'e si $d$ est impair ou si $d$ est pair et $\eta=1$. Il est isomorphe \`a $SO(Q_{iso})$ si $d$ est pair et $\eta\not=1$.
 
 Consid\'erons maintenant le cas des espaces quadratiques d\'efinis sur ${\mathbb F}_{q}$. Dans ce cas, $d$ \'etant fix\'e, la classe d'isomorphie de $(V,Q)$ est d\'etermin\'ee par $\eta(Q)$. Si $d$ est impair, le groupe $O(Q)$ est ind\'ependant de $\eta$. On le note simplement ${\bf O}(d)$ et on note ses deux composantes connexes ${\bf SO}(d)$, ou ${\bf O}^+(d)$, et ${\bf O}^-(d)$. Le groupe ${\bf SO}(d)$ est d\'eploy\'e. Si $d$ est pair, on note ${\bf O}(d)_{iso}$ le groupe orthogonal de l'espace $(V,Q)$ tel que $\eta(Q)\in {\mathbb F}_{q}^{\times2}$ et ${\bf O}(d)_{an}$ celui de l'espace $(V,Q)$ tel que $\eta(Q)\not\in {\mathbb F}_{q}^{\times 2}$. Pour un indice $\sharp=iso$ ou $an$, on note aussi ${\bf SO}(d)_{\sharp}$, ou ${\bf O}^+(d)_{\sharp}$,  et ${\bf O}^-(d)_{\sharp}$ les deux composantes connexes de ${\bf O}(d)_{\sharp}$. Le groupe ${\bf SO}(d)_{iso}$ est d\'eploy\'e et le groupe ${\bf SO}(d)_{an}$ ne l'est pas.

Dans cet article, nous fixons un entier $n\geq1$. On suppose

$p>6n+4$.

Cette borne est reprise de \cite{MW}. Nous consid\'erons la construction ci-dessus (sur $F$) pour $d=2n+1$ et $\eta=1$. On  a donc deux couples $(V_{iso},Q_{iso})$ et $(V_{an},Q_{an})$ d\'efinis sur $F$. Concr\`etement, fixons un \'el\'ement $\xi\in \mathfrak{o}^{\times}-\mathfrak{o}^{\times 2}$. Pour  un indice $\sharp=iso$ ou $an$, il y a une base $\{v_{1},...,v_{2n+1}\}$ de $V_{\sharp}$ telle qu'en notant $v=\sum_{i=1,...,2n+1}x_{i}v_{i}$ la d\'ecomposition d'un \'el\'ement $v\in V_{\sharp}$ dans cette base, on ait

si $\sharp=iso$, $Q_{iso}(v,v)=x_{n+1}^2+2\sum_{i=1,...,n}x_{i}x_{2n+2-i}$,

si $\sharp=an$, $Q_{an}(v,v)=\varpi x_{1}^2+\xi x_{n+1}^2-\xi \varpi^{-1}x_{2n+1}^2+2\sum_{i=2,...,n}x_{i}x_{2n+2-i}$.

On pose simplement $G_{\sharp}=SO(Q_{\sharp})$.

 Nous introduirons divers objets relatifs \`a notre entier  $n$  fix\'e. On aura parfois besoin des objets analogues relatifs \`a d'autres entiers.  On les notera simplement en ajoutant l'entier en question dans la notation. Par exemple $G_{\sharp}$ peut aussi bien \^etre not\'e $G_{n,\sharp}$.

\bigskip

\subsection{Sous-groupes parahoriques maximaux}

    Notons $D(n)$ l'ensemble des couples $(n',n'')\in {\mathbb N}^2$ tels que $n'+n''=n$. Posons $D_{iso}(n)=D(n)$ et $D_{an}(n)=\{(n',n'')\in D(n); n''\geq1\}$. Soit $\sharp=iso$ ou $an$.  Pour $ (n',n'')\in D_{\sharp}(n)$, on d\'efinit le r\'eseau $L_{n',n''}\subset V_{\sharp}$ par
$$L_{n',n''}=\mathfrak{o}v_{1}\oplus...\oplus\mathfrak{o}v_{2n+1-n''}\oplus \varpi\mathfrak{o}v_{2n+2-n''}\oplus...\oplus \varpi\mathfrak{o}v_{2n+1},$$
o\`u $\{v_{1},...,v_{2n+1}\}$ est la base introduite dans le paragraphe pr\'ec\'edent. 
 On a
$$L_{n',n''}^*=\varpi^{-1}\mathfrak{o}v_{1}\oplus...\oplus \varpi^{-1}\mathfrak{o}v_{n''}\oplus \mathfrak{o}v_{i+1}\oplus...\oplus \mathfrak{o}v_{2n+1}$$
et les inclusions
$$\varpi L_{n',n''}^*\subset L_{n',n''}\subset L_{n',n''}^*.$$
On pose $l_{2n'+1}=L_{n',n''}/\varpi L_{n',n''}^*$, $l_{2n''}=L_{n',n''}^*/L_{n',n''}$. Comme on l'a dit en 1.1, ces espaces sur ${\mathbb F}_{q}$ sont munis de formes bilin\'eaires sym\'etriques et non d\'eg\'en\'er\'ees. Leurs groupes orthogonaux sont respectivement ${\bf O}(2n'+1)$ et ${\bf O}(2n'')_{\sharp}$.

{\bf Remarque.} Si $\sharp=iso$ et $n''=0$, $l_{0}=\{0\}$. Dans ce cas, on supprime les constructions relatives \`a cet espace.
\bigskip

 On note $K_{n',n''}^{\pm}$ le sous-groupe des $g\in G_{\sharp}(F)$ tels que $g(L_{n',n''})\subset L_{n',n''}$ (ce qui entra\^{\i}ne aussi $g(L_{n',n''}^*)\subset L_{n',n''}^*$). C'est un groupe compact et on note $K_{n',n''}^{u}$ son radical pro-$p$-unipotent. Si $\sharp=iso$ et $n''=0$, on pose simplement $K_{n,0}^+=K_{n,0}^{\pm}$ et on a l'isomorphisme $K_{n,0}^+/K_{n,0}^{u}={\bf SO}(2n+1;{\mathbb F}_{q})$. Hormis ce cas, $K_{n',n''}^{\pm}/K_{n',n''}^{u}$ s'identifie au sous-groupe des \'el\'ements $(g',g'')\in {\bf O}(2n'+1;{\mathbb F}_{q})\times {\bf O}(2n'')_{\sharp}({\mathbb F}_{q})$ tels que $det(g')det(g'')=1$. On identifie ce groupe \`a ${\bf SO}(2n'+1;{\mathbb F}_{q})\times {\bf O}(2n'')_{\sharp}({\mathbb F}_{q})$ par l'application $(g',g'')\mapsto (g'det(g''),g'')$. On note $K_{n',n''}^+$, resp. $K_{n',n''}^-$, l'image r\'eciproque dans $K_{n',n''}^{\pm}$ du sous-groupe ${\bf SO}(2n'+1;{\mathbb F}_{q})\times {\bf SO}(2n'')_{\sharp}({\mathbb F}_{q})$, resp. du sous-ensemble ${\bf SO}(2n'+1;{\mathbb F}_{q})\times {\bf O}^-(2n'')_{\sharp}({\mathbb F}_{q})$. 
 
 Un \'el\'ement de $G_{\sharp}(F)$ est dit compact si et seulement si le sous-groupe qu'il engendre est d'adh\'erence compacte. On sait qu'un \'el\'ement $g\in G_{\sharp}(F)$ est compact si et seulement si il existe $(n',n'')\in D_{\sharp}(n)$  et $h\in G_{\sharp}(F)$ de sorte que $h^{-1}gh\in K_{n',n''}^{\pm}$.
 
 \subsection{Repr\'esentations de r\'eduction unipotente}
 
  Pour $\sharp=iso$ ou $an$, notons $Irr_{t,\sharp}$ l'ensemble  des (classes d'isomorphismes de) repr\'esentations admissibles irr\'eductibles temp\'er\'ees de $G_{\sharp}(F)$. Notons $Irr_{t}$ la r\'eunion disjointe de $Irr_{t,iso}$ et $Irr_{t,an}$. 
 
 On sait conjecturalement  classifier l'ensemble $Irr_{t}$ de la fa\c{c}on suivante. On note $W_{F}$ le groupe de Weyl de $F$ et $Sp(2n,{\mathbb C})$ le groupe symplectique complexe d'un espace de dimension $2n$. Pour un homomorphisme $\psi:W_{F}\times SL(2,{\mathbb C})\to Sp(2n,{\mathbb C})$, on note $S(\psi)$ le groupe des composantes du centralisateur dans $Sp(2n,{\mathbb C})$ de l'image de $\psi$. C'est un produit fini de groupes ${\mathbb Z}/2{\mathbb Z}$ et on note $S(\psi)^{\wedge}$ son groupe de caract\`eres. On note $z(\psi)$ l'image naturelle dans $S(\psi)$ de l'\'el\'ement central $-1$ de $Sp(2n,{\mathbb C})$.  Alors $Irr_{t}$ est en bijection avec les classes de conjugaison par $Sp(2n,{\mathbb C})$ de couples $(\psi,\epsilon)$, o\`u  
 
  $\psi:W_{F}\times SL(2,{\mathbb C})\to Sp(2n,{\mathbb C})$ est un homomorphisme dont la restriction \`a $W_{F}$ est semi-simple et d'image born\'ee et dont la restriction \`a $SL(2,{\mathbb C})$ est alg\'ebrique;
  
   $\epsilon$ est un \'el\'ement de $S(\psi)^{\wedge}$. 
   
   On note $\pi(\psi,\epsilon)$ la repr\'esentation param\'etr\'ee par $(\psi,\epsilon)$. 
  Celle-ci appartient \`a $Irr_{t,, iso}$, resp. $Irr_{t,an}$, si et seulement si $\epsilon(z(\psi))=1$, resp. $\epsilon(z(\psi))=-1$. 
  
  Les repr\'esentations $\pi(\psi,\epsilon)$ sont caract\'eris\'ees par leur comportement par endoscopie et endoscopie tordue. Nous y reviendrons en 2.1. Dans le cas o\`u $\sharp=iso$, la classification et ces propri\'et\'es relatives \`a l'endoscopie ne sont plus conjecturales: elles ont \'et\'e \'etablies par Arthur,  cf. \cite{Ar} th\'eor\`eme 2.2.1. La situation pr\'esente du cas $\sharp=an$ est peu claire.

  Soit $\sharp=iso$ ou $an$ et soit $\pi$ une repr\'esentation admissible irr\'eductible de $G_{\sharp}(F)$ dans un espace complexe $E$. Pour $(n',n'')\in D_{\sharp}(n)$, le sous-espace d'invariants $E^{K_{n',n''}^{u}}$ est de dimension finie et est stable par l'action du groupe $K_{n',n''}^{\pm}$. Notons $\pi_{n',n''}$ la repr\'esentation de ce groupe dans cet espace d'invariants. Elle est triviale sur $K_{n',n''}^{u}$ et se descend en une repr\'esentation du groupe ${\bf SO}(2n'+1;{\mathbb F}_{q})\times {\bf O}(2n'')_{\sharp}({\mathbb F}_{q})$, que l'on note encore $\pi_{n',n''}$. On sait d\'efinir la notion de repr\'esentation unipotente de ce groupe fini. On dit que $\pi$ est de r\'eduction unipotente si et seulement s'il existe $ (n',n'')\in D_{\sharp}(n)$ tel que $E^{K_{n',n''}^{u}}$ soit non nul et que $\pi_{n',n''}$ soit unipotente. On sait qu'alors, pour tout autre $(n',n'')\in D_{\sharp}(n)$ tel  que $E^{K_{n',n''}^{u}}$ soit non nul, $\pi_{n',n''}$ est unipotente. On note $Irr_{unip,\sharp}$ l'ensemble des (classes d'isomorphisme de) repr\'esentations admissibles irr\'eductibles  de r\'eduction unipotente de $G_{\sharp}(F)$. On note $Irr_{tunip,\sharp}=Irr_{t,\sharp}\cap Irr_{unip,\sharp}$. On note $Irr_{unip}$, resp. $Irr_{tunip}$, la r\'eunion disjointe de $Irr_{unip,iso}$ et de $Irr_{unip,an}$, resp. de $Irr_{tunip,iso}$ et de $Irr_{tunip,an}$.
  
  Conjecturalement, l'ensemble $Irr_{tunip}$ est classifi\'e par les (classes de conjugaison des) couples $(\psi,\epsilon)$ comme ci-dessus v\'erifiant de plus la condition: la restriction de $\psi$ \`a $W_{F}$ est non ramifi\'ee. 
  Lusztig a effectivement classifi\'e l'ensemble $Irr_{tunip}$ par de tels couples $(\psi,\epsilon)$, cf. \cite{L} th\'eor\`eme 5.21.   
    
  On va d\'ecrire  de fa\c{c}on plus combinatoire l'ensemble des couples $(\psi,\epsilon)$ qui param\'etrisent l'ensemble $Irr_{tunip}$. Pour cela, introduisons quelques d\'efinitions. On appelle partition une classe d'\'equivalence de suites d\'ecroissantes  finies  de nombres entiers positifs ou nuls, deux suites \'etant \'equivalentes si elles ne diff\`erent que par des termes nuls. Pour une telle partition $\lambda=(\lambda_{1}\geq \lambda_{2}\geq...\geq \lambda_{r})$, on pose $S(\lambda)=\sum_{j=1,...,r}\lambda_{j}$ et on note $l(\lambda)$ le plus grand entier $j$ tel que $\lambda_{j}\not=0$. Cas particulier: on note $\emptyset$ la partition $(0,...)$ et on pose $l(\emptyset)=0$. On note $mult_{\lambda}$ la fonction sur ${\mathbb N}-\{0\}$ telle que, pour tout $i$ dans cet ensemble, $mult_{\lambda}(i)$ est le nombre d'entiers $j$ tels que $\lambda_{j}=i$. On note $Jord(\lambda)$ l'ensemble des $i\geq1$ tels que $mult_{\lambda}(i)\geq1$.
Pour $N\in {\mathbb N}$, on note ${\cal P}(N)$ l'ensemble des partitions $\lambda$ telles que $S(\lambda)=N$. On note ${\cal P}^{symp}(2N)$ l'ensemble des partitions symplectiques de $2N$, c'est-\`a-dire les $\lambda\in {\cal P}(2N)$ telles que $mult_{\lambda}(i)$ est pair pour tout entier $i$ impair. Pour une telle partition, on note $Jord_{bp}(\lambda)$ l'ensemble des entiers $i\geq2$ pairs tels que $mult_{\lambda}(i)\geq1$.  On note $\boldsymbol{{\cal P}^{symp}}(2N)$ l'ensemble des couples $(\lambda,\epsilon)$ o\`u $\lambda\in {\cal P}^{symp}(2N)$ et $\epsilon\in \{\pm 1\}^{Jord_{bp}(\lambda)}$.

  On sait que les classes de conjugaison d'homomorphismes alg\'ebriques $\rho:SL(2;{\mathbb C})\to Sp(2n;{\mathbb C})$ s'identifient aux orbites unipotentes dans $Sp(2n;{\mathbb C})$: \`a $\rho$ on associe l'orbite de l'image par $\rho$ d'un \'el\'ement unipotent  r\'egulier  de $SL(2;{\mathbb C})$. Ces orbites unipotentes sont elles-m\^emes classifi\'ees par ${\cal P}^{symp}(2n)$. Ainsi, \`a $\rho$, on associe une partition $\lambda\in {\cal P}^{symp}(2n)$. Inversement, pour toute telle partition $\lambda$, fixons un homomorphisme $\rho_{\lambda}$ classifi\'e  par $\lambda$. On note $Z(\lambda)$ le commutant de $\rho_{\lambda}$ dans $Sp(2n;{\mathbb C})$. 

Consid\'erons un couple $(\psi,\epsilon)$ param\'etrisant un \'el\'ement de $Irr_{tunip}$. A la restriction $\rho$ de $\psi$ \`a $SL(2;{\mathbb C})$ est associ\'ee une partition $\lambda\in {\cal P}^{symp}(2n)$. 
 Puisque la restriction de $\psi$ \`a $W_{F}$ est non ramifi\'ee, celle-ci est d\'etermin\'ee par l'image $s$ d'un \'el\'ement de Frobenius. C'est un \'el\'ement semi-simple de $Z(\lambda)$. Parce que la restriction de $\psi$ \`a $W_{F}$ est d'image born\'ee,  les valeurs propres de $s$ (consid\'er\'e comme un \'el\'ement de $GL(2n;{\mathbb C})$) sont de valeurs absolues $1$. On dit que $s$ est compact.  Notons $Z(\lambda,s)$ le commutant de $s$ dans $Z(\lambda)$, ${\bf Z}(\lambda,s)$ son groupe de composantes connexes et   ${\bf Z}(\lambda,s)^{\wedge}$ le groupe des caract\`eres de ${\bf Z}(\lambda,s)$. On a les \'egalit\'es $S(\psi)={\bf Z}(\lambda,s)$, $S(\psi)^{\wedge}={\bf Z}(\lambda,s)^{\wedge}$. Le terme $\epsilon$ appartient \`a  ${\bf Z}(\lambda,s)^{\wedge}$. La partition $\lambda$ \'etant fix\'ee, il y a une notion \'evidente de conjugaison par $Z(\lambda)$ du couple $(s,\epsilon)$. On voit que l'ensemble des classes de conjugaison de couples $(\psi,\epsilon)$ s'identifie aux classes de conjugaison au sens que l'on vient d'indiquer des triplets $(\lambda,s,\epsilon)$ v\'erifiant les conditions ci-dessus. On note $\mathfrak{Irr}_{tunip}$ cet ensemble de classes de conjugaison de triplets $(\lambda,s,\epsilon)$. Pour un tel triplet, on note $\pi(\lambda,s,\epsilon)$ la repr\'esentation associ\'ee par Lusztig au couple $(\psi,\epsilon)$ associ\'e au triplet. 
 
 Remarquons que, pour $\lambda\in {\cal P}^{symp}(2n)$ et $s\in Z(\lambda)$, l'\'el\'ement $s$ d\'efinit une d\'ecomposition de $\lambda$. En effet, consid\'erons les valeurs propres de $s$. Ce sont 
 $+1$ intervenant avec une multiplicit\'e paire $2n^+\geq0$, $-1$ intervenant avec une multiplicit\'e paire $2n^-\geq0$ et un ensemble de couples $(s_{j},s_{j}^{-1})$, $j$ parcourant un ensemble fini d'indices $J$, chaque $s_{j}$ \'etant un nombre complexe diff\'erent de $\pm 1$, et $s_{j}$ comme $s_{j}^{-1}$ intervenant avec une multiplicit\'e $m_{j}\geq1$.   Le commutant d'un tel $s$ dans $Sp(2n,{\mathbb C})$ est
$$Sp(2n^+;{\mathbb C})\times Sp(2n^-;{\mathbb C})\times \prod_{j\in J}GL(m_{j};{\mathbb C}).$$
L'homomorphisme $\rho_{\lambda}$  prend ses valeurs 
  dans ce commutant. Pour un groupe $H=GL(m_{j};{\mathbb C})$ ou $Sp(2n^+;{\mathbb C})$ ou $Sp(2n^-;{\mathbb C})$,  la classe de conjugaison d'un homomorphisme de $SL(2,{\mathbb C})$  \`a valeurs dans $H$ est d\'etermin\'ee comme ci-dessus par   une partition $\lambda_{j}\in {\cal P}(m_{j})$ si $H=GL(m_{j};{\mathbb C})$, resp. $\lambda^+\in {\cal P}^{symp}(2n^+)$ si $H=Sp(2n^+;{\mathbb C})$ et $\lambda^-\in {\cal P}^{symp}(2n^-)$ si $H=Sp(2n^-;{\mathbb C})$. Ainsi, \`a $\rho_{\lambda}$ sont associ\'ees des partitions  $\lambda^+$, $\lambda^-$ et  $\lambda_{j}$ pour $j\in J$. On a l'\'egalit\'e 
  $$\lambda=\lambda^+\cup\lambda^-\cup_{j\in J}(\lambda_{j}\cup \lambda_{j}),$$
  o\`u il s'agit de l'union usuelle des partitions. 
  On v\'erifie facilement que ${\bf Z}(\lambda,s)^{\wedge}$ s'identifie \`a $\{\pm 1\}^{Jord_{bp}(\lambda^+)}\times  \{\pm 1\}^{Jord_{bp}(\lambda^-)}$. Si $\epsilon\in {\bf Z}(\lambda,s)^{\wedge}$ s'identifie ainsi \`a un \'el\'ement $(\epsilon^+,\epsilon^-)\in \{\pm 1\}^{Jord_{bp}(\lambda^+)}\times  \{\pm 1\}^{Jord_{bp}(\lambda^-)}$, on v\'erifie l'\'egalit\'e
 $$(1) \qquad \epsilon(z(\psi))=(\prod_{i\in Jord_{bp}(\lambda^+)}\epsilon^+(i)^{mult_{\lambda^+}(i)})(\prod_{i\in Jord_{bp}(\lambda^-)}\epsilon^-(i)^{mult_{\lambda^-}(i)}),$$
 o\`u $\psi$ est l'homomorphisme associ\'e \`a $(\lambda,s)$. 
  Il r\'esulte des constructions de Lusztig que  $\pi(\lambda,s,\epsilon)$  est bien, comme on l'attend,  une repr\'esentation de $G_{iso}(F)$ si le produit ci-dessus vaut $1$ et de $G_{an}(F)$ s'il vaut $-1$. 

 Un cas particuli\`erement int\'eressant est celui o\`u $s$ n'a pour valeurs propres que $+1$ et $-1$, c'est-\`a-dire $s^2=1$. On $\mathfrak{Irr}_{unip-quad}$  le sous-ensemble des $(\lambda,s,\epsilon)\in \mathfrak{Irr}_{tunip}$ tels que $s^2=1$. La construction ci-dessus l'identifie \`a celui des  quadruplets $(\lambda^+,\epsilon^+,\lambda^-,\epsilon^-)$ tels que
 
 il existe $(n^+,n^-)\in D(n)$ de sorte $(\lambda^+,\epsilon^+)\in \boldsymbol{{\cal P}}^{symp}(2n^+)$ et $(\lambda^-,\epsilon^-)\in \boldsymbol{{\cal P}}^{symp}(2n^-)$.
 
 Si $(\lambda,s,\epsilon)$ correspond ainsi \`a $(\lambda^+,\epsilon^+,\lambda^-,\epsilon^-)$, on note 
  $\pi(\lambda^+,\epsilon^+,\lambda^-,\epsilon^-)=\pi(\lambda,s,\epsilon)$. On note $Irr_{unip-quad}$ l'ensemble de  ces repr\'esentations. 
    
  \bigskip
   
    \subsection{Repr\'esentations elliptiques}
    
       {\bf Notation.} Pour tout ensemble $X$, on note ${\mathbb C}[X]$ l'espace vectoriel sur ${\mathbb C}$ de base $X$.
   
   \bigskip
 
 Pour toute partition symplectique $\lambda$, on a d\'efini l'ensemble $Jord_{bp}(\lambda)$ des entiers pairs $i\geq2$ tels que $mult_{\lambda}(i)\geq1$. Pour tout entier $k\geq1$, notons plus pr\'ecis\'ement $Jord_{bp}^k(\lambda)$ l'ensemble des entiers pairs $i\geq2$ tels que $mult_{\lambda}(i)=k$.
 
 Notons $\mathfrak{Ell}_{unip}$ l'ensemble des quadruplets $(\lambda^+,\epsilon^+,\lambda^-,\epsilon^-)$ v\'erifiant les conditions suivantes:
 
 il existe des entiers $n^+,n^-\in {\mathbb N}$ tels que $n^++n^-=n$, $\lambda^+\in {\cal P}^{symp}(2n^+)$  et $\lambda^-\in {\cal P}^{symp}(2n^-)$;  
  
  pour $\zeta=\pm$ et $i\geq1$, $mult_{\lambda^{\zeta}}(i)=0$ si $i$ est impair et  $mult_{\lambda^{\zeta}}(i)\leq 2$ si $i$ est pair;
 
 pour $\zeta=\pm$, $\epsilon^{\zeta}$ est un \'el\'ement de $\{-1\}^{Jord_{bp}^1(\lambda^{\zeta})}$.
 
 A un tel quadruplet, on associe l'\'el\'ement $\pi_{ell}(\lambda^+,\epsilon^+,\lambda^-,\epsilon^-)\in {\mathbb C}[Irr_{tunip}]$ d\'efini par
 $$\pi_{ell}(\lambda^+,\epsilon^+,\lambda^-,\epsilon^-)=2^{-\vert Jord_{bp}^2(\lambda^+)\vert -\vert Jord_{bp}^2(\lambda^-)\vert }\sum_{\epsilon^{_{'}+},\epsilon^{_{'}-}}\left(\prod_{i\in Jord^2_{bp}(\lambda^{+})}\epsilon^{_{'}+}(i)\right)$$
 $$\left(\prod_{i\in Jord^2_{bp}(\lambda^{-})}\epsilon^{_{'}-}(i)\right)\pi(\lambda^+,\epsilon^{_{'}+},\lambda^-,\epsilon^{_{'}-}),$$
 o\`u la somme porte sur les $(\epsilon^{_{'}+},\epsilon^{_{'}-})\in \{\pm 1\}^{Jord_{bp}(\lambda^+)}\times  \{\pm 1\}^{Jord_{bp}(\lambda^-)}$ tels que, pour $\zeta=\pm$ et $i\in Jord^1_{bp}(\lambda^{\zeta})$, on ait $\epsilon^{_{'}\zeta}(i)=\epsilon^{\zeta}(i)$. On note $Ell_{unip}$ l'ensemble de ces repr\'esentations $\pi_{ell}(\lambda^+,\epsilon^+,\lambda^-,\epsilon^-)$ pour $(\lambda^+,\epsilon^+,\lambda^-,\epsilon^-)\in \mathfrak{Ell}_{unip}$.
 
 Une repr\'esentation  $\pi_{ell}(\lambda^+,\epsilon^+,\lambda^-,\epsilon^-)$ est  elliptique au sens d'Arthur, cf. \cite{Ar2} paragraphe 3. Plus pr\'ecis\'ement, $Ell_{unip}$  est exactement l'ensemble  des repr\'esentations elliptiques de r\'eduction unipotente des groupes $G_{iso}(F)$ et $G_{an}(F)$ (la d\'efinition des repr\'esentations elliptiques d\'epend de certains choix; on veut dire que l'on peut effectuer ceux-ci de sorte que l'assertion ci-dessus soit vraie).

\bigskip

\subsection{L'application de restriction, les espaces ${\cal R}^{par}$ et ${\cal R}^{par,glob}$}

 Pour tout $n'\in {\mathbb N}$, on note $C'_{n'}$ le sous-espace de l'espace des fonctions sur ${\bf SO}(2n'+1;{\mathbb F}_{q})$ (\`a valeurs  complexes) engendr\'e lin\'eairement par les traces de repr\'esentations unipotentes de ce groupe. Pour $\sharp=iso$ ou $an$ et pour tout entier $n''\geq1$, introduisons  l'espace des fonctions sur ${\bf O}(2n'')_{\sharp}({\mathbb F}_{q})$   engendr\'e par les traces de repr\'esentations unipotentes. Pour $\zeta=\pm$, on note $C_{n'',\sharp}^{''\zeta}$ l'espace des restrictions \`a ${\bf O}^{\zeta}(2n'')_{\sharp}({\mathbb F}_{q})$ des fonctions de l'espace pr\'ec\'edent (ou encore le sous-espace des \'el\'ements de cet espace pr\'ec\'edent qui sont nuls sur l'autre composante ${\bf O}^{-\zeta}(2n'')_{\sharp}({\mathbb F}_{q})$). On pose
 $$C''_{n''}=C_{n'',iso}^+\oplus C_{n'',iso}^-\oplus C_{n'',an}^+\oplus C_{n'',an}^-.$$
 Dans le cas o\`u $n''=0$, on pose formellement $C''_{0}=C_{0,iso}^{''+}={\mathbb C}$. On pose
 $${\cal R}^{par}=\oplus_{(n',n'')\in D(n)}C'_{n'}\otimes C''_{n''}.$$
 
 Les espaces $C'_{n'}$ etc... sont naturellement munis de produits hermitiens d\'efinis positifs. On en d\'eduit un tel produit sur ${\cal R}^{par}$.

 Soit $\sharp=iso$ ou $an$, soit $(n',n'')\in D_{\sharp}(n)$, soit $\zeta=\pm$ (avec $\zeta=+$ si $\sharp=iso$ et $n''=0$) et soit $\pi\in Irr_{unip,\sharp}$. On a d\'efini la repr\'esentation $\pi_{n',n''}$  de ${\bf SO}(2n'+1;{\mathbb F}_{q})\times {\bf O}(2n'')_{\sharp}({\mathbb F}_{q})$. On note $Res^{\zeta}_{n',n''}(\pi)$ la restriction    \`a la composante ${\bf SO}(2n'+1;{\mathbb F}_{q})\times {\bf O}^{\zeta}(2n'')_{\sharp}({\mathbb F}_{q})$
   de la trace de $\pi_{n',n''}$. Cette fonction s'identifie \`a un \'el\'ement de $C'_{n'}\otimes C_{n'',\sharp}^{_{''}\zeta}$, donc \`a un \'el\'ement de ${\cal R}^{par}$. On note $Res(\pi)$  la somme des $Res_{n',n''}^{\zeta}(\pi)$  sur les triplets $(n',n'',\zeta)$ soumis aux restrictions indiqu\'ees ci-dessus.   Cela d\'efinit une application $Res:Irr_{unip} \to {\cal R}^{par}$.
 Cette application   se  prolonge en une application lin\'eaire $Res:{\mathbb C}[Irr_{unip}]\to {\cal R}^{par}$. 

Consid\'erons un entier $m\in \{1,...,n\}$, posons $n_{0}=n-m$. Notons $C^{GL(m)}$ l'espace de fonctions sur le groupe fini ${\bf GL}(m;{\mathbb F}_{q})$ engendr\'e par les traces de repr\'esentations unipotentes.  
 On d\'efinit deux applications lin\'eaires
 $$res_{m}',res_{m}'':{\cal R}^{par}\to C^{GL(m)}\otimes {\cal R}^{par}_{n_{0}}$$
 de la fa\c{c}on suivante. Consid\'erons une composante $C'_{n'}\otimes C_{n'',\sharp}^{_{''}\zeta}$ de ${\cal R}^{par}$. Si $n'<m$, $res_{m}'$ est nulle sur cette composante. Si $n'\geq m$, on introduit un sous-groupe parabolique ${\bf P}$ de ${\bf SO}(2n'+1)$ dont une composante de Levi est isomorphe \`a ${\bf GL}(m)\times {\bf SO}(2n'-2m+1)$. L'application module de Jacquet envoie $C'_{n'}$ dans $C^{GL(m)}\otimes C'_{n'-m}$ et elle ne d\'epend pas du choix de ${\bf P}$. Alors $res_{m}'$ est \'egale sur  $C'_{n'}\otimes C_{n'',\sharp}^{_{''}\zeta}$ \`a l'application
 $$C'_{n'}\otimes C_{n'',\sharp}^{_{''}\zeta}\to C^{GL(m)}\otimes C'_{n'-m}\otimes C_{n'',\sharp}^{_{''}\zeta}$$
 produit tensoriel de l'application pr\'ec\'edente et de l'identit\'e de $C_{n'',\sharp}^{_{''}\zeta}$. Si $\sharp=iso$ et $n''<m$ ou si $\sharp=an$ et $n''\leq m$, $res_{m}''$ est nulle sur notre composante. Sinon, on introduit comme ci-dessus un sous-groupe parabolique ${\bf P}$ de ${\bf O}(2n'')_{\sharp}$ dont une composante de Levi est isomorphe \`a ${\bf GL}(m)\times {\bf O}(2n''-2m)_{\sharp}$. On a de nouveau une application module de Jacquet  $C_{n'',\sharp}^{_{''}\zeta}\to C^{GL(m)}\otimes C_{n''-m,\sharp}^{_{''}\zeta}$, dont on d\'eduit l'application
  $$res_{m}'':C'_{n'}\otimes C_{n'',\sharp}^{_{''}\zeta}\to C^{GL(m)}\otimes C'_{n'}\otimes C_{n''-m,\sharp}^{_{''}\zeta}.$$ 
  
  On note ${\cal R}^{par,glob}$  le sous-espace des \'el\'ements $\phi\in {\cal R}^{par}$ tels que $res'_{m}(\phi)=res''_{m}(\phi)$ pour tout $m\in \{1,...,n\}$. On note $res_{m}$ la restriction \`a ce sous-espace de l'une ou l'autre des applications $res'_{m}$ ou $res''_{m}$.  C'est une application lin\'eaire
  $$res_{m}:{\cal R}^{par,glob}\to C^{GL(m)}\otimes {\cal R}^{par}_{n_{0}}.$$
  Elle prend ses valeurs dans $C^{GL(m)}\otimes {\cal R}^{par,glob}_{n_{0}}$. En effet, pour $r\in \{1,...,n_{0}\}$, notons $res'_{m,r}$ la compos\'ee de $res_{m}$ et de l'application $res'_{r}$ d\'efinie sur ${\cal R}^{par}_{n_{0}}$. On d\'efinit de fa\c{c}on similaire l'application $res''_{m,r}$. On voit que ces deux applications sont compos\'ees de $res'_{m+r}$, resp. $res''_{m+r}$, et d'une application module de Jacquet de $C^{GL(m+r)}$ dans $C^{GL(m)}\otimes C^{GL(r)}$. L'\'egalit\'e $res'_{m+r}=res''_{m+r}$ sur ${\cal R}^{par,glob}$ entra\^{\i}ne l'\'egalit\'e $res'_{m,r}=res''_{m,r}$, d'o\`u notre assertion.

On v\'erifie que l'application $Res$ introduite  plus haut prend ses valeurs dans ${\cal R}^{par,glob}$.  Plus pr\'ecis\'ement, soient $m\in \{1,...,n\}$, $\sharp=iso$ ou $an$ et $\pi\in Irr_{unip,\sharp}$. Introduisons un sous-groupe parabolique $P$ de $G_{\sharp}$ dont une composante de Levi soit isomorphe \`a $GL(n)\times G_{n_{0},\sharp}$ (en supposant $m<n$ si $\sharp=an$). Notons $\pi_{M}$ le module de Jacquet de $\pi$ relatif \`a $P$. On a une application  $Res^M:{\mathbb C}[Irr_{unip,\sharp}]\to C^{GL(m)}\otimes {\cal R}^{par}_{n_{0}}$ similaire \`a $Res$. On a l'\'egalit\'e

(1) $res_{m}\circ Res(\pi)=Res^M(\pi_{M})$. 

Cela r\'esulte directement de \cite{MP} proposition 6.7. Dans le cas o\`u $\sharp=an$ et $m=n$, on a $res_{m}\circ Res(\pi)=0$.

Dans les espaces $C'_{n'}$ et $C_{n'',\sharp}^{_{''}\zeta}$, on sait d\'efinir le sous-espace des fonctions cuspidales et la projection orthogonale sur ce sous-espace, que l'on note $proj_{cusp}$. 
 
 {\bf Remarque.} Dans le cas particulier de l'espace $C_{1,iso}^{_{''}+}$, correspondant au groupe ${\bf SO}(2)_{iso}\simeq {\bf GL}(1)$, cette projection est nulle.
 \bigskip

On note ${\cal R}_{cusp}^{par}$   le sous-espace des \'el\'ements de ${\cal R}$ dont toutes les composantes sont cuspidales.   Des projections pr\'ec\'edentes se d\'eduit une projection $proj_{cusp}:{\cal R}^{par}\to {\cal R}_{cusp}^{par}$.  L'espace ${\cal R}_{cusp}$  est inclus dans ${\cal R}^{par,glob}$. Pour tout entier $m\geq1$, on note aussi $C_{cusp}^{GL(m)}$ le sous-espace des \'el\'ements cuspidaux de $C^{GL(m)}$. On sait qu'il est de dimension $1$.  Notons ${\cal P}(\leq n)$ l'ensemble des partitions ${\bf m}=(m_{1},...,m_{t})$ telles que $S({\bf m})\leq n$. Soit  ${\bf m}=(m_{1},...,m_{t}) \in {\cal P}(\leq n)$, posons $n_{0}=n-S({\bf m})$.  En it\'erant la construction
pr\'ec\'edente, on obtient une application lin\'eaire
$$res_{{\bf m}}:{\cal R}^{par,glob}\to   {\cal R}^{par,glob}_{{\bf m}}:=C^{GL(m_{1})}\otimes...\otimes C^{GL(m_{t})}\otimes {\cal R}^{par,glob}_{n_{0}}.$$
On pose
$${\cal R}_{{\bf m},cusp}^{par}=C_{cusp}^{GL(m_{1})}\otimes...\otimes C_{cusp}^{GL(m_{t})}\otimes {\cal R}^{par}_{n_{0},cusp}\simeq {\cal R}^{par}_{n_{0},cusp}.$$
 On dispose aussi des projections cuspidales de ${\cal R}_{n_{0}}^{par}$ sur ${\cal R}_{n_{0},cusp}^{par}$ et de $C^{GL(m_{j})}$ sur $C_{cusp}^{GL(m_{j})}$. On note $proj_{cusp}$ ces projections ainsi que les produits tensoriels de telles projections. On a donc une application lin\'eaire
$$(2) \qquad {\cal R}^{par,glob}\to \oplus_{{\bf m}} {\cal R}^{par}_{{\bf m},cusp},$$
qui est la somme  des applications $proj_{cusp}\circ res_{{\bf m}}$ sur toutes les partitions ${\bf m}\in {\cal P}(\leq n)$. On v\'erifie facilement que c'est un isomorphisme. 

 La restriction de l'application $proj_{cusp}\circ Res$ au sous-espace ${\mathbb C}[Ell_{unip}]$ de ${\mathbb C}[Irr_{tunip}]$ est un isomorphisme de ce sous-espace sur l'espace ${\cal R}^{par}_{cusp}$, cf. \cite{MW} 4.2 et 5.4. 

{\bf Remarque.} On peut d\'efinir l'application $res_{{\bf m}}$  pour toute suite d'entiers positifs ${\bf m}=(m_{1},...,m_{t})$ telle que  $m_{1}+...+m_{t}\leq n$ (on n'a pas besoin que ${\bf m}$ soit une partition, c'est-\`a-dire que $m_{1}\geq...\geq m_{t}$). Mais, pour une telle suite, $res_{{\bf m}}$ se d\'eduit par permutation des facteurs de $res_{{\bf m}'}$, o\`u ${\bf m}'$ est la partition associ\'ee \`a ${\bf m}$, c'est-\`a-dire celle qui a les m\^emes termes que ${\bf m}$, ordonn\'es de fa\c{c}on d\'ecroissante. 
 
\bigskip

\subsection{Egalit\'e de restrictions aux \'el\'ements compacts} 

On va consid\'erer deux situations auxquelles s'appliqueront le lemme ci-dessous.

Dans le cas (A), on consid\`ere un triplet $(\lambda,s,\epsilon)\in \mathfrak{Irr}_{tunip}$. On pose $\pi=\pi(\lambda,s,\epsilon)$.  On sait d\'ecrire le groupe $Z(\lambda)$. C'est le produit des $Sp(mult_{\lambda}(i);{\mathbb C})$ pour $i\in Jord(\lambda)$ impair et des $O(mult_{\lambda}(i);{\mathbb C})$ pour $i\in Jord_{bp}(\lambda)$. Ainsi, l'\'el\'ement $s\in Z(\lambda)$ se d\'ecompose en produit de $s_{i}$ pour $i\in Jord(\lambda)$. Fixons $i_{0}\in Jord(\lambda)$ et supposons que $s_{i_{0}}$ ait une valeur propre $z\not=\pm 1$ (il a donc aussi la valeur propre $z^{-1}$). On peut alors fixer une d\'ecomposition
$${\mathbb C}^{2n}={\mathbb C}^{2i_{0}}\oplus {\mathbb C}^{2n-2i_{0}}$$
qui soit stable par $s$ et par l'homomorphisme $\rho_{\lambda}$, de sorte que les composantes $\underline{s}$ et $\underline{\rho}$ de ces termes \`a valeurs dans la premi\`ere composante ${\mathbb C}^{2i_{0}}$ v\'erifient:

$\underline{s}$ a pour valeurs propres $z$ et $z^{-1}$, chacune avec multiplicit\'e $i_{0}$;

$\underline{\rho}$ est param\'etr\'e par la partition $(i_{0},i_{0})$.

On note $s_{0}$ et $\rho_{0}$ les composantes de $s$ et $\rho_{\lambda}$ dans l'autre composante ${\mathbb C}^{2n-2i_{0}}$ et on note $\lambda_{0}$ la partition associ\'ee \`a $\rho_{0}$. On note $\bar{s}$ l'\'el\'ement qui agit comme $s_{0}$ dans cette deuxi\`eme composante mais qui agit par l'identit\'e sur la premi\`ere ${\mathbb C}^{2i_{0}}$. Si $i_{0}$ est impair ou si $i_{0}$ et pair et $s_{i_{0}}$ admet la valeur propre $1$, on voit que le groupe ${\bf Z}(\lambda,\bar{s})$ est isomorphe \`a ${\bf Z}(\lambda,s)$. L'\'el\'ement $\epsilon$  s'identifie \`a un \'el\'ement de ${\bf Z}(\lambda,\bar{s})^{\wedge}$, le triplet $(\lambda,\bar{s},\epsilon)$ appartient \`a $\mathfrak{Irr}_{tunip}$ et on pose $\bar{\pi}=\pi(\lambda,\bar{s},\epsilon)$. Si $i_{0}$ est pair et $s_{i_{0}}$ n'a pas $1$ pour valeur propre, le groupe ${\bf Z}(\lambda,\bar{s})$ est plus gros que ${\bf Z}(\lambda,s)$. Plus pr\'ecis\'ement, le deuxi\`eme s'identifie \`a un sous-groupe d'indice $2$ du premier. Il y a deux prolongements du caract\`ere $\epsilon$ au groupe ${\bf Z}(\lambda,\bar{s})$, que l'on note $\bar{\epsilon}'$ et $\bar{\epsilon}''$. On pose $\bar{\pi}=\pi(\lambda,\bar{s},\bar{\epsilon}')+\pi(\lambda,\bar{s},\bar{\epsilon}'')$. 
 C'est un \'el\'ement de ${\mathbb C}[Irr_{tunip}]$. 

Dans le cas (B), on consid\`ere un entier $i_{0}\geq1$, deux partitions symplectiques $\bar{\lambda}^+$ et $\lambda^-$ et des \'el\'ements $\bar{\epsilon}^+\in \{\pm 1\}^{Jord_{bp}(\bar{\lambda}^+)}$ et $\epsilon^-\in \{\pm 1\}^{Jord_{bp}(\lambda^-)}$. On suppose $2n=2i_{0}+S(\bar{\lambda}^+)+S(\lambda^-)$. On pose $\lambda^+=\bar{\lambda}^+\cup\{i_{0},i_{0}\}$. 
  Si $i_{0}$ est impair ou si $i_{0}$ est pair et $mult_{\bar{\lambda}^+}(i_{0})\geq1$, on a $Jord_{bp}(\lambda^+)=Jord_{bp}(\bar{\lambda}^+)$ et $\bar{\epsilon}^+$ est aussi un \'el\'ement de $Jord_{bp}(\lambda^+)$. On pose $\pi=\pi(\lambda^+,\bar{\epsilon}^+,\lambda^-,\epsilon^-)$. Si $i_{0}$ est pair et $mult_{\bar{\lambda}^+}(i_{0})=0$, on a $Jord_{bp}(\lambda^+)=Jord_{bp}(\bar{\lambda}^+)\cup\{i_{0}\}$.  Il y a deux prolongements de $\bar{\epsilon}^+$ en des \'el\'ements de $Jord_{bp}(\lambda^+)$, que l'on note $\epsilon^{'+}$ et $\epsilon^{''+}$. On pose   $\pi=\pi(\lambda^+,\epsilon^{'+},\lambda^-,\epsilon^-)+\pi(\lambda^+,\epsilon^{''+},\lambda^-,\epsilon^-)$.   On pose $\bar{\lambda}^-=\lambda^-\cup\{i_{0},i_{0}\}$.  Si $i_{0}$ est impair ou si $i_{0}$ est pair et $mult_{\lambda^-}(i_{0})\geq1$, on a $Jord_{bp}(\bar{\lambda}^-)=Jord_{bp}(\lambda^-)$ et $\epsilon^-$ appara\^{\i}t comme un \'el\'ement de $\{\pm 1\}^{Jord_{bp}(\bar{\lambda}^-)}$. On pose $\bar{\pi}=\pi(\bar{\lambda}^+,\bar{\epsilon}^+,\bar{\lambda}^-,\epsilon^-)$. Si $i_{0}$ est pair et $mult_{\lambda^-}(i_{0})=0$, on a $Jord_{bp}(\bar{\lambda}^-)=Jord_{bp}(\lambda^-)\cup\{i_{0}\}$ et il y a de nouveau deux   prolongements de $\epsilon^-$ en un \'el\'ement de $\{\pm 1\}^{Jord_{bp}(\bar{\lambda}^-)}$, que l'on note $\bar{\epsilon}^{'-}$ et $\bar{\epsilon}^{''-}$. On pose $\bar{\pi}=\pi(\bar{\lambda}^+,\bar{\epsilon}^+,\bar{\lambda}^-,\bar{\epsilon}^{'-})+\pi(\bar{\lambda}^+,\bar{\epsilon}^+,\bar{\lambda}^-,\bar{\epsilon}^{''-})$.

\ass{Lemme}{Dans les deux  cas (A) et (B) ci-dessus, on a l'\'egalit\'e $Res(\pi)=Res(\bar{\pi})$.}

{\bf Remarque.} On peut consid\'erer des situations similaires \`a celles ci-dessus mais o\`u l'on \'echange les r\^oles des  valeurs propres $1$ et $-1$ dans le  cas (A), ou des exposants $+$ et $-$ dans le cas (B). Evidemment, le lemme vaut aussi dans ces cas.

\bigskip
 
Preuve. On consid\`ere d'abord la premi\`ere situation. 
  Notons $\sharp$ l'indice tel que nos repr\'esentations soient des repr\'esentations de $G_{\sharp}(F)$.  Introduisons un sous-groupe parabolique $P$ de $G_{\sharp}$ de composante de Levi
 $$M=GL( i_{0})\times G_{n-i_{0},\sharp}.$$
 On a introduit les termes $s_{0}$ et $\lambda_{0}$. Le groupe $Z(\lambda_{0},s_{0})$ est naturellement un sous-groupe de $Z(\lambda,s)$, d'o\`u un homomorphisme de ${\bf Z}(\lambda_{0},s_{0})$ dans ${\bf Z}(\lambda,s)$. On v\'erifie que c'est un isomorphisme (parce que $z\not=\pm 1$). Ainsi, $\epsilon$ s'identifie \`a un \'el\'ement de ${\bf Z}(\lambda_{0},s_{0})$ et on d\'efinit la repr\'esentation $\pi_{0}=\pi(\lambda_{0},s_{0},\epsilon)$ de $G_{n-i_{0},\sharp}(F)$. Notons $st_{i_{0}}$ la repr\'esentation de Steinberg de $GL(i_{0};F)$. 
 Pour un \'el\'ement $y\in {\mathbb C}$, notons $st_{i_{0}}(\vert .\vert_{F}^y\circ det)\times \pi_{0}$ la repr\'esentation de $G_{\sharp}(F)$ induite \`a l'aide de $P$ de la repr\'esentation $st_{i_{0}}(\vert .\vert_{F} ^y\circ det)\otimes \pi _{0}$ de $M(F)$. Soit $y $ un nombre complexe tel que $q^{-y}=z$. Il est connu que $\pi=st_{i_{0}}(\vert .\vert_{F} ^{y}\circ det)\times \pi _{0}$. La repr\'esentation $st_{i_{0}}(\vert .\vert_{F} ^0\circ det)\times \pi_{0}$ est irr\'eductible si $i_{0}$ est impair ou si $i_{0}$ est pair et $s_{i_{0}}$ a pour valeur propre $1$; elle est r\'eductible si $i_{0}$ est pair et $ s_{i_{0}}$ n'a pas $1$ pour valeur propre. Dans les deux cas, il est connu que son image dans le groupe de Grothendieck des repr\'esentations lisses de longueur finie de $G_{\sharp}(F)$ est \'egale \`a $\bar{\pi}$.  Il suffit donc de prouver que, pour $y\in {\mathbb C}$, $Res(st_{i_{0}}(\vert .\vert_{F} ^y\circ det)\times \pi _{0})$ ne d\'epend pas de $y$. Ce terme se d\'eduit des restrictions de la repr\'esentation $st_{i_{0}}(\vert .\vert _{F}^y\circ det)\times \pi _{0}$ aux diff\'erents groupes $K_{n',n''}^{\pm}$. 
 Ces restrictions s'identifient \`a des repr\'esentations de dimension finie de groupes finis, que l'on peut \'ecrire $\sum_{\rho}m(y,\rho)\rho$, o\`u $\rho$ parcourt les repr\'esentations irr\'eductibles du groupe en question et les multiplicit\'es $m(y,\rho)$ sont des entiers positifs ou nuls. Il est clair que les caract\`eres de ces repr\'esentations varient contin\^ument en $y$, ce qui entra\^{\i}ne que les multiplicit\'es  $m(y,\rho)$ aussi. Comme ce sont des entiers, il sont constants en $z$. Cela d\'emontre le lemme dans  le cas (A) .

Consid\'erons maintenant  le cas (B). On note $(\lambda_{0},s_{0},\epsilon)$ l'\'el\'ement de $\mathfrak{Irr}_{n-i_{0},tunip}$ auquel s'identifie le quadruplet $(\bar{\lambda}^+,\bar{\epsilon}^+,\lambda^-,\epsilon^-)$. On fixe un nombre complexe $z\not=\pm 1$ de valeur absolue $1$. On introduit la partition $\underline{\lambda}=(i_{0},i_{0})$ et un \'el\'ement $\underline{s}\in Z(\underline{\lambda})$ ayant deux valeurs propres $z$ et $z^{-1}$. On note $(\lambda,s)$ la somme directe, en un sens \'evident, de $(\underline{\lambda},\underline{s})$ et $(\lambda_{0},s_{0})$. On a l'\'egalit\'e ${\bf Z}(\lambda,s)={\bf Z}(\lambda_{0},s_{0})$ et $\epsilon$ peut \^etre consid\'er\'e comme un \'el\'ement de ${\bf Z}(\lambda,s)^{\wedge}$. On pose alors $\Pi=\pi(\lambda,s,\epsilon)$. 
 Appliquons le cas (A) aux donn\'ees $\lambda$, $s$, $\epsilon$, \`a l'entier $i_{0}$ et \`a la valeur propre $z$ de $s_{i_{0}}$. On en d\'eduit un \'el\'ement   de ${\mathbb C}[Irr_{tunip}]$  que l'on note $\bar{\Pi}^+$ (pour le distinguer du $\bar{\pi}$ dont on dispose d\'ej\`a). On vient de d\'emontrer que $Res(\Pi)=Res(\bar{\Pi}^+)$. Mais on voit que, par construction, $\bar{\Pi}^+$ n'est autre que le pr\'esent $\pi$. Donc $Res(\Pi)=Res(\pi)$. Maintenant, on applique le cas (A) aux m\^emes donn\'ees, mais en \'echangeant les r\^oles des  valeurs propres $1$ et $-1$, ce qui est loisible ainsi qu'on l'a remarqu\'e. On obtient un autre \'el\'ement $\bar{\Pi}^-$ et l'\'egalit\'e $Res(\Pi)=Res(\bar{\Pi}^-)$. De nouveau, par construction, on a l'\'egalit\'e $\bar{\Pi}^-=\bar{\pi}$. Donc $Res(\Pi)=Res(\bar{\pi})$, puis l'\'egalit\'e $Res(\pi)=Res(\bar{\pi})$, qui ach\`eve la d\'emonstration. $\square$  

\bigskip

 \subsection{L'involution de Aubert-Zelevinsky}
 
 Soit $\sharp=iso$ ou $an$. On sait d\'efinir une involution de l'ensemble des classes de repr\'esentations admissibles irr\'eductibles, qui g\'en\'eralise l'involution introduite par Zelevinsky dans le cas du groupe $GL(n)$, cf. \cite{Au} ou \cite{SS} paragraphe III.3. On la note $D$. Il est connu qu'elle conserve l'ensemble $Irr_{unip,\sharp}$.  Pour un groupe ${\bf SO}(2n'+1)$ ou ${\bf O}(2n'')_{\sharp}$ d\'efini sur ${\mathbb F}_{q}$, on sait d\'efinir une involution similaire. Par produit tensoriel et sommation, on en d\'eduit une involution $D^{par}$ de ${\cal R}^{par}$. On sait que les involutions en question commutent en un sens convenable \`a l'application module de Jacquet. Il en r\'esulte que $D^{par}$ conserve l'espace ${\cal R}^{par,glob}$.
 
 \ass{Lemme}{On a l'\'egalit\'e $D^{par}\circ Res=Res\circ D$.}
 
 Preuve. Fixons $(n',n'')\in D_{\sharp}(n)$  et une repr\'esentation $\pi\in Irr_{unip,\sharp}$. On a d\'efini  en 1.3 la repr\'esentation $\pi_{n',n''}$ de ${\bf G}_{\sharp}({\mathbb F}_{q})$, o\`u ${\bf G}_{\sharp}={\bf SO}(2n'+1) \times {\bf O}(2n'')_{\sharp} $. Notons  ${\bf D}$ l'involution de l'ensemble des repr\'esentations de ${\bf G}_{\sharp}({\mathbb F}_{q})$. On doit prouver que ${\bf D}(\pi_{n',n''})=(D(\pi))_{n',n''}$.

 Notons $P_{min}$ le sous-groupe parabolique minimal de $G_{\sharp}$ form\'e des \'el\'ements qui, avec les notations de 1.1, stabilisent les drapeaux de sous-espaces 
 
 $Fv_{1}$, $Fv_{1}\oplus Fv_{2}$,...,$Fv_{1}\oplus...\oplus Fv_{n}$, si $\sharp=iso$,
 
 $Fv_{2}$, $Fv_{2}\oplus Fv_{3}$,...,$Fv_{2}\oplus...\oplus Fv_{n}$, si $\sharp=an$. 
 
 On note $M_{min}$ sa composante de Levi ''\'evidente''.

 Les sous-groupes paraboliques standard de $G_{\sharp}$, c'est-\`a-dire contenant $P_{min}$, sont en bijection avec les multiplets d'entiers ${\bf m}=(m_{1},...,m_{t},n_{0})$, o\`u $t\in {\mathbb N}$, $m_{j}\geq1$ pour tout $j\in \{1,...,t\}$, $n_{0}+\sum_{j=1,...,t}m_{j}=n$ et $n_{0}\geq0$ si $\sharp=iso$, $n_{0}\geq1$ si $\sharp=an$. On note $P_{{\bf m}}$ le sous-groupe parabolique standard associ\'e \`a ${\bf m}$ et $M_{{\bf m}}$ sa composante de Levi standard, c'est-\`a-dire contenant $M_{min}$.  On a un isomorphisme 
$$M_{{\bf m}}\simeq GL(m_{1})\times...\times GL(m_{t})\times G_{n_{0},\sharp} .$$
Par d\'efinition, on a l'\'egalit\'e suivante, que l'on expliquera ci-dessous:
$$(1) \qquad (-1)^{n(\pi)}D(\pi)=\sum_{{\bf m}}(-1)^{t({\bf m})}Ind_{P_{{\bf m}}}^{G_{\sharp}}\circ res_{P_{{\bf m}}}^{G_{\sharp}}(\pi).$$
Par ''\'egalit\'e'', on veut dire ici que les deux membres ont m\^eme image dans le groupe de Grothendieck des repr\'esentations admissibles de longueur finie de $G_{\sharp}(F)$. On a not\'e $Ind_{P_{{\bf m}}}^{G_{\sharp}}$ le foncteur d'induction et $res_{P_{{\bf m}}}^{G_{\sharp}}$ le foncteur module de Jacquet. On a not\'e $t({\bf m})$ l'entier $t$ qui figure dans la donn\'ee ${\bf m}$.  Enfin, $n(\pi)$ est un entier d\'ependant de $\pi$ sur lequel on reviendra plus loin. Fixons ${\bf =}(m_{1},...,m_{t},n_{0})$, notons 
$$\pi_{{\bf m}}=Ind_{P_{{\bf m}}}^{G_{\sharp}}\circ res_{P_{{\bf m}}}^{G_{\sharp}}(\pi)$$
et calculons $(\pi_{{\bf m}})_{n',n''}$. Cette repr\'esentation est une somme index\'ee par les doubles classes $P_{{\bf m}}(F)\backslash G_{\sharp}(F)/ K_{n',n''}^{\pm}$. On voit que  cet ensemble de doubles classes est index\'e par l'ensemble $X_{{\bf m}}$ des couples de multiplets d'entiers ${\bf m}'=(m'_{1},...,m'_{t},n'_{0})$, ${\bf m}''=(m''_{1},...,m''_{t},n''_{0})$ v\'erifiant les conditions suivantes:

 pour tout $j=1,...,t$, $m'_{j},m''_{j}\geq0$ et $m_{j}=m'_{j}+m''_{j}$;

 $n_{0}=n'_{0}+n''_{0}$, $n'_{0}\geq0$, $n''_{0}\geq0$ si $\sharp=iso$ et $n''_{0}\geq1$ si $\sharp=an$;

$n'=n'_{0}+\sum_{j=1,...,t}m'_{j}$, $n''=n''_{0}+\sum_{j=1,...,t}m''_{j}$. 

Ainsi, on a une \'egalit\'e $(\pi_{{\bf m}})_{n',n''}=\sum_{({\bf m}',{\bf m}'')\in X_{{\bf m}}}\pi_{{\bf m}',{\bf m}''}$.

Fixons un  couple $({\bf m}',{\bf m}'')\in X_{{\bf m}}$. Il d\'etermine plusieurs sous-groupes. D'une part un sous-groupe parabolique ${\bf P}'_{{\bf m}'}$ de ${\bf SO}(2n'+1)$, dont une composante de Levi ${\bf M}'_{{\bf m}'}$ est isomorphe \`a
$${\bf GL}(m'_{1})\times ...\times {\bf GL}(m'_{t})\times {\bf SO}(2n'_{0}+1),$$
et  un sous-groupe parabolique ${\bf P}''_{{\bf m}''}$ de ${\bf O}(2n'')_{\sharp}$, dont une composante de Levi ${\bf M}''_{{\bf m}''}$ est isomorphe \`a
$${\bf GL}(m''_{1})\times ...\times {\bf GL}(m''_{t})\times {\bf O}(2n''_{0})_{\sharp}.$$
On pose ${\bf P}_{{\bf m}',{\bf m}''}={\bf P}'_{{\bf m}'}\times {\bf P}''_{{\bf m}''}$ et ${\bf M}_{{\bf m}',{\bf m}''}={\bf M}'_{{\bf m}'}\times {\bf M}''_{{\bf m}''}$
Le couple $({\bf m}',{\bf m}'')$ d\'etermine d'autre part un sous-groupe compact $K_{{\bf m}',{\bf m}''}$ de $M_{{\bf m}}(F)$, qui est produit de sous-groupes des diverses composantes. Le sous-groupe de $G_{n_{0},\sharp}(F)$ est  (\`a conjugaison pr\`es) $K_{n'_{0},n''_{0}}^{\pm}$. Pour $j=1,...,t$, le sous-groupe de $GL(m_{j};F)$ est un sous-groupe parahorique  dont le groupe r\'esiduel associ\'e est isomorphe \`a
$${\bf GL}(m'_{j})\times {\bf GL}(m''_{j}).$$
La repr\'esentation $\pi_{{\bf m}',{\bf m}''}$ de ${\bf SO}(2n'+1)({\mathbb F}_{q})\times {\bf O}(2n'')_{\sharp}({\mathbb F}_{q})$ associ\'ee \`a $({\bf m}',{\bf m}'')$ s'obtient de la fa\c{c}on suivante. Par un proc\'ed\'e analogue \`a celui  de 1.3, le groupe $K_{n',n''}^{\pm}$ de ce paragraphe \'etant remplac\'e par $K_{{\bf m}',{\bf m}''}^{\pm}$, on d\'eduit de la repr\'esentation $res_{P_{{\bf m}}}^{G_{\sharp}}(\pi)$ de $M_{{\bf m}}(F)$ une repr\'esentation du groupe r\'esiduel de $K_{{\bf m}',{\bf m}''}^{\pm}$. Par permutation des facteurs, cette repr\'esentation devient une repr\'esentation $\sigma_{{\bf m}',{\bf m}''}$ de ${\bf M}_{{\bf m}',{\bf m}''}({\mathbb F}_{q}) $. Alors $\pi_{{\bf m}',{\bf m}''}=Ind_{{\bf P}_{{\bf m}',{\bf m}''}}^{{\bf G}_{\sharp}}(\sigma_{{\bf m}',{\bf m}''})$.
D'apr\`es  \cite{MP} proposition 6.7, $\sigma_{{\bf m}',{\bf m}''}$ n'est autre que  l'image de $\pi_{n',n''}$ par le foncteur module de Jacquet $res_{{\bf M}_{{\bf m}',{\bf m}''}}^{{\bf G}_{\sharp}}$. On obtient
$$(2) \qquad (\pi_{{\bf m}})_{n',n''}=\sum_{({\bf m}',{\bf m}'')\in X_{{\bf m}}}Ind_{{\bf P}_{{\bf m}',{\bf m}''}}^{{\bf G}_{\sharp}}\circ res_{{\bf M}_{{\bf m'},{\bf m}''}}^{{\bf G}_{\sharp}}(\pi_{n',n''}).$$

Notons $X$ l'ensemble des couples de multiplets d'entiers ${\bf r}'=(r'_{1},...,r'_{t'},n'_{0})$, ${\bf r}''=(r''_{1},...,r''_{t''},n''_{0})$ tels que

 $r'_{j}\geq1$ pour tout $j=1,...,t'$ et $r''_{j}\geq1$ pour tout $j=1,...,t''$;

 $n'_{0}\geq0$, $n''_{0}\geq0$ si $\sharp=iso$ et $n''_{0}\geq1$ si $\sharp=an$;

 $n'=n'_{0}+\sum_{j=1,...,t'}r'_{j}$,   $n''=n''_{0}+\sum_{j=1,...,t'}r''_{j}$.

Evidemment, \`a tout tel couple, on peut associer comme ci-dessus un sous-groupe parabolique ${\bf P}_{{\bf r}',{\bf r}''}$ et sa composante de Levi ${\bf M}_{{\bf r}',{\bf r}''}$. Un couple $({\bf m}',{\bf m}'')\in X_{{\bf m}}$  n'appartient pas toujours \`a l'ensemble $X$ car des composantes $m'_{j}$ ou $m''_{j}$ peuvent \^etre nulles. Mais, en supprimant les termes nuls, on obtient un \'el\'ement $({\bf r}',{\bf r}'')\in X$ et, \'evidemment, les paraboliques et groupes de Levi associ\'es \`a $({\bf m}',{\bf m}'')$ et $({\bf r}',{\bf r}'')$ sont les m\^emes. On peut donc r\'ecrire (2) sous la forme
$$(\pi_{{\bf m}})_{n',n''}=\sum_{({\bf r}',{\bf r}'')\in X }x_{{\bf m}}({\bf r}',{\bf r}'')Ind_{{\bf P}_{{\bf r}',{\bf r}''}}^{{\bf G}_{\sharp}}\circ res_{{\bf M}_{{\bf r}',{\bf r}''}}^{{\bf G}_{\sharp}}(\pi_{n',n''}),$$
o\`u $x_{{\bf m}}({\bf r}',{\bf r}'')$ est le nombre des \'el\'ements $({\bf m}',{\bf m}'')\in X_{{\bf m}}$  \'egaux \`a $({\bf r}',{\bf r}'')$, \`a des termes $m'_{j}$ ou $m''_{j}$ nuls pr\`es. En utilisant (1), on en d\'eduit
$$(-1)^{n(\pi)}D(\pi)_{n',n''}=\sum_{({\bf r}',{\bf r}'')\in X }x ({\bf r}',{\bf r}'')Ind_{{\bf P}_{{\bf r}',{\bf r}''}}^{{\bf G}_{\sharp}}\circ res_{{\bf M}_{{\bf r}',{\bf r}''}}^{{\bf G}_{\sharp}}(\pi_{n',n''}),$$
o\`u
$$x({\bf r}',{\bf r}'')=\sum_{{\bf m}}(-1)^{t({\bf m})}x_{{\bf m}}({\bf r}',{\bf r}'').$$
On prouvera plus loin l'\'egalit\'e
$$(3) \qquad x({\bf r}',{\bf r}'')=(-1)^{t'({\bf r}')+t''({\bf r}'')},$$
$t'({\bf r}')$ et $t''({\bf r}'')$ \'etant \'evidemment les entiers $t'$ et $t''$ figurant dans les donn\'ees ${\bf r}'$, ${\bf r}''$. En utilisant cela, on obtient
$$(4) \qquad (-1)^{n(\pi)}D(\pi)_{n',n''}=\sum_{({\bf r}',{\bf r}'')\in X } (-1)^{t'({\bf r}')+t''({\bf r}'')}Ind_{{\bf P}_{{\bf r}',{\bf r}''}}^{{\bf G}_{\sharp}}\circ res_{{\bf M}_{{\bf r}',{\bf r}''}}^{{\bf G}_{\sharp}}(\pi_{n',n''}).$$

Soit $\rho$ une composante irr\'eductible de $\pi_{n',n''}$. On a alors une formule similaire \`a (1):
$$(5) \qquad (-1)^{n(\rho)}{\bf D}(\rho)=\sum_{({\bf r}',{\bf r}'')\in X } (-1)^{t'({\bf r}')+t''({\bf r}'')}Ind_{{\bf P}_{{\bf r}',{\bf r}''}}^{{\bf G}_{\sharp}}\circ res_{{\bf M}_{{\bf r}',{\bf r}''}}^{{\bf G}_{\sharp}}(\rho).$$
Expliquons cela. La formule r\'esulte d'une formule similaire pour les groupes ${\bf SO}(2n'+1)$ et ${\bf O}(2n'')_{\sharp}$. Pour le premier et aussi le second quand $\sharp=an$, il n'y a pas de probl\`eme, c'est bien la formule de d\'efinition de l'involution.  Consid\'erons le cas d'un groupe ${\bf O}(2n'')_{iso}$, avec $n''>0$ (sinon le groupe dispara\^{\i}t). Traitons chacune des composantes, en commen\c{c}ant par ${\bf SO}(2n'')_{iso}$. Les sous-groupes paraboliques standard de ce groupe ne sont pas en bijection avec nos donn\'ees ${\bf r}''$. Pour param\'etrer ces sous-groupes, on doit d'une part ne consid\'erer que des ${\bf r}''=(r''_{1},...,r''_{t''},n''_{0})$ telles que $n''_{0}\not=1$ (car les sous-groupes paraboliques d\'efinis par $(r''_{1},...,r''_{t''},n''_{0}=1)$ et $(r''_{1},...,r''_{t''},1,n''_{0}=0)$ sont les m\^emes). D'autre part, si $n''_{0}=0$ et $r''_{t''}\geq2$, il y a deux sous-groupes paraboliques associ\'es \`a ${\bf r}''$, qui sont conjugu\'es par un \'el\'ement de ${\bf O}^{-}(2n'')_{iso}({\mathbb F}_{q})$ mais pas par un \'el\'ement de ${\bf SO}(2n'')_{iso}({\mathbb F}_{q})$. Consid\'erons le second probl\`eme, soit ${\bf r}''=(r''_{1},...,r''_{t''},n''_{0}=0)$ avec $r''_{t''}\geq2$. Dans la d\'efinition de l'involution interviennent les induites \`a ${\bf SO}(2n'')_{iso}({\mathbb F}_{q})$ \`a partir des deux sous-groupes paraboliques en question. D'apr\`es leur d\'efinition, les repr\'esentations que l'on induit sont conjugu\'ees par l'\'el\'ement de ${\bf O}^{-}(2n'')_{iso}({\mathbb F}_{q})$ qui \'echange les deux paraboliques. Dans la formule similaire \`a (5) intervient la restriction \`a ${\bf SO}(2n'')_{iso}({\mathbb F}_{q})$ de l'induite de ${\bf P}''_{{\bf r}''}({\mathbb F}_{q})$ \`a ${\bf O}(2n'')_{iso}({\mathbb F}_{q})$ de l'une ou l'autre de ces repr\'esentations. Mais, parce que $n''_{0}=0$,  ${\bf P}''_{{\bf r}''}$ ne coupe pas la composante ${\bf O}^-(2n'')_{iso}$ et  la restriction \`a ${\bf SO}(2n'')_{iso}({\mathbb F}_{q})$ de cette derni\`ere induite se d\'ecompose naturellement en deux induites qui ne sont autres que les pr\'ec\'edentes. Consid\'erons maintenant le cas d'une donn\'ee ${\bf r}''=(r''_{1},...,r''_{t''},n''_{0}=1)$. Introduisons l'autre donn\'ee $\bar{{\bf r}}''=(r''_{1},...,r''_{t''},1,n''_{0}=0)$. Dans la d\'efinition de l'involution intervient  le produit de $(-1)^{t''+1}$ et d'une induite de ${\bf P}''_{\bar{{\bf r}}''}({\mathbb F}_{q})$ \`a ${\bf SO}(2n'')_{iso}({\mathbb F}_{q})$ d'une certaine repr\'esentation. Par construction de celle-ci et parce que le normalisateur de ${\bf P}''_{\bar{{\bf r}}''}$ dans ${\bf O}^-(2n'')_{iso}({\mathbb F}_{q})$ est non vide, cette repr\'esentation est invariante par ce normalisateur.  Dans la formule similaire \`a (5) intervient d'une part le produit de $(-1)^{t''+1}$ et de la restriction \`a ${\bf SO}(2n'')_{iso}({\mathbb F}_{q})$ de l'induite de ${\bf P}''_{\bar{{\bf r}}''}({\mathbb F}_{q})$ \`a ${\bf O}(2n'')_{iso}({\mathbb F}_{q})$ de la m\^eme repr\'esentation. D'apr\`es la propri\'et\'e d'invariance ci-dessus et le fait que ${\bf P}''_{\bar{{\bf r}}''}$ ne coupe pas ${\bf O}^{-}(2n'')_{iso}$, on obtient $2$ fois la contribution de $\bar{{\bf r}}''$ \`a la  formule de l'involution. Mais, dans la formule similaire \`a (5), il y a aussi le produit de $(-1)^{t''}$ et de la restriction \`a ${\bf SO}(2n'')_{iso}({\mathbb F}_{q})$ de l'induite de ${\bf P}''_{{\bf r}''}({\mathbb F}_{q})$ \`a ${\bf O}(2n'')_{iso}({\mathbb F}_{q})$ de la m\^eme repr\'esentation, prolong\'ee en une repr\'esentation de ${\bf M}''_{{\bf r}''}({\mathbb F}_{q})$. Cette fois, ${\bf P}''_{{\bf r}''}$ coupe ${\bf O}^-(2n'')_{iso}$ et on obtient $-1$ fois  la contribution de $\bar{{\bf r}}''$ \`a la formule de l'involution. Comme $2-1=1$, les deux formules co\"{\i}ncident. La comparaison des formules pour la composante ${\bf O}^-(2n'')_{iso}$ est similaire. Cette fois, l'involuion est d\'efinie en induisant \`a partir des normalisateurs dans ${\bf O}^{-}(2n'')_{iso}({\mathbb F}_{q})$ des sous-groupes paraboliques de ${\bf SO}(2n'')_{iso}$ pour lesquels ce normalisateur est non vide. Les donn\'ees ${\bf r}''=(r''_{1},...,r''_{t''},n''_{0}=0)$ avec $r''_{t''}\geq2$ disparaissent, leurs deux sous-groupes paraboliques associ\'es ne v\'erifiant pas cette condition. Une telle donn\'ee intervient dans la formule similaire \`a (5) par la restriction \`a ${\bf O}^-(2n'')_{iso}({\mathbb F}_{q})$ d'une certaine induite de ${\bf P}''_{{\bf r}''}({\mathbb F}_{q})$ \`a ${\bf O}(2n'')_{iso}({\mathbb F}_{q})$. Parce que ${\bf P}''_{{\bf r}''}$ ne coupe pas ${\bf O}^-(2n'')_{iso}$, on voit que l'action de ${\bf O}^-(2n'')_{iso}({\mathbb F}_{q})$ permute deux sous-${\bf SO}(2n'')_{iso}({\mathbb F}_{q})$-modules de l'induite et une telle repr\'esentation de  ${\bf O}^-(2n'')_{iso}({\mathbb F}_{q})$ est de trace nulle. Enfin, pour une donn\'ee ${\bf r}''=(r''_{1},...,r''_{t''},n''_{0}=1)$, il intervient dans la formule de l'involution une induite \`a partir du normalisateur dans ${\bf O}^-(2n'')_{iso}({\mathbb F}_{q})$ de ${\bf P}''_{\bar{{\bf r}}''}$, avec la m\^eme notation que plus haut. Or ce normalisateur est justement l'intersection ${\bf O}^-(2n'')_{iso}({\mathbb F}_{q})\cap {\bf P}''_{{\bf r}''}({\mathbb F}_{q})$. On voit que les deux formules se r\'econcilient encore. Ces consid\'erations justifient la formule (5). 

En vertu de (4) et (5), pour prouver l'\'egalit\'e cherch\'ee ${\bf D}(\pi_{n',n''})=D(\pi)_{n',n''}$, il suffit de prouver que, pour toute composante irr\'eductible $\rho$ de $\pi_{n',n''}$, on a l'\'egalit\'e $(-1)^{n(\pi)}=(-1)^{n(\rho)}$. Cela r\'esulte de la d\'efinition par la m\'ethode de Lusztig des \'el\'ements de $Irr_{unip,\sharp}$. En fait, ces signes se calculent, ils valent $(-1)^n$ si $\sharp=iso$ et $(-1)^{n-1}$ si $\sharp=an$, cf. \cite{MW} corollaire 5.7.

Cela ach\`eve la d\'emonstration, \`a ceci pr\`es qu'il nous reste \`a prouver l'\'egalit\'e (3). Consid\'erons un couple ${\bf r}'=(r'_{1},...,r'_{t'},n'_{0})$, ${\bf r}''=(r''_{1},...,r''_{t''},n''_{0})$ appartenant \`a $X$. Les donn\'ees $t$, ${\bf m}'$ et ${\bf m}''$ intervenant dans la d\'efinition de $x({\bf r}',{\bf r}'')$ se construisent de la fa\c{c}on suivante. On consid\`ere un entier $t\geq sup(t',t'')$ et deux sous-ensembles $Y',Y''\subset \{1,...,t\}$. On suppose que $Y'$, resp. $Y''$, a $t-t'$ \'el\'ements, resp. $t-t''$. Notons $i_{1},...,i_{t'}$ les \'el\'ements de $\{1,...,t\}-Y'$. On d\'efinit ${\bf m}'$ par $m'_{i_{k}}=r'_{k}$ pour $k=1,...,t'$ et $m'_{j}=0$ pour $j\in Y'$. On d\'efinit de fa\c{c}on similaire ${\bf m}''$, en utilisant l'ensemble $Y''$ \`a la place de $Y'$. Pour tout $j=1,...,t$ le terme $m_{j}=m'_{j}+m''_{j}$ doit \^etre non nul. Cela \'equivaut \'evidemment \`a la condition $Y'\cap Y''=\emptyset$. On voit qu'elle ne peut \^etre r\'ealis\'ee que si $t\leq t'+t''$. Pour $t\in \{sup(t',t''),...,t'+t''\}$, on obtient que
$$\sum_{{\bf m}; l({\bf m})=t} x_{{\bf m}}({\bf r}',{\bf r}'')$$
est le nombre de couples $(Y',Y'')$ comme ci-dessus. Ce nombre est le produit du nombre de sous-ensembles $Y'$ \`a $t-t'$ \'el\'ements de $\{1,...,t\}$ et du nombre de sous-ensembles $Y''$ \`a $t-t''$ \'el\'ements de $\{1,...,t\}-Y'$. C'est-\`a-dire
$$\frac{t!}{t'!(t-t')!}\frac{t'!}{(t-t'')!(t'+t''-t)!}=\frac{t!}{(t-t')!(t-t'')!(t'+t''-t)!}.$$
D'o\`u
$$x({\bf r}',{\bf r}'')=\sum_{t=sup(t',t'')...,t'+t''}(-1)^t\frac{t!}{(t-t')!(t-t'')!(t'+t''-t)!}.$$
Supposons $t'\geq t''$ pour fixer la notation. On remplace $t$ par $t'+u$. Alors $x({\bf r}',{\bf r}'')$ est la valeur en $X=1$ du polyn\^ome
$$P(X)=\sum_{u=0,...,t''}(-1)^{t'+u}\frac{(t'+u)!}{u!(t'-t''+u)!(t''-u)!}X^{t'-t''+u}.$$
En posant 
$$Q(X)=\sum_{u=0,...,t''}(-1)^{t'+u}\frac{ 1}{u! (t''-u)!}X^{t'+u},$$
on voit que $P(X)=(\frac{d}{dX})^{t''}Q(X)$. Or on calcule
$$Q(X)=\frac{(-1)^{t'}}{t''!}X^{t'}(1-X)^{t''}.$$
On voit que $(\frac{d}{dX})^{t''}Q(X)$ est la somme de $(-1)^{t'+t''}$ et d'un polyn\^ome divisible par $(1-X)$. Donc la valeur en $1$ de $(\frac{d}{dX})^{t''}Q(X)$, c'est-\`a-dire de $P(X)$, vaut $(-1)^{t'+t''}$. Cela prouve (3). $\square$

\subsection{Les espaces ${\cal R}$ et ${\cal R}^{glob}$}

Pour tout groupe fini $W$, on note $\hat{W}$ l'ensemble des classes de repr\'esentations irr\'eductibles de $W$. En identifiant une telle repr\'esentation \`a son caract\`ere, l'espace ${\mathbb C}[\hat{W}]$ s'identifie \`a celui des fonctions de $W$ dans ${\mathbb C}$ qui sont invariantes par conjugaison. 

Soit $N\in {\mathbb N}$. On note $\mathfrak{S}_{N}$ le groupe des permutations de l'ensemble $\{1,...,N\}$.   On note $sgn$ le caract\`ere signe usuel de $\mathfrak{S}_{N}$.    Les classes de conjugaison dans $\mathfrak{S}_{N}$ sont param\'etr\'ees par ${\cal P}(N)$. On note $w_{\alpha}$ un \'el\'ement de la classe ainsi param\'etr\'ee.

Pour un entier $k\geq1$, notons ${\cal P}_{k}(N)$ l'ensemble des familles $(\alpha_{1},...,\alpha_{k})$ de partitions telles que $S(\alpha_{1})+...+S(\alpha_{k})=N$.  Dans la suite, on utilisera des variantes de cette notation, par exemple ${\cal P}_{k}^{symp}(N)$. On note $W_{N}$ le groupe de Weyl d'un syst\`eme de racines de type $B_{N}$ ou $C_{N}$ (avec la convention $W_{0}=\{1\}$).  On note $sgn$ le caract\`ere signe  usuel de $W_{N}$ et $sgn_{CD}$ le caract\`ere dont le noyau est le sous-groupe $W_{N}^D$ d'un syst\`eme de racines de type $D_{N}$.      Les classes de conjugaison dans $W_{N}$ sont param\'etr\'ees par les couples de partitions $(\alpha,\beta) \in {\cal P}_{2}(N)$. On note $w_{\alpha\beta}$ un \'el\'ement de la classe ainsi param\'etr\'ee.   On a $sgn(w_{\alpha\beta})=(-1)^{N+l(\alpha)}$ et $sgn_{CD}(w_{\alpha\beta})=(-1)^{l(\beta)}$.

  On note $\Gamma$ l'ensemble des quadruplets $\gamma=(r',r'',N',N'')$ tels que
$$r'\in {\mathbb N},\, r''\in {\mathbb Z},\, N'\in {\mathbb N}, \,N''\in {\mathbb N},\, r^{_{'}2}+r'+N'+r^{_{''}2}+N''=n.$$
Pour un tel $\gamma $, on pose
$${\cal R}(\gamma)={\mathbb C}[\hat{W}_{N'}]\otimes {\mathbb C}[\hat{W}_{N''}].$$
On d\'efinit
$${\cal R}=\oplus_{\gamma\in \Gamma}{\cal R}(\gamma).$$

 Chaque espace ${\mathbb C}[\hat{W}_{N}]$ est naturellement muni d'un produit hermitien d\'efini positif. On en d\'eduit un tel produit sur l'espace ${\cal R}$. 

 Soit $m\in {\mathbb N}$ avec $1\leq m\leq n$.  Pour un entier $N \geq n$, le groupe $\mathfrak{S}_{m}\times W_{N-m}$ appara\^{\i}t comme sous-groupe "de Levi" de $W_{N}$ et il y a une application lin\'eaire de restriction
 $$res_{m}:{\mathbb C}[\hat{W}_{N}]\to {\mathbb C}[\hat{\mathfrak{S}}_{m}]\otimes {\mathbb C}[\hat{W}_{N-m}].$$
 Pour $\varphi\in {\mathbb C}[\hat{W}_{N}]$, $\mu\in {\cal P}(m)$ et $(\alpha,\beta)\in {\cal P}_{2}(N-m)$, on a l'\'egalit\'e 

$res_{m}(\varphi)(w_{\mu}\times w_{\alpha,\beta})=\varphi(w_{\mu\cup \alpha,\beta})$. 

On construit des applications lin\'eaires
$$res_{m}',res_{m}'':{\cal R}\to {\mathbb C}[\hat{\mathfrak{S}}_{m}]\otimes {\cal R}_{n-m}$$
de la fa\c{c}on suivante. Soit $\gamma=(r',r'',N',N'')\in \Gamma_{n}$. Si $N'< m$, $res_{m}'$ est nulle sur ${\cal R}(\gamma)$. Si $N'\geq m$, $\gamma'=(r',r'',N'-m,N'')$ appartient \`a $\Gamma_{n-m}$. L'application $res'_{m}$ envoie ${\cal R}(\gamma)$ dans ${\cal R}(\gamma')$ et co\"{\i}ncide sur ${\cal R}(\gamma)$ avec le produit tensoriel de 
$$res_{m}:{\mathbb C}[\hat{W}_{N'}]\to {\mathbb C}[\hat{\mathfrak{S}}_{m}]\otimes {\mathbb C}[\hat{W}_{N'-m}]$$
et de l'identit\'e de ${\mathbb C}[\hat{W}_{N''}]$. L'application $res''_{m}$ est similaire, en permutant les r\^oles des $'$ et $''$. On d\'efinit ${\cal R}^{glob}$ comme l'ensemble des \'el\'ements de ${\cal R}$ qui, pour tout $m$,  ont m\^eme image par les deux applications $res'_{m}$, $res''_{m}$. On note 
$$res_{m}:{\cal R}^{glob}\to  {\mathbb C}[\hat{\mathfrak{S}}_{m}]\otimes{\cal R}_{n-m}$$
l'application commune $res'_{m}=res''_{m}$.  On voit que $res_{m}$ envoie ${\cal R}^{glob}$ dans $ {\mathbb C}[\hat{\mathfrak{S}}_{m}]\otimes {\cal R}^{glob}_{n-m}$ (la preuve est similaire \`a celle de 1.5 concernant l'espace ${\cal R}^{par,glob}$). 

On note ${\cal R}_{cusp}$ le sous-espace des \'el\'ements de ${\cal R}$ annul\'es par toutes les applications $res'_{m}$ et $res''_{m}$. Il est inclus dans ${\cal R}^{glob}$. On note $proj_{cusp}$ la projection orthogonale de ${\cal R}$ sur ${\cal R}_{cusp}$. L'espace ${\cal R}_{cusp}$ se d\'ecrit de la fa\c{c}on suivante. Pour tout $N\in {\mathbb N}$, notons ${\mathbb C}[\hat{W}_{N}]_{cusp}$ le sous-espace des fonctions sur $W_{N}$, invariantes par conjugaison, et \`a support  dans les classes de conjugaison param\'etr\'ees par des couples de partitions  de la forme $(\emptyset,\beta)$. Pour tout $\gamma\in \Gamma$, on pose
$${\cal R}_{cusp}(\gamma)={\mathbb C}[\hat{W}_{N'}]_{cusp}\otimes {\mathbb C}[\hat{W}_{N''}]_{cusp}.$$
Alors
$${\cal R}_{cusp}=\oplus_{\gamma\in \Gamma}{\cal R}_{cusp}(\gamma).$$

Donnons une autre pr\'esentation de l'espace ${\cal R}$. 
Notons $\boldsymbol{\Gamma}$ l'ensemble des triplets $(r',r'',N)$ tels que
$$r'\in {\mathbb N},\,\, r''\in {\mathbb Z},\,\, N\in {\mathbb N},\,\, r^{'2}+r'+r^{''2}+N=n.$$
Il y a une application \'evidente $\Gamma\to \boldsymbol{\Gamma}$ qui, \`a $(r',r'',N',N'')\in \Gamma$, associe $(r',r'',N'+N'')$. On la note $\gamma\mapsto \boldsymbol{\gamma}$. Pour $\boldsymbol{\gamma}=(r',r'',N)\in \boldsymbol{\Gamma}$, posons 
$${\cal R}(\boldsymbol{\gamma})=\oplus_{\gamma\in \Gamma; \gamma\mapsto \boldsymbol{\gamma}}{\cal R}(\gamma)=\sum_{(N',N'')\in D(N)}{\mathbb C}[\hat{W}_{N'}]\otimes {\mathbb C}[\hat{W}_{N''}].$$
Alors 
$${\cal R}=\oplus_{\boldsymbol{\gamma}\in \boldsymbol{\Gamma}}{\cal R}(\boldsymbol{\gamma}).$$

Soit $\boldsymbol{\gamma}=(r',r'',N)\in \boldsymbol{\Gamma}$ et soit $\varphi\in {\cal R}(\boldsymbol{\gamma})$. On peut identifier $\varphi$ \`a une fonction sur l'ensemble des $w_{\alpha',\beta'}\times w_{\alpha'',\beta''}$ pour $(\alpha',\beta',\alpha'',\beta'')\in {\cal P}_{4}(N)$. D'apr\`es la formule explicite \'ecrite plus haut pour l'application $res_{m}$, la condition $res'_{m}(\varphi)=res''_{m}(\varphi)$ signifie que, pour tout $\mu\in {\cal P}(m)$ et tout 
$(\alpha',\beta',\alpha'',\beta'')\in {\cal P}_{4}(N-m)$, on a l'\'egalit\'e $\varphi(w_{\mu\cup \alpha',\beta'}\times w_{\alpha'',\beta''})=\varphi(w_{\alpha',\beta'}\times w_{\mu\cup \alpha'',\beta''})$. En faisant varier $m$, on voit que $\varphi\in {\cal R}^{glob}$ si et seulement si, pour tout $(\alpha',\beta',\alpha'',\beta'')\in {\cal P}_{4}(N)$, $\varphi(w_{\alpha',\beta'}\times w_{\alpha'',\beta''})$ ne d\'epend que du triplet $(\alpha'\cup \alpha'',\beta',\beta'')$. 
On note ${\cal R}^{glob}(\boldsymbol{\gamma})$ l'espace des $\varphi\in {\cal R}(\boldsymbol{\gamma})$ qui v\'erifient cette condition. Pour $(\alpha,\beta',\beta'')\in {\cal P}_{3}(N)$, on note $w_{\alpha,\beta',\beta''}$ un \'el\'ement quelconque $w_{\alpha',\beta'}\times w_{\alpha'',\beta''}$ o\`u $\alpha'\cup \alpha''=\alpha$. On peut consid\'erer que ${\cal R}^{glob}(\boldsymbol{\gamma})$ est l'espace des fonctions sur l'ensemble de ces \'el\'ements $w_{\alpha,\beta',\beta''}$. On a l'\'egalit\'e
$${\cal R}^{glob}=\oplus_{\boldsymbol{\gamma}\in \boldsymbol{\Gamma}}{\cal R}^{glob}(\boldsymbol{\gamma}).$$

  \bigskip
  
  \subsection{L'involution de Lusztig}
  En suivant Lusztig, on a d\'efini en \cite{MW} 3.16 deux isomorphismes 
$$Rep:{\cal R}\to {\cal R}^{par},\,\,k:{\cal R}\to {\cal R}^{par}.$$
Donnons seulement une id\'ee des d\'efinitions, en renvoyant \`a \cite{MW} 2.6, 2.7, 2.9, 2.10 pour plus de pr\'ecision. Lusztig a classifi\'e de fa\c{c}on combinatoire les repr\'esentations irr\'eductibles unipotentes d'un groupe ${\bf SO}(2n'+1;{\mathbb F}_{q})$. Elles sont param\'etr\'ees par les couples $(r',N')\in {\mathbb N}^2$ tels que $r^{_{'}2}+r'+N'=n'$. De m\^eme, la r\'eunion disjointe des ensembles de repr\'esentations irr\'eductibles unipotentes de ${\bf O}(2n'')_{iso}({\mathbb F}_{q})$ et de  ${\bf O}(2n'')_{an}({\mathbb F}_{q})$ sont param\'etr\'ees par les couples $(r'',N'')\in {\mathbb Z}\times {\mathbb N}$ tels que $r^{_{''}2}+N''=n''$. Pour $\gamma=(r',r'',N',N'')\in \Gamma$, pour $\rho'\in \hat{W}_{N'}$ et $\rho''\in \hat{W}_{N''}$, l'isomorphisme $Rep$ envoie l'\'el\'ement $\rho'\otimes \rho''\in {\cal R}(\gamma)$ sur un \'el\'ement de $C'_{n'}\otimes C''_{n''}$ o\`u $n'=r^{_{'}2}+r'+N'=n'$ et  $n''=r^{_{''}2}+N''$. Cet \'el\'ement est le produit tensoriel des traces des repr\'esentations irr\'eductibles unipotentes param\'etr\'ees par $(r',\rho')$ et $(r'',\rho'')$. 

Lusztig a  aussi introduit la notion de faisceau-caract\`ere. Pour un tel faisceau d\'efini par exemple sur un groupe ${\bf SO}(2n'+1)$, la fonction trace associ\'ee \`a ce faisceau est une fonction sur ${\bf SO}(2n'+1;{\mathbb F}_{q})$, invariante par conjugaison. Lusztig a classifi\'e les faisceaux-caract\`eres  dont la fonction trace est unipotente, c'est-\`a-dire appartient \`a l'espace $C'_{n'}$. De nouveau, ils sont param\'etr\'es par les couples $(r',N')\in {\mathbb N}^2$ tels que $r^{_{'}2}+r'+N'=n'$. Une construction analogue vaut pour les groupes ${\bf O}(2n'')_{iso}$ ou ${\bf O}(2n'')_{an}$. Pour $\gamma=(r',r'',N',N'')\in \Gamma$, pour $\rho'\in \hat{W}_{N'}$ et $\rho''\in \hat{W}_{N''}$, l'isomorphisme $k$ envoie l'\'el\'ement $\rho'\otimes \rho''\in {\cal R}(\gamma)$ sur un \'el\'ement de $C'_{n'}\otimes C''_{n''}$ o\`u $n'$ et $n''$ sont comme ci-dessus. Cet \'el\'ement est le produit tensoriel des traces des faisceaux-caract\`eres param\'etr\'es par $(r',\rho')$ et $(r'',\rho'')$. 

{\bf Remarque.} La  fonction-trace associ\'ee \`a un faisceau-caract\`ere n'est vraiment canonique qu'\`a homoth\'etie pr\`es. Pour la d\'efinir pr\'ecis\'ement, on doit faire des choix de normalisations. Nous utilisons ceux de \cite{MW} 2.7 et 2.10. 

\bigskip
Notons ${\cal F}^L$ l'automorphisme de ${\cal R}$ tel que $Rep\circ {\cal F}^L=k$. C'est une isom\'etrie.  Lusztig a prouv\'e que c'\'etait une involution, qui peut se d\'ecrire de fa\c{c}on combinatoire. On transporte ${\cal F}^L$ en une involution ${\cal F}^{par}$ de ${\cal R}^{par}$. Autrement dit, ${\cal F}^{par}$ est l'involution de ${\cal R}^{par}$ telle que ${\cal F}^{par}\circ Rep=k$. Elle est isom\'etrique.

Pour tout entier $m\geq1$,  on a de m\^eme des isomorphismes $Rep,k:{\mathbb C}[\hat{\mathfrak{S}}_{{\bf m}}]\to C^{GL(m)}$. Cette fois, ils sont \'egaux et les analogues de ${\cal F}^L$ et ${\cal F}^{par}$ sont les identit\'es. Supposons $m\leq n$, posons $n_{0}=n-m$.  On a construit des applications lin\'eaires
$$res'_{m},res''_{m}:{\cal R}^{par}\to  C^{GL(m)}\otimes {\cal R}^{par}_{n_{0}},$$
$$res_{ m}',res_{m}'':{\cal R}\to {\mathbb C}[\hat{\mathfrak{S}}_{m}]\otimes {\cal R}_{n_{0}}.$$
 Les applications $Rep$ et $k$ sont compatibles en un sens plus ou moins \'evident \`a ces applications. Il en r\'esulte que ${\cal F}^L$  et ${\cal F}^{par}$ le sont aussi.

Les applications $Rep$ et $k$ se restreignent en des isomorphismes ${\cal R}^{glob}\to {\cal R}^{par,glob}$. L'involution ${\cal F}^L$ conserve ${\cal R}^{glob}$ et ${\cal F}^{par}$ conserve ${\cal R}^{par,glob}$. Les propri\'et\'es de compatibilit\'e ci-dessus  de l'involution isom\'etrique ${\cal F}^{par}$ entra\^{\i}nent que:

(1) ${\cal F}^{par}\circ proj_{cusp}=proj_{cusp}\circ {\cal F}^{par}$.

\bigskip

\subsection{L'induction endoscopique}
 On a d\'efini en \cite{MW} 3.18 une application lin\'eaire $\rho\iota:{\cal R}\to {\cal R}^{glob}$.  Rappelons sa d\'efinition. Soit $\gamma=(r',r'',N^+,N^-)\in \Gamma$ et $\varphi\in {\cal R}(\gamma) $. Posons $N=N^++N^-$. L'\'el\'ement $\rho\iota(\varphi)$ appartient \`a ${\cal R}(\boldsymbol{\delta})$, o\`u $\boldsymbol{\delta}=(r',(-1)^{r'}r'',N)\in \boldsymbol{\Gamma}$.   Soit $\delta=(r',(-1)^{r'},N_{1},N_{2})\in \Gamma$. Nous allons d\'ecrire la composante $\rho\iota(\varphi)_{\delta}$ de $\rho\iota(\varphi)$ dans ${\cal R}(\delta)$.
    
  On d\'efinit un quadruplet d'entiers ${\bf a}=(a_{1}^+,a_{1}^-,a_{2}^+,a_{2}^-)$ par les formules suivantes:

${\bf a}=(0,0,0,1)$  si $0< r''\leq r'$ ou si  $r''=0$ et $r'$ est pair;

${\bf a}=(0,0,1,0)$ si  $-r'\leq r''<0$ ou si  $r''=0$ et $r'$ est impair;

${\bf a}=(0,1,0,0)$ si $r'<r''$;

${\bf a}=(1,0,0,0)$ si   $r''<-r'$.
 
  Notons ${\cal N}$ l'ensemble des quadruplets ${\bf N}=(N^+_{1},N^-_{1},N^+_{2},N^-_{2})$ d'entiers positifs ou nuls tels que
  $$N^+= N^+_{1}+N^+_{2},\,\,N^-=N^-_{1}+N^-_{2},\,\,N_{1}=N^+_{1}+N^-_{1},\,\,N_{2}=N^+_{2}+N^-_{2}.$$
 Pour un tel quadruplet, 
posons $W_{{\bf N}}=W_{N^+_{1}}\times W_{N^-_{1}}\times W_{N^+_{2}}\times W_{N^-_{2}}$. Ce groupe se plonge de fa\c{c}on \'evidente dans $W_{N_{1}}\times W_{N_{2}}$, resp. $W_{N^+}\times W_{N^-}$, et ces plongements sont bien d\'efinis \`a conjugaison pr\`es. On a donc des foncteurs de restriction $res^{W_{N^+}\times W_{N^-}}_{W_{{\bf N}}}$ et d'induction $ind^{W_{N_{1}}\times W_{N_{2}}}_{W_{{\bf N}}}$. On note $sgn_{CD}^{{\bf a}}$ le caract\`ere de $W_{{\bf N}}$ qui est le produit tensoriel des caract\`eres $sgn_{CD}^{a^+_{1}}$, $sgn_{CD}^{a^-_{1}}$, $sgn_{CD}^{a^+_{2}}$, $sgn_{CD}^{a^-_{2}}$ sur chacun des facteurs de $W_{{\bf N}}$. Alors 
$$\rho\iota(\varphi)_{\delta}=\sum_{{\bf N}\in {\cal N}}ind^{W_{N_{1}}\times W_{N_{2}}}_{W_{{\bf N}}}\left(sgn_{CD}^{{\bf a}}\otimes res^{W_{N^+}\times W_{N^-}}_{W_{{\bf N}}}( \varphi )\right).$$
 
Donnons une formule plus concr\`ete. Fixons plut\^ot $\boldsymbol{\gamma}=(r',r'',N)\in \boldsymbol{\Gamma}$ et supposons $\varphi\in {\cal R}(\boldsymbol{\gamma})$. On a $\rho\iota(\varphi)\in {\cal R}(\boldsymbol{\delta})$, o\`u $\boldsymbol{\delta}$ est comme ci-dessus.  Soit $(\alpha,\beta_{1},\beta_{2})\in {\cal P}_{3}(N)$. Notons ${\cal I}$ l'ensemble des sextuplets ${\bf I}=(I^+,I^-,J^+_{1},J^-_{1},J^+_{2},J^-_{2})$, dont les composantes sont des ensembles tels que
$$\{1,...,l(\alpha)\}=I^+\sqcup I^-,\,\, \{1,...,l(\beta_{1})\}=J^+_{1}\sqcup J^-_{1}, \,\{1,...,l(\beta_{2})\}=J^+_{2}\sqcup J^-_{2}.$$
A tout tel sextuplet, on associe six partitions $\alpha^+({\bf I}), \alpha^-({\bf I}),\beta^+_{1}({\bf I}),\beta^-_{1}({\bf I}),\beta^+_{2}({\bf I}),\beta^-_{2}({\bf I}) $. La partition $\alpha^+({\bf I})$ est form\'ee des $\alpha_{j}$ pour $j\in I^+$. Les autres sont construites de fa\c{c}on similaire. 
On v\'erifie la formule suivante
$$(1) \qquad  \rho\iota(\varphi)(w_{\alpha,\beta_{1},\beta_{2}})=\sum_{{\bf I}\in {\cal I}}(-1)^{a^+_{1}\vert J^+_{1}\vert +a^-_{1}\vert J^-_{1}\vert +a^+_{2}\vert J^+_{2}\vert +a^-_{2}\vert J^-_{2}\vert }$$
$$ \varphi(w_{\alpha^+({\bf I}),\beta^+_{1}({\bf I})\cup \beta^+_{2}({\bf I})}\times w_{\alpha^-({\bf I}),\beta^-_{1}({\bf I})\cup\beta^-_{2}({\bf I})}).$$

 On voit que notre application $\rho\iota$ commute aux projections $proj_{cusp}$.

\bigskip

\subsection{Description de $Res\circ D(\pi)$ pour $\pi\in Irr_{unip-quad}$}

On d\'efinit un espace  
$${\cal S}_{n}=\oplus_{k,N}{\mathbb C}[\hat{W}_{N}],$$
o\`u l'on somme sur les couples $(k,N)\in {\mathbb N}^2$ tels que $k(k+1)+2N=2n$. On a l'\'egalit\'e
$$(1) \qquad \oplus_{(n^+,n^-)\in D(n)}{\cal S}_{n^+}\otimes {\cal S}_{n^-}= \oplus_{k^+,N^+,k^-,N^-}{\mathbb C}[\hat{W}_{N^+}]\otimes {\mathbb C}[\hat{W}_{N^-}],$$
o\`u l'on somme sur l'ensemble $\Gamma_{{\cal S}}$ des quadruplets $(k^+,N^+,k^-,N^-)\in {\mathbb N}^4$ tels que  $k^+(k^++1)+2N^++k^-(k^-+1)+2N^-=2n$. Pour un tel quadruplet, d\'efinissons des entiers $r'\in {\mathbb N}$, $r''\in {\mathbb Z}$   par les formules suivantes:

- si $k^+\equiv k^-\,mod\,2{\mathbb Z}$, $r'=\frac{k^++k^-}{2}$, $r''=\frac{k^+-k^-}{2}$;

- si $k^+ \not\equiv k^-\,mod\,2{\mathbb Z}$ et $k^+>k^-$, $r'=\frac{k^+-k^--1}{2}$, $r''=\frac{k^++k^-+1}{2}$;

- si $k^+\not\equiv k^-\,mod\,2{\mathbb Z}$ et $k^+< k^-$, $r'=\frac{k^--k^+-1}{2}$, $r^-=-\frac{k^++k^-+1}{2}$.

On v\'erifie que $(r',r'',N^+,N^-)\in \Gamma$. Inversement, pour $(r',r'',N^+,N^-)\in \Gamma$, d\'efinissons deux entiers $k^+,k^-\in {\mathbb N}$ par les formules suivantes:

- si $\vert r''\vert \leq r'$, $k^+=r'+r''$, $k^-=r'-r''$;

- si $r'< r'' $, $k^+=r'+r''$, $k^-=r''-r'-1$;

- si $r''< -r' $, $k^+=-r''-r'-1$, $k^-=r'- r''$.

Alors $(k^+,k^-,N^+,N^-)\in \Gamma_{{\cal S}}$. Ces deux applications sont inverses l'une de l'autre et d\'efinissent donc des bijections entre $\Gamma_{{\cal S}}$ et $\Gamma$. On peut donc remplacer l'ensemble de sommation  $\Gamma_{{\cal S}}$  par $\Gamma$ dans la formule (1) et on obtient un isomorphisme
  $$j:\oplus_{(n^+,n^-)\in D(n)}{\cal S}_{n^+}\otimes {\cal S}_{n^-}\to {\cal R}.$$
  
 Puisque ${\cal S}_{0}={\mathbb C}$, l'espace ${\cal S}_{n}={\cal S}_{n}\otimes {\cal S}_{0}$ se plonge dans l'espace de d\'epart de $j$. On note $j_{n,0}$ la restriction de $j$ \`a ce sous-espace.

Soit $(\lambda^+,\epsilon^+,\lambda^-,\epsilon^-)\in \mathfrak{Irr}_{unip-quad}$. Posons $S(\lambda^+)=2n^+$, $S(\lambda^-)=2n^-$.  Soit $\zeta=\pm$. Par la correspondance de Springer g\'en\'eralis\'ee,  le couple $(\lambda^{\zeta},\epsilon^{\zeta})$ d\'etermine deux entiers $k^{\zeta} ,N^{\zeta}\in {\mathbb N}$ tel que $k^{\zeta}(k^{\zeta}+1)+2N^{\zeta}= 2n^{\zeta}$ et une repr\'esentation irr\'eductible $\rho_{\lambda^{\zeta},\epsilon^{\zeta}}$ de $W_{N^{\zeta}}$. On d\'efinit  une autre repr\'esentation $\boldsymbol{\rho}_{\lambda^{\zeta},\epsilon^{\zeta}}$ de $W_{N^{\zeta}}$, cf. \cite{W1} 5.1. Nous ne rappelons pas sa d\'efinition. Disons seulement que $\rho_{\lambda^{\zeta},\epsilon^{\zeta}}$ est l'action de $W_{N}$ dans un certain sous-espace d\'etermin\'e par $\epsilon^{\zeta}$ de l'espace de cohomologie de plus haut degr\'e d'une certaine vari\'et\'e alg\'ebrique, tandis que $\boldsymbol{\rho}_{\lambda^{\zeta},\epsilon^{\zeta}}$ est l'action de $W_{N}$ dans un sous-espace analogue de la somme de tous les espaces de cohomologie de la m\^eme vari\'et\'e. Cette repr\'esentation  n'est  pas irr\'eductible en g\'en\'eral. Elle est de la forme
$$\boldsymbol{\rho}_{\lambda^{\zeta},\epsilon^{\zeta}}=\oplus_{(\lambda_{1}^{\zeta},\epsilon^{\zeta}_{1})}c(\lambda^{\zeta},\epsilon^{\zeta};\lambda_{1}^{\zeta},\epsilon^{\zeta}_{1})\rho_{\lambda^{\zeta}_{1},\epsilon^{\zeta}_{1}},$$
o\`u $(\lambda_{1}^{\zeta},\epsilon^{\zeta}_{1})$ parcourt les \'el\'ements de $\boldsymbol{{\cal P}}^{symp}(2n^{\zeta})$ tels que les entiers $k_{1}^{\zeta},N_{1}^{\zeta}$ qui leur sont associ\'es sont \'egaux \`a $k^{\zeta},N^{\zeta}$ et o\`u  $c(\lambda^{\zeta},\epsilon^{\zeta};\lambda_{1}^{\zeta},\epsilon^{\zeta}_{1})\in {\mathbb N}$ est une multiplicit\'e. On sait que si $c(\lambda^{\zeta},\epsilon^{\zeta};\lambda_{1}^{\zeta},\epsilon^{\zeta}_{1})\geq1$, on a soit $\lambda_{1}^{\zeta}> \lambda^{\zeta}$ pour l'ordre usuel des partitions, soit $(\lambda_{1}^{\zeta},\epsilon_{1}^{\zeta})=(\lambda^{\zeta},\epsilon^{\zeta})$ et, dans ce dernier cas, on a $c(\lambda^{\zeta},\epsilon^{\zeta};\lambda^{\zeta},\epsilon^{\zeta})=1$.

  On identifie le  produit $\boldsymbol{\rho}_{\lambda^+,\epsilon^+}\otimes \boldsymbol{\rho}_{\lambda^-,\epsilon^-}$ \`a un \'el\'ement du membre de droite de (1), plus pr\'ecis\'ement de la composante index\'ee par $k^+,N^+,k^-,N^-$. On a donc $j(\boldsymbol{\rho}_{\lambda^+,\epsilon^+}\otimes \boldsymbol{\rho}_{\lambda^-,\epsilon^-})\in {\cal R}$ et on d\'efinit l'\'el\'ement
  $\rho\iota\circ j(\boldsymbol{\rho}_{\lambda^+,\epsilon^+}\otimes \boldsymbol{\rho}_{\lambda^-,\epsilon^-})\in {\cal R}$.

\ass{Proposition}{ On a l'\'egalit\'e 
$$Res\circ D(\pi(\lambda^+,\epsilon^+,\lambda^-,\epsilon^-))=Rep\circ \rho\iota\circ j(\boldsymbol{\rho}_{\lambda^+,\epsilon^+}\otimes \boldsymbol{\rho}_{\lambda^-,\epsilon^-}).$$}

Cf. \cite{W1} 5.3. 

\bigskip

\subsection{Un lemme technique}

 Soit  ${\bf m}=(m_{1},...,m_{t}>0)\in {\cal P}(\leq n)$. Posons $n_{0}({\bf m})=n-S({\bf m})$. En it\'erant les constructions de 1.8, on d\'efinit l'application lin\'eaire
$$res_{{\bf m}}:{\cal R}\to {\mathbb C}[\hat{\mathfrak{S}}_{m_{1}}]\otimes...\otimes {\mathbb C}[\hat{\mathfrak{S}}_{m_{t}}]\otimes {\cal R}_{n_{0}({\bf m})}.$$
Pour $i=1,...,t$, ${\mathbb C}[\hat{\mathfrak{S}}_{m_{i}}]_{cusp}$ est de dimension $1$. On choisit pour base la fonction caract\'eristique de la classe de conjugaison dans   $\mathfrak{S}_{m_{i}}$ param\'etr\'ee par la partition $(m_{i})$. On peut alors consid\'erer que $proj_{cusp}\circ res_{{\bf m}}$ prend ses valeurs dans ${\cal R}_{n_{0}({\bf m}),cusp}$. 

Soit $(n^+,n^-)\in D(n)$. Pour $\zeta=\pm$, soit $\phi^{\zeta}\in {\cal S}_{n^{\zeta}}$, posons $\varphi^{\zeta}=\rho\iota\circ j_{n^{\zeta},0}(\phi^{\zeta})$, cf. 1.11 pour la d\'efinition de $j_{n^{\zeta},0}$. Dans ce paragraphe, l'application $\rho\iota$ est  relative \`a divers entiers, ici c'est   $n^{\zeta}$.  Pour toute partition   ${\bf m}\in {\cal P}(\leq n^{\zeta})$, soit $\phi^{\zeta}[{\bf m}]\in {\cal S}_{n^{\zeta}_{0}({\bf m})}$, o\`u $n^{\zeta}_{0}({\bf m})=n^{\zeta}-S({\bf m})$. Posons $\varphi^{\zeta}[{\bf m}]=\rho\iota\circ j_{n^{\zeta}_{0}({\bf m}),0}(\phi^{\zeta}[{\bf m}])$.  Posons $\varphi=\rho\iota\circ j(\phi^+\otimes \phi^-)$. Fixons maintenant une partition ${\bf m}=(m_{1},...,m_{t}>0)\in {\cal P}(\leq n)$  et posons $n_{0}=n_{0}({\bf m})$.  Consid\'erons un couple d'ensembles $(I^+,I^-)$ tels que $I^+\sqcup I^-=\{1,...,t\}$. Pour tout tel couple, on d\'efinit les partitions ${\bf m}(I^+)$, resp.  ${\bf m}(I^-)$, form\'ees des $m_{i}$ pour $i\in I^+$, resp. $i\in I^-$. On note $\mathfrak{I}$ l'ensemble des couples $(I^+,I^-)$ comme ci-dessus tels que $S({\bf m}(I^+))\leq n^+$ et $S({\bf m}(I^-))\leq n^-$. Pour un tel couple, posons
$$\varphi[I^+,I^-]=\rho\iota\circ j_{n_{0}}(\phi^{+}[{\bf m}(I^+)]\otimes \phi^-[{\bf m}(I^-)]).$$

\ass{Lemme} {Supposons que, pour tout $\zeta=\pm$ et pour toute partition ${\bf m}'\in {\cal P}(\leq n^{\zeta})$, on ait l'\'egalit\'e
$$proj_{cusp}\circ res_{{\bf m}'}(\varphi^{\zeta})=proj_{cusp}(\varphi^{\zeta}[{\bf m}']).$$
Alors on a l'\'egalit\'e
$$proj_{cusp}\circ res_{{\bf m}}(\varphi)=\sum_{(I^+,I^-)\in \mathfrak{I}}proj_{cusp} (\varphi[I^+,I^-]).$$}

Preuve. Pour tout $f\in {\cal R}$ et tout   $\boldsymbol{\gamma}\in \boldsymbol{\Gamma}$, on note   $f_{\boldsymbol{\gamma}}$ la composante de $f$ dans  ${\cal R}(\boldsymbol{\gamma})$.  Posons $f=proj_{cusp}\circ res_{{\bf m}}(\varphi)$. Fixons $\boldsymbol{\gamma}_{0}=(r',r'',N)\in \boldsymbol{\Gamma}_{n_{0}}$, posons $\boldsymbol{\delta}_{0}=(r',(-1)^{r'}r'',N)$. Nous allons calculer    $f_{\boldsymbol{\delta}_{0}}$.  Puisque $f$ est cuspidale, il suffit de calculer $f_{\boldsymbol{\delta}_{0}}(w_{\emptyset,\beta_{1},\beta_{2}})$ pour tout couple $(\beta_{1},\beta_{2})\in {\cal P}_{2}(N)$. Posons $\boldsymbol{\delta}=(r',(-1)^{r'}r'',N+S({\bf m}))\in \boldsymbol{\Gamma}$ et $\boldsymbol{\gamma}=(r',r'',N+S({\bf m})$. Par d\'efinition de $res_{{\bf m}}$, on a $f_{\boldsymbol{\delta}_{0}}(w_{\emptyset,\beta_{1},\beta_{2}})=\varphi_{\boldsymbol{\delta}}(w_{{\bf m},\beta_{1},\beta_{2}})$. De $(r',r'')$ se d\'eduit   un quadruplet ${\bf a}$ comme en 1.10 et la formule (1) de ce paragraphe nous dit que
$$\varphi_{\boldsymbol{\delta}}(w_{{\bf m},\beta_{1},\beta_{2}})=\sum_{{\bf I}\in {\cal I}}(-1)^{a^+_{1}\vert J_{1}^+\vert +a^-_{1}\vert J_{1}^-\vert+a^+_{2}\vert J_{2}^+\vert+a^-_{2}\vert J_{2}^-\vert}$$
$$  j(\phi^+\otimes \phi^-)_{\boldsymbol{\gamma}}(w_{{\bf m}^+({\bf I}),\beta^+_{1}({\bf I})\cup \beta^+_{2}({\bf I})}\times 
 w_{{\bf m}^-({\bf I}),\beta^-_{1}({\bf I})\cup \beta^-_{2}({\bf I})}).$$
De $(r',r'')$ se d\'eduit aussi un couple d'entiers $(k^+,k^-)$ comme en 1.11. Par d\'efinition de $j$, on a $ j(\phi^+\otimes \phi^-)_{\boldsymbol{\gamma}}=0$  sauf si 

(1) $k^+(k^++1)\leq 2n^+$ et $k^-(k^-+1)\leq 2n^-$. 

On en d\'eduit d'abord
que  $f_{\boldsymbol{\delta}_{0}}(w_{\emptyset,\beta_{1},\beta_{2}})=0$ si (1) n'est pas v\'erifi\'ee.

Supposons (1) v\'erifi\'ee et d\'efinissons $N^+$ et $N^-$ par $k^+(k^++1)+2N^+=2n^+$, $k^-(k^-+1)+2N^-=2n^-$.  Notons $\phi^+_{k_{+}}$ la composante de $\phi^+$ dans le composant de ${\cal S}_{n^+}$ index\'e par $(k^+,N^+)$. On d\'efinit de m\^eme $\phi^-_{k^-}$. Toujours par d\'efinition de $j$, un terme 
$$ j(\phi^+\otimes \phi^-)_{\boldsymbol{\gamma}}(w_{{\bf m}^+({\bf I}),\beta^+_{1}({\bf I})\cup \beta^+_{2}({\bf I})}\times 
 w_{{\bf m}^-({\bf I}),\beta^-_{1}({\bf I})\cup \beta^-_{2}({\bf I})})$$
  est nul sauf si 
  
  (2) $S({\bf m}^{\zeta}({\bf I}))+S(\beta^{\zeta}_{1}({\bf I}))+S(\beta^{\zeta}_{2}({\bf I}))=N^{\zeta}$ pour $\zeta=\pm$. 
  
  Si ces conditions sont v\'erifi\'ees, il vaut
  $$\phi^+_{k^+}(w_{{\bf m}^+({\bf I}),\beta^+_{1}({\bf I})\cup \beta^+_{2}({\bf I})})\phi^-_{k^-}(w_{{\bf m}^-({\bf I}),\beta^-_{1}({\bf I})\cup \beta^-_{2}({\bf I})}).$$
Remarquons que, pour ${\bf I}=(I^+,I^-,J_{1}^+,J_{1}^-,J_{2}^+,J_{2}^-)$, les partitions ${\bf m}^+({\bf I})$ et ${\bf m}^-({\bf I})$ co\^{\i}ncident avec les partitions ${\bf m}[I^+]$ et ${\bf m}[I^-]$ d\'efinies plus haut.  Les ${\bf I} $ qui v\'erifient (2) sont les sextuplets  ${\bf I}=(I^+,I^-,J_{1}^+,J_{1}^-,J_{2}^+,J_{2}^-)$ tels que
 $$(3) \qquad I^+\sqcup I^-=\{1,...,t\},\,\, J_{1}^+\sqcup J_{1}^-=\{1,...,l(\beta_{1})\},\,\, J_{2}^+\sqcup J_{2}^-=\{1,...,l(\beta_{2})\},$$
 et
 $$(4) \qquad k^{\zeta}(k^{\zeta}+1)+2S({\bf m} (I^{\zeta}))+2S(\beta^{\zeta}_{1}({\bf I}))+2S(\beta^{\zeta}_{2}({\bf I}))=2n^{\zeta}$$
 pour $\zeta=\pm$. Cette derni\`ere condition implique
 
  (5) $k^+(k^++1)\leq 2n^+-2S({\bf m}[I^+])$, $k^-(k^-+1)\leq 2n^--2S({\bf m}[I^-])$.
  
  On note $\mathfrak{I}_{0}$ l'ensemble des $(I^+,I^-)$ v\'erifiant la premi\`ere condition de (3) ainsi que (5). Notons que c'est un sous-ensemble de $\mathfrak{I}$. Pour $(I^+,I^-)\in \mathfrak{I}_{0}$, notons ${\cal J}[I^+,I^-]$ l'ensemble des quadruplets $(J_{1}^+,J_{1}^-,J_{2}^+,J_{2}^-)$ v\'erifiant les deux derni\`eres conditions de (3) ainsi que (4). 
   On remplace les notations $\beta^+_{1}({\bf I})$ etc... par $\beta^+_{1}({\bf J})$ etc... On obtient
 $$(6)\qquad f_{\boldsymbol{\delta}_{0}}(w_{\emptyset,\beta_{1},\beta_{2}})=\sum_{(I^+,I^-)\in \mathfrak{I}_{0}}\sum_{{\bf J}\in {\cal J}[I^+,I^-]}(-1)^{a^+_{1}\vert J_{1}^+\vert +a^-_{1}\vert J_{1}^-\vert+a^+_{2}\vert J_{2}^+\vert+a^-_{2}\vert J_{2}^-\vert} 
 \phi^+_{k^+}(w_{{\bf m}[I^+],\beta^+_{1}({\bf J})\cup \beta^+_{2}({\bf J})})$$
 $$\phi^-_{k^-}(w_{{\bf m}[I^-],\beta^-_{1}({\bf J})\cup \beta^-_{2}({\bf J})}).$$
 Remarquons que, si (1) n'est pas v\'erifi\'ee, $\mathfrak{I}_{0}$ est vide. L'\'egalit\'e  (6) est donc vraie que cette condition (1) soit v\'erifi\'ee ou pas.

 Fixons $(I^+,I^-)\in \mathfrak{I}$
  et posons $f[I^+,I^-]=proj_{cusp}(\varphi[I^+,I^-])$. On calcule de la m\^eme fa\c{c}on   $f[I^+,I^-]_{\boldsymbol{\delta}_{0}}(w_{\emptyset,\beta_{1},\beta_{2}})$. L'application $res_{{\bf m}}$ dispara\^{\i}t. On obtient que $f[I^+,I^-]_{\boldsymbol{\delta}_{0}}(w_{\emptyset,\beta_{1},\beta_{2}})=0$ sauf si l'analogue de (1) est v\'erifi\'ee. Mais cette analogue est pr\'ecis\'ement la condition (5), autrement dit la condition $(I^+,I^-)\in \mathfrak{I}_{0}$. 
 Si $(I^+,I^-)\in \mathfrak{I}_{0}$,  l'ensemble ${\cal I}$ est directement remplac\'e par ${\cal J}[I^+,I^-]$. Le quadruplet ${\bf a}$ est inchang\'e. On obtient
 $$(7)\qquad  f[I^+,I^-]_{\boldsymbol{\delta}_{0}}(w_{\emptyset,\beta_{1},\beta_{2}})=\sum_{{\bf J}\in {\cal J}[I^+,I^-]}(-1)^{a^+_{1}\vert J_{1}^+\vert +a^-_{1}\vert J_{1}^-\vert+a^+_{2}\vert J_{2}^+\vert+a^-_{2}\vert J_{2}^-\vert}$$
 $$ 
 \phi^+[{\bf m}(I^+)]_{k^+}(w_{ \emptyset,\beta^+_{1}({\bf J})\cup \beta^+_{2}({\bf J})})\phi^-[{\bf m}(I^-)]_{k^-}(w_{ \emptyset,\beta^-_{1}({\bf J})\cup \beta^-_{2}({\bf J})}).$$
 
  En comparant les formules (4) et (6), on voit que, pour achever la preuve du lemme, il suffit de fixer $(I^+,I^-)\in \mathfrak{I}_{0}$, ${\bf J}\in {\cal J}[I^+,I^-]$  et de d\'emontrer  que
  le terme  
 $$ \phi^+_{k^+}(w_{{\bf m}[I^+],\beta^+_{1}({\bf J})\cup \beta^+_{2}({\bf J})})\phi^-_{k^-}(w_{{\bf m}[I^-],\beta^-_{1}({\bf J})\cup \beta^-_{2}({\bf J})})$$
 de l'expression (6) est \'egal au terme 
 $$ 
 \phi^+[{\bf m}(I^+)]_{k^+}(w_{ \emptyset,\beta^+_{1}({\bf J})\cup \beta^+_{2}({\bf J})})\phi^-[{\bf m}(I^-)]_{k^-}(w_{ \emptyset,\beta^-_{1}({\bf J})\cup \beta^-_{2}({\bf J})})$$
 de l'expression (7). Cela va r\'esulter de deux assertions parall\`eles, l'une pour les termes affect\'es d'un  exposant $+$, l'autre pour les termes affect\'es d'un exposant $-$, dont nous n'\'enoncerons que la premi\`ere.  Soit ${\bf m}^+\in {\cal P}(\leq N^+)$ et soit $\beta\in {\cal P}(N^+-S({\bf m}^+))$. Alors on a l'\'egalit\'e
 $$(8) \qquad \phi^+_{k^+}(w_{{\bf m}^+,\beta})=\phi^+[{\bf m}^+]_{k^+}(w_{\emptyset,\beta}).$$
 
Du couple $(k^+,0)$ se d\'eduit un couple $(r^{'+},r^{''+})$. On pose   $\boldsymbol{\delta}_{0}^+=(r^{'+},(-1)^{r^{'+}}r^{''+},N^+-S({\bf m}^+))\in \boldsymbol{\Gamma}_{n^+-S({\bf m}^+)}$. Posons aussi $f^+=proj_{cusp}\circ res_{{\bf m}^+}(\varphi^+)$. Le calcul du  terme $f^+_{\boldsymbol{\delta}_{0}^+}(w_{\emptyset,\beta,\emptyset})$  est un cas particulier du calcul ci-dessus. Il se simplifie grandement car les composantes affect\'ees d'un exposant $-$ disparaissent, ainsi que le $\beta_{2}$. L'analogue de l'ensemble ${\bf I}$ est r\'eduit \`a l'\'el\'ement $(I^+,I^-,J_{1}^+,J_{1}^-,J_{2}^+,J_{2}^-)$ tel que $I^+=\{1,...,l({\bf m}^+)\}$, $J_{1}^+=\{1,...,l(\beta)\}$, $I^-=J_{1}^-=J_{2}^+=J_{2}^-=\emptyset$. D'apr\`es la recette de 1.11, on a $r^{''+}\geq0$ donc $a_{1}^+=0$ d'apr\`es la d\'efinition de 1.10.  Comme analogue de la relation (6), on obtient simplement
$$f^+_{\boldsymbol{\delta}_{0}^+}(w_{\emptyset,\beta,\emptyset})=\phi^+_{k^+}(w_{{\bf m}^+,\beta}).$$
Posons $f_{0}^+=proj_{cusp}(\varphi^+[{\bf m}^+])$. On calcule de m\^eme  
$$f^+_{0,\boldsymbol{\delta}_{0}^+}(w_{\emptyset,\beta,\emptyset}) =\phi^+[{\bf m}^+]_{k^+}(w_{\emptyset,\beta}).$$
L'hypoth\`ese de l'\'enonc\'e est que $f^+=f^+_{0}$. Avec les \'egalit\'es ci-dessus, cela entra\^{\i}ne (8), ce qui ach\`eve la d\'emonstration. $\square$

\section{Endoscopie; l'involution ${\cal F}$}

\bigskip

\subsection{Endoscopie}
Pour $\pi\in {\mathbb C}[Irr_{tunip}]={\mathbb C}[Irr_{tunip,iso}]\oplus{\mathbb C}[Irr_{tunip,an}]$ et $\sharp=iso$ ou $an$, notons $\pi_{\sharp}$ la composante de $\pi$ dans ${\mathbb C}[Irr_{tunip,\sharp}]$.

On sait d\'efinir la notion de distribution stable sur $G_{iso}(F)$. Un \'el\'ement $\pi\in {\mathbb C}[Irr_{tunip,iso}]$ est dit stable si son caract\`ere-distribution (c'est-\`a-dire la forme lin\'eaire $f\mapsto trace(\pi(f))$) est stable. 

Pour travailler commod\'ement avec nos deux groupes $G_{iso}$ et $G_{an}$, il convient de modifier la d\'efinition usuelle des donn\'ees endoscopiques. Nous n'entrerons pas dans les d\'etails th\'eoriques. En pratique, une donn\'ee endoscopique elliptique est d\'etermin\'ee ici par la classe de conjugaison dans $Sp(2n;{\mathbb C})$ d'un \'el\'ement $h$ de ce groupe tel que $h^2=1$. En notant $2n_{1}$ la multiplicit\'e de $1$ parmi les valeurs propres de $h$ et $2n_{2}$ celle de $-1$, la donn\'ee endoscopique est aussi bien associ\'ee au couple $(n_{1},n_{2})\in D(n)$ (il serait plus logique de noter $n_{-1}$ au lieu de $n_{2}$ mais cela compliquerait l'\'ecriture). Le groupe endoscopique associ\'e \`a la donn\'ee est $H=G_{n_{1},iso}\times G_{n_{2},iso}$.  Pour des groupes sp\'eciaux orthogonaux impairs, il y a un choix canonique de facteurs de transfert, cf. \cite{W2} 1.10.  
  On doit faire attention au point suivant. Un tel facteur est une fonction $\Delta_{h,\sharp}$ d\'efinie presque partout sur $H(F)\times G_{\sharp}(F)$. Si on multiplie $h$ par $-1$, on remplace $H$ par le groupe isomorphe $H'=G_{n_{2},iso}\times G_{n_{1},iso}$. Notons $\iota$ l'isomorphisme de $H$ sur $H'$ qui permute les deux facteurs. Pour $h\in H(F)$ et $g\in G_{\sharp}(F)$, on a alors 
  
  - si $\sharp=iso$, $\Delta_{-h,iso}(\iota(h),g)=\Delta_{h,iso}(h,g)$;
  
  - si $\sharp=an$, $\Delta_{-h,an}(\iota(h),g)=-\Delta_{h,an}(h,g)$.
  
  En tout cas, on dispose d'un tranfert qui envoie une combinaison lin\'eaire stable de repr\'esentations temp\'er\'ees de $H(F)$ sur une combinaison lin\'eaire de repr\'esentations temp\'er\'ees de $G_{\sharp}(F)$. Notons-le simplement $transfert_{h,\sharp}$. 
  
  En particulier, pour $h=1$, on a $H=G_{iso}$. Evidemment, le transfert de $H$ vers $G_{iso}$ est l'identit\'e. Mais il y a aussi un transfert de $H=G_{iso}$ vers $G_{an}$.

Introduisons l'automorphisme ext\'erieur  $\theta:g\mapsto {^tg}^{-1}$ de $GL(2n)$, le groupe $GL(2n)\rtimes\{1,\theta\}$ et sa composante connexe $\tilde{GL}(2n)=GL(2n)\theta$.   Par endoscopie tordue, il y a aussi une notion de transfert qui envoie une combinaison lin\'eaire stable de repr\'esentations temp\'er\'ees de $G_{iso}(F)$ sur une combinaison lin\'eaire  de repr\'esentations temp\'er\'ees tordues de $\tilde{GL}(2n;F)$.

Introduisons l'ensemble $\mathfrak{Endo}_{tunip}$ form\'e des classes de conjugaison de triplets $(\lambda,s,h)$ tels que $\lambda\in {\cal P}^{symp}(2n)$, $s,h\in Z(\lambda)$, $s$ et $h$ commutent entre eux ($sh=hs$),  $s$ est compact   et $h^2=1$. La notion de conjugaison est la conjugaison de $s$ et $h$ par un m\^eme \'el\'ement de $Z(\lambda)$. Pour un tel triplet, $h$ appartient \`a $Z(\lambda,s)$, notons ${\bf h}$ son image dans ${\bf Z}(\lambda,s)$. Pour $\epsilon\in {\bf Z}(\lambda,s)^{\wedge}$, on peut \'evaluer $\epsilon$ au point ${\bf h}$. On note simplement $\epsilon(h)$ cette valeur. On pose
$$\Pi(\lambda,s,h)=\sum_{\epsilon\in {\bf Z}(\lambda,s)^{\wedge}}\pi(\lambda,s,\epsilon)\epsilon(h).$$
Par lin\'earit\'e, on obtient une application
$$\Pi:{\mathbb C}[\mathfrak{Endo}_{tunip}]\to {\mathbb C}[Irr_{tunip}].$$

Notons $\mathfrak{St}_{tunip}$ le sous-ensemble des $(\lambda,s,h)\in \mathfrak{Endo}_{tunip}$ tels que $h=1$. 
Les r\'esultats principaux de \cite{MW} sont que 
  
 pour $(\lambda,s,1)\in \mathfrak{St}_{tunip}$,  $\Pi_{iso}(\lambda,s,1)$ est stable et $-\Pi_{an}(\lambda,s,1)$ est son  transfert \`a la forme non d\'eploy\'ee $G_{an}$;
  
  tout \'el\'ement de $\pi\in{\mathbb C}[Irr_{tunip}]$  tel que $\pi_{iso}$ soit stable et que $-\pi_{an}$ soit le transfert de $\pi_{iso}$ est combinaison lin\'eaire des $ \Pi(\lambda,s,1)$ quand $(\lambda ,s,1)$ parcourt $\mathfrak{St}_{tunip}$.

  Pour tout entier $m\geq1$, notons $st_{m}$ la repr\'esentation de Steinberg de $GL(m;F)$. Pour une partition $\mu=(\mu_{1}\geq...\geq\mu_{t}>0)$, notons $st_{\mu}$ la repr\'esentation de $GL(S(\mu);F)$ qui est induite, en un sens \'evident, de $st_{\mu_{1}}\otimes...\otimes st_{\mu_{t}}$. Soit $(\lambda,s,1)\in \mathfrak{St}_{tunip}$. D\'ecomposons $\lambda$ en 
  $$\lambda^+\cup \lambda^-\cup\bigcup_{j\in J}(\lambda_{j}\cup\lambda_{j})$$
  selon les valeurs propres de $s$, cf. 1.3. Pour tout $j\in J$, notons $\chi_{j}$ l'unique caract\`ere non ramifi\'e de $F^{\times}$ tel que $\chi_{j}(\varpi)=s_{j}$. Notons $\chi^-$ l'unique caract\`ere non ramifi\'e de $F^{\times}$ tel que $\chi^-(\varpi)=-1$.   Notons $\Pi^{GL}(\lambda,s)$ l'induite de 
  $st_{\lambda^+}\otimes (\chi^-\circ det)st_{\lambda^-}\otimes \otimes_{j\in J}((\chi_{j}\circ det)st_{\lambda_{j}})\otimes (\chi_{j}^{-1}\circ det)st_{\lambda_{j}}))$.  Cette repr\'esentation se prolonge en une repr\'esentation tordue $\tilde{\Pi}^{GL}(\lambda,s)$ de $\tilde{GL}(2n;F)$ (il y a diff\'erents prolongements possibles, il faut normaliser le prolongement de fa\c{c}on ad\'equate). 
Le r\'esultat principal de \cite{W3} est que le transfert par endoscopie tordue de $\Pi_{iso}(\lambda,s,1)$ \`a $GL(2n;F)$ est $\tilde{\Pi}^{GL}(\lambda,s)$. 

{\bf Remarque.} Ceci n'est pas tout-\`a-fait exact. Dans \cite{W3}, on n'a \'enonc\'e le r\'esultat que dans le cas o\`u $s^2=1$. Mais le r\'esultat g\'en\'eral en r\'esulte par induction.
  
  \bigskip
  
  Ces propri\'et\'es caract\'erisent enti\`erement $\Pi_{iso}(\lambda,s,1)$  et $\Pi_{an}(\lambda,s,1)$.

 Soit maintenant $h\in Sp(2n;{\mathbb C})$ un \'el\'ement tel que $h^2=1$.   Il d\'etermine une donn\'ee endoscopique elliptique dont le groupe est $G_{n_{1},iso}\times G_{n_{2},iso}$ avec les notations ci-dessus. Soient $(\lambda_{1},s_{1},1)\in \mathfrak{St}_{n_{1},tunip}$ et $(\lambda_{2},s_{2},1)\in \mathfrak{St}_{n_{2},tunip}$.  Pour $j=1,2$, $\lambda_{j}$ param\`etre 
un homomorphisme $\rho_{j}:SL(2;{\mathbb C})\to Sp(2n_{j};{\mathbb C})$ 
  et $s_{j}$ est un \'el\'ement du   commutant de l'image de $\rho_{j}$.
    Le groupe $Sp(2n_{1};{\mathbb C})\times Sp(2n_{2};{\mathbb C})$ est le commutant de $h$ dans 
$Sp(2n;{\mathbb C})$ et est donc plong\'e dans ce groupe. L'homomorphisme $\rho_{1}\times \rho_{2}$ devient un homomorphisme 
$\rho$ de $SL(2;{\mathbb C})$ dans $Sp(2n;{\mathbb C})$. Evidemment, $\rho$ est param\'etr\'e par $\lambda=\lambda_{1}\cup \lambda_{2}$.  L'\'el\'ement $(s_{1},s_{2})$ devient un \'el\'ement $s$ du commutant de  l'image de $\rho$, c'est-\`a-dire de $Z(\lambda)$. On a aussi $h\in Z(\lambda)$ par construction. Les \'el\'ements $s$ et $h$ commutent. On conclut que $(\lambda,s,h)$ appartient \`a $\mathfrak{Endo}_{tunip}$.   

\ass{Th\'eor\`eme}{Sous ces hypoth\`eses, on a les \'egalit\'es

$transfert_{h,iso}(\Pi_{iso}(\lambda_{1},s_{1},1)\otimes \Pi_{iso}(\lambda_{2},s_{2},1))=\Pi_{iso}(\lambda,s,h)$;

$transfert_{h,an}(\Pi_{iso}(\lambda_{1},s_{1},1)\otimes \Pi_{iso}(\lambda_{2},s_{2},1))=-\Pi_{an}(\lambda,s,h)$.}

A l'aide d'une transformation de Fourier sur ${\bf Z}(\lambda,s)$,  les repr\'esentations $\pi( \lambda,s,\epsilon)$  pour $( \lambda,s,\epsilon)\in \mathfrak{Irr}_{tunip}$,  s'expriment comme combinaisons lin\'eaires de repr\'esentations $\Pi(\lambda,s,h)$. Alors, les propri\'et\'es de l'\'enonc\'e ci-dessus caract\'erisent enti\`erement nos repr\'esentations $\pi( \lambda,s,\epsilon)$.   Le th\'eor\`eme sera d\'emontr\'e  dans l'article suivant.

{\bf Remarque.} Pour d\'efinir un transfert endoscopique entre fonctions, il est n\'ecessaire de fixer des mesures de Haar sur tous les groupes consid\'er\'es. Mais le transfert dual entre repr\'esentations est ind\'ependant de ces choix.

\bigskip

\subsection{Le cas "quadratique-unipotent"}

 Notons $\mathfrak{Endo}_{unip-quad}$    le sous-ensemble des $(\lambda,s,h)\in \mathfrak{Endo}_{tunip}$ tels que $s^2=1$.  
 Soit $(\lambda,s,h)\in \mathfrak{Endo}_{unip-quad}$ et soit $i\in Jord(\lambda)$.  On note $s_{i}$ et $h_{i}$ les composantes de $s$ et $h$ dans la composante $Sp(mult_{\lambda}(i);{\mathbb C})$ ou $O(mult_{\lambda}(i);{\mathbb C})$  du groupe $Z(\lambda)$. Cette composante agit naturellement dans un espace ${\mathbb C}^{mult_{\lambda}(i)}$ et $s_{i}$ et $h_{i}$ sont des automorphismes de cet espace, de carr\'es $1$ et qui commutent. Ainsi, l'espace se d\'ecompose en sous-espaces propres communs pour $s_{i}$ et $h_{i}$. Pour $\zeta,\xi\in \{\pm\}\simeq \{\pm 1\}$, on note $m_{s,h}^{\zeta\xi}(i)$ la dimension du sous-espace propre o\`u $s_{i}$ agit par $\zeta$ et $h_{i}$ agit par $\xi$. D'autre part, on a vu qu'\`a $s$ est associ\'ee une d\'ecomposition $\lambda=\lambda^+\cup \lambda^-$. On calcule alors
$$(1) \qquad \Pi(\lambda,s,h)=\sum_{\epsilon^+,\epsilon^-}\pi(\lambda^+,\epsilon^+,\lambda^-,\epsilon^-)\left(\prod_{i\in Jord_{bp}(\lambda^+)}\epsilon^+(i)^{m_{s,h}^{+-}(i)}\right)\left(\prod_{i\in Jord_{bp}(\lambda^-)}\epsilon^-(i)^{m_{s,h}^{--}(i)}\right),$$
o\`u $(\epsilon^+,\epsilon^-)$ parcourt l'ensemble $\{\pm 1\}^{Jord_{bp}(\lambda^+)}\times \{\pm 1\}^{Jord_{bp}(\lambda^-)}$.

Pour $(\lambda,s,h)$ et $(\bar{\lambda},\bar{s},\bar{h})\in  \mathfrak{Endo}_{unip-quad}$, disons que ces deux triplets sont $Res$-\'equivalents si et seulement si $\lambda=\bar{\lambda}$ et $m_{s,h}^{\zeta\xi}(i)\equiv m_{\bar{s},\bar{h}}^{\zeta\xi}(i)\,\,mod\,\,2{\mathbb Z}$ pour tout $i\in Jord(\lambda)$ et tous $\zeta,\xi=\pm$. Notons $\mathfrak{K}$ le sous-espace de ${\mathbb C}[\mathfrak{Endo}_{unip-quad}]$ engendr\'e par les diff\'erences $(\lambda,s,h)-(\lambda',s',h')$ de deux \'el\'ements $Res$-\'equivalents. 

\ass{Lemme}{L'application $Res\circ\Pi:{\mathbb C}[\mathfrak{Endo}_{unip-quad}]\to {\cal R}^{par,glob}$ est surjective. Son noyau est l'espace $\mathfrak{K}$.}

{\bf Remarque.} En vertu du lemme 1.7, l'\'enonc\'e reste vrai si l'on y remplace l'application $Res\circ \Pi$ par $Res\circ D\circ\Pi$.

\bigskip

Preuve.  
Pour $(\lambda,s,h)$ et $(\bar{\lambda},\bar{s},\bar{h})\in  \mathfrak{Endo}_{unip-quad}$, disons que ces deux triplets sont li\'es si, quitte \`a les permuter, les conditions suivantes sont v\'erifi\'ees:

$\bar{\lambda}=\lambda$;

il existe $i\in Jord(\lambda)$ et il existe $\zeta,\xi,\bar{\zeta},\bar{\xi}=\pm$, avec $(\zeta,\xi)\not=(\bar{\zeta},\bar{\xi})$, de sorte que $m_{\bar{s},\bar{h}}^{\zeta,\xi}(i)=m_{s,h}^{\zeta,\xi}(i)-2$, $m_{\bar{s},\bar{h}}^{\bar{\zeta}\bar{\xi}}(i)=m_{s,h}^{\bar{\zeta}\bar{\xi}}(i)+2$ et $m_{\bar{s},\bar{h}}^{\zeta'\xi'}(i')=m_{s,h}^{\zeta'\xi'}(i')$ pour tout triplet $(i',\zeta',\xi')$ diff\'erent de $(i,\zeta,\xi)$ et $(i,\bar{\zeta},\bar{\xi})$. Consid\'erons deux tels triplets li\'es (et abandonnons la notation $\bar{\lambda}$ inutile puisque $\bar{\lambda}=\lambda$). Montrons que
$$(2) \qquad Res\circ\Pi(\lambda,\bar{s},\bar{h})=Res\circ\Pi(\lambda,s,h).$$
On utilise les donn\'ees $i$, $\zeta$, $\xi$, $\bar{\zeta}$, $\bar{\xi}$ fournies par la d\'efinition ci-dessus. Supposons d'abord $\bar{\zeta}=\zeta$. Les valeurs propres des termes $s$ et $\bar{s}$ sont alors les m\^emes et on peut supposer $\bar{s}=s$. La d\'ecomposition $\lambda=\lambda^+\cup \lambda^-$ est la m\^eme pour nos deux triplets. On utilise la formule (1) pour chacun de nos triplets. Les seuls termes qui changent d'un triplet \`a l'autre sont les  exposants des termes $\epsilon^+(i)$ et $\epsilon^-(i)$. Mais ces exposants $m_{s,h}^{+-}(i)$ etc... n'interviennent que par leur parit\'e et, par d\'efinition, celle-ci ne change pas quand on remplace $h$ par $\bar{h}$. L'\'egalit\'e (2) en r\'esulte.  Supposons maintenant que $\bar{\zeta}\not=\zeta$. Pour simplifier la notation, supposons $\zeta=+$ et $\bar{\zeta}=-$. Notons $\lambda=\lambda^+\cup \lambda^-$ la d\'ecomposition associ\'ee \`a $s$ et $\lambda=\bar{\lambda}^+\cup \bar{\lambda}^-$ celle associ\'ee \`a $\bar{s}$. On a $\lambda^+=\bar{\lambda}^+\cup\{i,i\}$ et $\bar{\lambda}^-=\lambda^-\cup\{i,i\}$. Il y a des inclusions naturelles $Jord_{bp}(\bar{\lambda}^+)\subset Jord_{bp}(\lambda^+)$ et $Jord_{bp}(\lambda^-)\subset Jord_{bp}(\bar{\lambda}^-)$. Pour $(\bar{\epsilon}^+,\epsilon^-)\in \{\pm 1\}^{Jord_{bp}(\bar{\lambda}^+)}\times \{\pm 1\}^{Jord_{bp}(\lambda^-)}$, posons
$$P_{s,h}^{i}(\bar{\epsilon}^+,\epsilon^-)=\left(\prod_{i'\in Jord_{bp}(\bar{\lambda}^+),i'\not=i}\bar{\epsilon}^+(i')^{m_{s,h}^{+-}(i')}\right)\left(\prod_{i'\in Jord_{bp}(\lambda^-),i'\not=i}\epsilon^-(i')^{m_{s,h}^{--}(i')}\right),$$
$$R_{s,h}(\bar{\epsilon}^+,\epsilon^-)=\sum_{\epsilon^+}\pi(\lambda^+,\epsilon^+,\lambda^-,\epsilon^-)\epsilon^+(i)^{m_{s,h}^{+-}(i)}\epsilon^-(i)^{m_{s,h}^{--}(i)},$$
o\`u $\epsilon^+$ parcourt les \'el\'ements de $Jord_{bp}(\lambda^+)$ qui prolongent $\bar{\epsilon}^+$ (il y en a $1$ ou $2$) et, par convention,  $\epsilon^+(i)=1$ si $i\not\in Jord_{bp}(\lambda^+)$ et $\epsilon^-(i)=1$ si  $i\not\in Jord_{bp}(\lambda^-)$. On d\'efinit $P_{\bar{s},\bar{h}}^{i}(\bar{\epsilon}^+,\epsilon^-)$ de la m\^eme fa\c{c}on que ci-dessus et on pose
$$R_{\bar{s},\bar{h}}(\bar{\epsilon}^+,\epsilon^-)=\sum_{\bar{\epsilon}^-}\pi(\bar{\lambda}^+,\bar{\epsilon}^+,\bar{\lambda}^-,\bar{\epsilon}^-)\bar{\epsilon}^+(i)^{m_{\bar{s},\bar{h}}^{+-}(i)}\bar{\epsilon}^-(i)^{m_{\bar{s},\bar{h}}^{--}(i)},$$
o\`u $\bar{\epsilon}^-$ parcourt les \'el\'ements de $Jord_{bp}(\bar{\lambda}^-)$ qui prolongent $\epsilon^-$ (il y en a $1$ ou $2$) et avec la m\^eme convention que ci-dessus. 
On peut r\'ecrire (1) sous la forme
$$\Pi(\lambda,s,h)=\sum_{\bar{\epsilon}^+,\epsilon^-}P_{s,h}^{i}(\bar{\epsilon}^+,\epsilon^-)R_{s,h}(\bar{\epsilon}^+,\epsilon^-),$$
et on a la formule similaire
$$\Pi(\lambda,\bar{s},\bar{h})=\sum_{\bar{\epsilon}^+,\epsilon^-}P_{\bar{s},\bar{h}}^{i}(\bar{\epsilon}^+,\epsilon^-)R_{\bar{s},\bar{h}}(\bar{\epsilon}^+,\epsilon^-).$$
Il nous suffit de fixer $(\bar{\epsilon}^+,\epsilon^-)$ et de d\'emontrer l'\'egalit\'e
$$(3) \qquad P_{s,h}^{i}(\bar{\epsilon}^+,\epsilon^-)Res(R_{s,h}(\bar{\epsilon}^+,\epsilon^-))=P_{\bar{s},\bar{h}}^{i}(\bar{\epsilon}^+,\epsilon^-)Res(R_{\bar{s},\bar{h}}(\bar{\epsilon}^+,\epsilon^-)).$$
Il est clair que $P_{s,h}^{i}(\bar{\epsilon}^+,\epsilon^-)=P_{\bar{s},\bar{h}}^{i}(\bar{\epsilon}^+,\epsilon^-)$: pour les $i'\in Jord_{bp}(\lambda)$ tels que $i'\not=i$, les multiplicit\'es ne changent pas quand on remplace $s,h$ par $\bar{s},\bar{h}$. Dans la d\'efinition de $R_{s,h}(\bar{\epsilon}^+,\epsilon^-)$, on peut remplacer le terme $\epsilon^+(i)^{m_{s,h}^{+-}(i)}$ par $\bar{\epsilon}^+(i)^{m_{\bar{s},\bar{h}}^{+-}(i)}$. En effet, si $i\in Jord_{bp}(\bar{\lambda}^+)$, il n'y a qu'un seul prolongement $\epsilon^+$ de $\bar{\epsilon}^+$, pour lequel $\epsilon^+(i)=\bar{\epsilon}^+(i)$, et les exposants $m_{s,h}^{+-}(i)$ et $m_{\bar{s},\bar{h}}^{+-}(i)$ sont de m\^eme parit\'e.  Si $i\not\in Jord_{bp}(\bar{\lambda}^+)$, on a $\bar{\epsilon}^+(i)=1$ par convention. Si $i$ est impair, on a aussi $\epsilon^+(i)=1$ par convention. Si $i$ est pair, l'hypoth\`ese $i\not\in Jord_{bp}(\bar{\lambda}^+)$ entra\^{\i}ne $m_{\bar{s},\bar{h}}^{+-}(i)=0$, donc $m_{s,h}^{+-}(i)$ est pair et $\epsilon^+(i)^{m_{s,h}^{+-}(i)}=1$ pour chacun des deux prolongements $\epsilon^+$. On obtient
$$R_{s,h}(\bar{\epsilon}^+,\epsilon^-)=\bar{\epsilon}^+(i)^{m_{\bar{s},\bar{h}}^{+-}(i)}\epsilon^-(i)^{m_{s,h}^{--}(i)}\sum_{\epsilon^+}\pi(\lambda^+,\epsilon^+,\lambda^-,\epsilon^-).$$
De m\^eme
$$R_{\bar{s},\bar{h}}(\bar{\epsilon}^+,\epsilon^-)=\bar{\epsilon}^+(i)^{m_{\bar{s},\bar{h}}^{+-}(i)}\epsilon^-(i)^{m_{s,h}^{--}(i)}\sum_{\bar{\epsilon}^-}\pi(\bar{\lambda}^+,\bar{\epsilon}^+,\bar{\lambda}^-,\bar{\epsilon}^-).$$
Les premiers facteurs sont les m\^emes et le lemme 1.6 affirme pr\'ecis\'ement que
$$Res(\sum_{\epsilon^+}\pi(\lambda^+,\epsilon^+,\lambda^-,\epsilon^-))=Res(\sum_{\bar{\epsilon}^-}\pi(\bar{\lambda}^+,\bar{\epsilon}^+,\bar{\lambda}^-,\bar{\epsilon}^-)).$$
Cela d\'emontre (3), d'o\`u (2). 

Pour deux \'el\'ements $(\lambda,s,h)$ et $(\bar{\lambda},\bar{s},\bar{h})$ de $\mathfrak{Endo}_{unip-quad}$ qui sont $Res$-\'equivalents, on voit qu'il existe une suite $(\lambda,s,h)=(\lambda_{1},s_{1},h_{1})$, $(\lambda_{2},s_{2},h_{2})$,...,$(\lambda_{k},s_{k},h_{k})=(\bar{\lambda},\bar{s},\bar{h})$ d'\'el\'ements de $\mathfrak{Endo}_{unip-quad}$ de sorte que, pour $j=1,...,k-1$, $(\lambda_{j},s_{j},h_{j})$ et $(\lambda_{j+1},s_{j+1},h_{j+1})$ soient li\'es. En appliquant (2) \`a chacun de ces couples,  on obtient
$$Res\circ\Pi(\lambda,s,h)=Res\circ\Pi(\bar{\lambda},\bar{s},\bar{h}).$$
Cela d\'emontre que $\mathfrak{K}$ est contenu dans le noyau de $Res\circ \Pi$.

Notons $\mathfrak{Endo}_{unip-quad}^{red}$ le sous-ensemble des $(\lambda,s,h)\in \mathfrak{Endo}_{unip-quad}$ tels que, pour tout $i\in Jord(\lambda)$,   on ait $m_{s,h}^{+-}(i)\leq1$, $m_{s,h}^{-+}(i)\leq 1$, $m_{s,h}^{--}(i)\leq 1$. 
 Il est clair que tout \'el\'ement de $\mathfrak{Endo}_{unip-quad}$ est $Res$-\'equivalent \`a un unique \'el\'ement de $\mathfrak{Endo}^{red}_{unip-quad}$. En cons\'equence, on a l'\'egalit\'e
$${\mathbb C}[\mathfrak{Endo}_{unip-quad}]={\mathbb C}[\mathfrak{Endo}_{unip-quad}^{red}]\oplus \mathfrak{K}.$$
Pour achever la preuve du lemme, il suffit de prouver que

(4)  la restriction de $Res\circ\Pi$ \`a ${\mathbb C}[\mathfrak{Endo}_{unip-quad}^{red}]$ est un isomorphisme de cet espace sur ${\cal R}^{par,glob}$. 

 Notons $\mathfrak{N}$ l'ensemble des quadruplets $(\lambda,\tilde{\gamma},m,\tilde{m})$ tels que
 
 $\lambda\in {\cal P}^{symp}(2n)$;

$\tilde{\gamma}$, $m$ et $\tilde{m}$ sont des \'el\'ements de $\{\pm 1\}^{Jord_{bp}(\lambda)}$;

pour $i\in Jord_{bp}(\lambda)$ tel que $mult_{\lambda}(i)=1$, on a $\tilde{\gamma}(i)=1$;

pour $i\in Jord_{bp}(\lambda)$ tel que $mult_{\lambda}(i)=2$, on a $(\tilde{\gamma}(i),m(i),\tilde{m}(i))\not=(-1,1,1)$.

Dans \cite{MW} 6.7, on a associ\'e \`a tout tel quadruplet un \'el\'ement $rea(\psi_{\lambda,\tilde{\gamma},m,\tilde{m}})\in {\mathbb C}[Irr_{unip-quad}]$ (\`a ceci pr\`es que, dans \cite{MW}, la partition $\lambda$ est remplac\'ee par l'orbite unipotente symplectique ${\cal O}$ qu'elle param\`etre). D'apr\`es \cite{MW} 6.9, quand $ (\lambda,\tilde{\gamma},m,\tilde{m})$ d\'ecrit $\mathfrak{N}$, les \'el\'ements $rea(\psi_{\lambda,\tilde{\gamma},m,\tilde{m}})$ sont lin\'eairement ind\'ependants. Notons ${\mathbb C}[Irr_{unip-quad}]^0$ le sous-espace de ${\mathbb C}[Irr_{unip-quad}]$ engendr\'e par ces \'el\'ements. La proposition 6.21 de \cite{MW} dit que la restriction de $Res\circ D$ \`a ce sous-espace est un isomorphisme de celui-ci sur ${\cal R}^{par,glob}$.  D'apr\`es la remarque suivant l'\'enonc\'e ci-dessus, cela revient \`a dire que la restriction de $Res$ \`a ce sous-espace est aussi un  isomorphisme de celui-ci sur ${\cal R}^{par,glob}$. Pour prouver (4), il suffit donc de d\'efinir  une bijection $\iota:\mathfrak{Endo}^{red}_{unip-quad}\to \mathfrak{N}$ de telle sorte que, pour tout $ (\lambda,s,h)\in \mathfrak{Endo}^{red}_{unip-quad}$, on ait l'\'egalit\'e
$$(5) \qquad Res\circ \Pi(\lambda,s,h)=Res(rea(\psi_{\lambda',\tilde{\gamma},m,\tilde{m}})),$$
o\`u $(\lambda',\tilde{\gamma},m,\tilde{m})=\iota( \lambda,s,h)$. 
Construisons cette bijection. Soit $ (\lambda,s,h)\in \mathfrak{Endo}^{red}_{unip-quad}$. On d\'efinit un triplet $(\tilde{\gamma},m,\tilde{m})$ de la fa\c{c}on suivante. Soit $i\in Jord_{bp}(\lambda)$. On pose
$$m(i)=(-1)^{m_{s,h}^{-+}(i)+m_{s,h}^{--}(i)},\, \tilde{m}(i)=(-1)^{m_{s,h}^{+-}(i)+m_{s,h}^{--}(i)}.$$
 Si $mult_{\lambda}(i)=1$, on pose
$\tilde{\gamma}(i)=1$.  
Si $mult_{\lambda}(i)\geq2$, on pose
$\tilde{\gamma}(i)=(-1)^{m_{s,h}^{--}(i)}$. On v\'erifie que le quadruplet  $(\lambda,\tilde{\gamma},m,\tilde{m})$ appartient \`a $\mathfrak{N}$. On pose $\iota(\lambda,s,h)=(\lambda,\tilde{\gamma},m,\tilde{m})$. Il est facile de v\'erifier que l'application $\iota$ ainsi d\'efinie est une bijection de $\mathfrak{Endo}^{red}_{unip-quad}$ sur $ \mathfrak{N}$. On doit prouver la relation (5). Utilisons les notations de cette relation (on a $\lambda'=\lambda$ par d\'efinition de $\iota$). Le lemme  6.7 de \cite{MW} calcule le terme $rea(\psi_{\lambda,\tilde{\gamma},m,\tilde{m}})$. Avec nos pr\'esentes notations, ce lemme dit que
$$(6) \qquad rea(\psi_{\lambda,\tilde{\gamma},m,\tilde{m}})=2^{-\vert {\cal E}\vert }\sum_{(s',h')\in {\cal E}}\Pi(\lambda, s',h'),$$
o\`u ${\cal E}$ est un certain ensemble de paires $(s',h')$ telles que $(\lambda,s',h')$ appartienne \`a $\mathfrak{Endo}_{unip-quad}$.    L'\'enonc\'e de ce lemme d\'ecrit l'ensemble ${\cal E}$. En utilisant cette description  et la d\'efinition de $\iota$ donn\'ee ci-dessus, on s'aper\c{c}oit que l'ensemble ${\cal E}$ n'est pas vide et que, pour tout $ (s',h')\in {\cal E}$, l'\'el\'ement $(\lambda,s',h')$ est  $Res$-\'equivalent \`a $(\lambda, s,h)$. Puisque deux \'el\'ements $Res$-\'equivalents ont m\^eme image par $Res\circ\Pi$, l'\'egalit\'e (5) r\'esulte de (6). Cela ach\`eve la d\'emonstration. $\square$

 \bigskip

\subsection{D\'efinition de l'involution}
On d\'efinit une involution ${\cal F}$ de l'ensemble $\mathfrak{Endo}_{unip-quad}$ par ${\cal F}(\lambda,s,h)=(\lambda,h,s)$ pour tout $(\lambda,s,h)\in \mathfrak{Endo}_{unip-quad}$.  Il r\'esulte de la d\'efinition de la $Res$-\'equivalence que, si deux \'el\'ements de $\mathfrak{Endo}_{unip-quad}$ sont $Res$-\'equivalents, leurs images par ${\cal F}$ sont encore $Res$-\'equivalents. 
On prolonge  l'involution  ${\cal F}$ en une involution lin\'eaire de ${\mathbb C}[\mathfrak{Endo}_{unip-quad}]$, encore not\'ee ${\cal F}$. La propri\'et\'e pr\'ec\'edente entra\^{\i}ne que ${\cal F}$ conserve le sous-espace $\mathfrak{K}$. D'apr\`es le lemme 2.2 et la remarque qui le suit, ${\cal F}$ se descend via l'application $Res\circ D\circ\Pi$ en une involution de ${\cal R}^{par,glob}$, notons-la $\mathfrak{F}^{par}$. On a donc par d\'efinition l'\'egalit\'e $\mathfrak{F}^{par}\circ Res\circ D\circ\Pi=Res\circ D\circ\Pi\circ{\cal F}$. On a aussi d\'efini au paragraphe 1.9 une involution ${\cal F}^{par}$ de l'espace ${\cal R}^{par,glob}$. En admettant le th\'eor\`eme 2.1, nous d\'emontrerons  en 2.7  que les deux involutions ${\cal F}^{par}$ et $\mathfrak{F}^{par}$ sont \'egales. Pour le moment, contentons-nous du lemme inconditionnel suivant.

\ass{Lemme}{On a l'\'egalit\'e $proj_{cusp}\circ \mathfrak{F}^{par}={\cal F}^{par}\circ proj_{cusp}$.}

Preuve. En \cite{MW} section 6, on a d\'efini une involution de ${\cal R}^{par,glob}$ similaire \`a $\mathfrak{F}^{par}$, par une construction l\'eg\`erement diff\'erente. Notons-la $\mathfrak{F}^{par}_{MW}$. Montrons que

(1) on a l'\'egalit\'e $\mathfrak{F}^{par}=\mathfrak{F}^{par}_{MW}$.

D'apr\`es la remarque suivant l'\'enonc\'e du lemme 2.2, on peut glisser des involutions $D$ dans les assertions (4) et (5) de la preuve de ce lemme. Gr\^ace \`a cette assertion (4), il suffit de prouver l'\'egalit\'e
$$(2) \qquad \mathfrak{F}^{par}\circ Res\circ D\circ\Pi(\lambda,s,h)=\mathfrak{F}_{MW}^{par}\circ Res\circ D\circ\Pi(\lambda,s,h)$$
pour tout $(\lambda,s,h)\in \mathfrak{Endo}_{unip-quad}^{red}$. Posons $\iota(\lambda,s,h)=(\lambda,\tilde{\gamma},m,\tilde{m})$. D'apr\`es 2.2(5), on a
$$\mathfrak{F}_{MW}^{par}\circ Res\circ D\circ\Pi(\lambda,s,h)=\mathfrak{F}_{MW}^{par}\circ Res\circ D(rea(\psi_{\lambda,\tilde{\gamma},m,\tilde{m}})).$$
Par d\'efinition de $\mathfrak{F}^{par}_{MW}$, cf. \cite{MW} 6.4, le membre de droite ci-dessus est \'egal \`a
$$Res\circ D(rea(\psi_{\lambda,\tilde{\gamma},\tilde{m},m})).$$
Mais il r\'esulte de la d\'efinition de la bijection $\iota$ (cf. 2.2) que $(\lambda,\tilde{\gamma},\tilde{m},m)=\iota(\lambda,h,s)=\iota\circ{\cal F}(\lambda,s,h)$ (remarquons que $(\lambda,h,s)$ appartient encore \`a $\mathfrak{Endo}_{unip-quad}^{red}$). En utilisant de nouveau 2.2(5), on  obtient  l'\'egalit\'e
$$\mathfrak{F}_{MW}^{par}\circ Res\circ D(rea(\psi_{\lambda,\tilde{\gamma},m,\tilde{m}}))=Res\circ D\circ \Pi\circ{\cal F}(\lambda,s,h).$$
Le membre de droite ci-dessus est \'egal au membre de gauche de (2) par d\'efinition de $\mathfrak{F}^{par}$. Cela d\'emontre (2), d'o\`u (1). 

Le th\'eor\`eme 6.26 de \cite{MW} affirme que l'\'egalit\'e de l'\'enonc\'e est vraie si l'on y remplace $\mathfrak{F}^{par}$ par $\mathfrak{F}^{par}_{MW}$. L'\'egalit\'e de ces deux  involutions entra\^{\i}ne donc l'\'enonc\'e. $\square$

\bigskip

\subsection{Le cas elliptique}
Notons $\mathfrak{Endo}_{unip,disc}$ le sous-ensemble des $(\lambda,s,h)\in \mathfrak{Endo}_{unip-quad}$ tels que le commutant commun de $s$ et $h$ dans $Z(\lambda)$ (c'est-\`a-dire le groupe des $x\in Z(\lambda)$ tels que $sx=xs$ et $hx=xh$) soit fini. Ce commutant commun se calcule facilement. C'est le produit des groupes $Sp(m^{\zeta\xi}_{s,h}(i); {\mathbb C})$ pour les $i\in Jord(\lambda)$ impairs et les $\zeta,\xi=\pm$ et des groupes $O(m^{\zeta\xi}_{s,h}(i);{\mathbb C})$ pour les $i\in Jord_{bp}(\lambda)$ et les $\zeta,\xi=\pm$. Que ce groupe soit fini est \'equivalent \`a ce que $m^{\zeta\xi}_{s,h}(i)\leq1$ pour tous $i$, $\zeta$, $\xi$.  Dans le cas o\`u $i$ est impair,  $m_{s,h}^{\zeta\xi}(i)$ est forc\'ement pair, la condition \'equivaut donc  \`a $m_{s,h}^{\zeta\xi}(i)=0$. Puisque $mult_{\lambda}(i)$ est la somme des $m_{s,h}^{\zeta\xi}(i)$ sur les $\zeta,\xi$, cela entra\^{\i}ne que $\lambda$ ne contient que des termes pairs. 

On v\'erifie facilement que l'image par l'application $\Pi$ de l'espace ${\mathbb C}[\mathfrak{Endo}_{unip,disc}]$ est \'egal \`a ${\mathbb C}[Ell_{unip}]$. On a une suite d'applications lin\'eaires 
$${\mathbb C}[\mathfrak{Endo}_{unip,disc}]\stackrel{\Pi}{\to}{\mathbb C}[Ell_{unip}]\stackrel{Res}{\to}{\cal R}^{par,glob}\stackrel{proj_{cusp}}{\to}{\cal R}^{par}_{cusp}.$$
L'ensemble $\mathfrak{Endo}_{unip,disc}$ est inclus dans l'ensemble $\mathfrak{Endo}_{unip-quad}^{red}$ introduit dans la d\'emonstration pr\'ec\'edente. Comme on l'a vu alors, l'application $Res\circ\Pi$ est injective sur  ${\mathbb C}[\mathfrak{Endo}_{unip-quad}^{red}]$. Il en r\'esulte que la premi\`ere application de la suite ci-dessus est bijective. Comme on l'a dit en 1.5, la compos\'ee $proj_{cusp}\circ Res$ des deux derni\`eres applications est bijective. La suite fournit donc un isomorphisme de ${\mathbb C}[\mathfrak{Endo}_{unip,disc}]$ sur ${\cal R}^{par}_{cusp}$.

L'involution ${\cal F}$ de $\mathfrak{Endo}_{unip-quad}$ pr\'eserve le sous-ensemble $\mathfrak{Endo}_{unip,disc}$. Donc l'involution ${\cal F}$ de ${\mathbb C}[\mathfrak{Endo}_{unip-quad}]$ pr\'eserve le sous-espace ${\mathbb C}[\mathfrak{Endo}_{unip,disc}]$. Via la bijection ci-dessus, on  peut la consid\'erer comme une involution de ${\mathbb C}[Ell_{unip}]$. Le lemme 2.3 entra\^{\i}ne que, pour $\pi\in {\mathbb C}[Ell_{unip}]$, on a l'\'egalit\'e
$$(1) \qquad {\cal F}^{par}\circ proj_{cusp}\circ Res\circ D(\pi)=proj_{cusp}\circ Res\circ D\circ{\cal F}(\pi).$$

 Notons $\mathfrak{St}_{unip-quad}$ le sous-ensemble des $(\lambda,s,1)\in \mathfrak{St}_{tunip}$ tels que $s^2=1$.   Pour $(\lambda,s,1)\in \mathfrak{St}_{unip-quad}$, $s$ d\'etermine une d\'ecomposition $\lambda=\lambda^+\cup \lambda^-$ et le couple $(\lambda^+,\lambda^-)$ appartient \`a ${\cal P}^{symp}_{2}(2n)$. L'application $(\lambda,s,1)\mapsto (\lambda^+,\lambda^-)$ est une bijection de $\mathfrak{St}_{unip-quad}$ sur ${\cal P}^{symp}_{2}(2n)$.
 
 {\bf Notation.} Pour $(\lambda^+,\lambda^-)\in {\cal P}_{2}^{symp}(2n)$, on note $\Pi^{st}(\lambda^+,\lambda^-)$ la repr\'esentation $\Pi(\lambda,s,1)$ o\`u $(\lambda,s,1)$ correspond \`a $(\lambda^+,\lambda^-)$.

Pour tout entier pair $N\in {\mathbb N}$, notons ${\cal P}^{symp,disc}(N)$ l'ensemble des partitions symplectiques de $N$ qui ne contiennent que des termes pairs, lesquels interviennent avec multiplicit\'e au plus $1$. 
Notons $\mathfrak{St}_{unip,disc}$ l'intersection de $\mathfrak{St}_{unip-quad}$ et de  $ \mathfrak{Endo}_{unip,disc}$. L'application $(\lambda,s,1)\mapsto (\lambda^+,\lambda^-)$ se restreint en  une bijection de $\mathfrak{St}_{unip,disc}$ sur  ${\cal P}_{2}^{symp,disc}(2n)$. Pour $(\lambda^+,\lambda^-)\in {\cal P}_{2}^{symp,disc}(2n)$ et pour  
  $(\epsilon^+,\epsilon^-)\in \{\pm 1\}^{Jord_{bp}(\lambda^+)}\times \{\pm 1\}^{Jord_{bp}(\lambda^+)}$, la repr\'esentation $\pi(\lambda^+,\epsilon^+,\lambda^-,\epsilon^-)$ est de la s\'erie discr\`ete. 
  D'apr\`es \cite{MW} proposition 8.2, on a

tout \'el\'ement $\pi\in {\mathbb C}[Ell_{unip}]$ tel que $\pi_{iso}$ soit stable et que $-\pi_{an}$ soit le transfert de $\pi_{iso}$ est combinaison lin\'eaire des  
$\Pi^{st}(\lambda^+,\lambda^-)$ quand $(\lambda^+,\lambda^-)$ d\'ecrit $ {\cal P}_{2}^{symp,disc}(2n)$.

\bigskip

  \subsection{Modules de Jacquet des repr\'esentations stables}
  Soit ${\bf m}=(m_{1}\geq...\geq m_{t}>0)\in {\cal P}(\leq n)$. Posons $n_{0}=n-S({\bf m})$. Pour $\sharp=iso$ ou $an$, la partition ${\bf m}$ d\'etermine un sous-groupe parabolique standard $P_{ \sharp}$  (standard voulant dire qu'il contient le groupe $P_{min}$ d\'efini en 1.7) et sa composante de Levi standard 
  $$M_{\sharp}=GL(m_{1})\times... \times GL(m_{t})\times G_{n_{0},\sharp}.$$
  Dans le cas o\`u $\sharp=an$ et $n_{0}=0$, il n'y a pas de tels sous-groupes et on les \'elimine de ce qui suit.   De m\^eme que l'on a consid\'er\'e simultan\'ement les deux groupes $G_{iso}$ et $G_{an}$, on consid\`ere simultan\'ement les groupes $M_{iso}$ et $M_{an}$. 
  Les notions introduites pr\'ec\'edemment pour les groupes $G_{iso}$ et $G_{an}$ se g\'en\'eralisent aux groupes $M_{iso}$ et $M_{an}$. On ajoute des indices ou exposants $M$ aux objets d\'efinis pour ces groupes.    En particulier, on d\'efinit l'ensemble $Irr_{unip,M_{\sharp}} $ des (classes d'isomorphismes de) repr\'esentations irr\'eductibles de r\'eduction unipotente de $M_{\sharp}(F)$ et on note $Irr_{unip,M}$ la r\'eunion disjointe de $Irr_{unip,M_{iso}}$ et $Irr_{unip,M_{an}}$.  On d\'efinit aussi  le sous-ensemble $Ell_{unip,M}\subset {\mathbb C}[Irr_{unip,M}]$ des repr\'esentations elliptiques. Un \'el\'ement de $Ell_{unip,M}$ est de la forme
  $$(1) \qquad st_{m_{1}}(\vert .\vert_{F}^{z_{1}}\circ det)\otimes...\otimes st_{m_{t}}(\vert .\vert_{F}^{z_{t}}\circ det)\otimes \sigma_{0},$$
  o\`u $z_{1},...,z_{t}$ sont des nombres complexes, $\sigma_{0}\in Ell_{unip, n_{0}}$. On rappelle que, pour tout $m\geq1$, $st_{m}$ est la repr\'esentation de Steinberg de $GL(m;F)$. On sait que tout \'el\'ement de ${\mathbb C}[Irr_{unip,M}]$ est somme d'un \'el\'ement bien d\'etermin\'e de ${\mathbb C}[Ell_{unip,M}]$ et d'une combinaison lin\'eaire d'induites \`a partir de sous-groupes paraboliques propres.

 En identifiant un \'el\'ement de ${\mathbb C}[Irr_{unip,M}]$ \`a sa distribution trace (ou plus exactement \`a la paire de distributions, l'une sur $M_{iso}(F)$, l'autre sur $M_{an}(F)$), on peut d\'efinir la restriction de cette distribution aux \'el\'ements elliptiques (fortement r\'eguliers) de $M(F)$, que l'on note $proj_{ell}$, ou aux \'el\'ements elliptiques et compacts, que l'on note $proj_{ell,comp}$. Remarquons qu'un \'el\'ement elliptique de $M(F)$ est compact si et seulement si ses composantes dans les groupes $GL(m_{j};F)$ sont de d\'eterminants de valeurs absolues $1$. 
 
 Montrons que 

(2) pour $\sigma\in {\mathbb C}[Irr_{unip,M}]$, on a  $proj_{ell,comp}(\sigma)=0$ si et seulement si $proj_{cusp}\circ Res^M(\sigma)=0$.

Les applications $proj_{ell,comp}$ comme $proj_{cusp}\circ Res^M$ annulent les induites \`a partir de sous-groupes paraboliques propres. On peut donc supposer que $\sigma\in {\mathbb C}[Ell_{unip,M}]$. Pour une repr\'esentation $\sigma$ de la forme (1), on voit que  $proj_{ell,comp}(\sigma)$ comme $proj_{cusp}\circ Res^M(\sigma)$ ne d\'ependent pas de $z_{1},...,z_{t}$. Modulo l'intersection des noyaux   de nos applications $proj_{ell,comp}$  et $proj_{cusp}\circ Res^M$, on peut donc supposer $\sigma$ de la forme
  $$st_{m_{1}} \otimes...\otimes st_{m_{t}}\otimes \sigma_{0},$$
o\`u $\sigma_{0}\in Ell_{unip, n_{0}}$. Pour tout $i=1,...t$, $st_{m_{i}}$ n'est annul\'e par aucune des analogues pour $GL(m_{i})$ des deux applications $proj_{ell,comp}$ et    $proj_{cusp}\circ Res$. Les conditions $proj_{ell,comp}(\sigma)=0$, resp.  $proj_{cusp}\circ Res^M(\sigma)=0$, \'equivalent donc \`a $proj_{ell}(\sigma_{0})=0$, resp. $proj_{cusp}\circ Res_{n_{0}}(\sigma_{0})=0$. Parce que $\sigma_{0}$ est elliptique, chacune de  ces conditions \'equivaut \`a $\sigma_{0}=0$. Cela prouve (2). 

Pour $\sharp=iso$ ou $an$ et pour $\pi\in Irr_{unip,\sharp}$, on note $\pi_{M}$ le semi-simplifi\'e du module de Jacquet de $\pi$ relatif au sous-groupe parabolique $P_{\sharp}$. C'est une repr\'esentation de $M_{\sharp}(F)$ dont on sait que toutes ses composantes irr\'eductibles sont de r\'eduction unipotente. Par lin\'earit\'e, on prolonge l'application $\pi\mapsto \pi_{M}$ en  une application lin\'eaire de ${\mathbb C}[Irr_{unip}]$ dans ${\mathbb C}[Irr_{unip,M}]$. Pour $i=1,...,t$, notons $\underline{st}_{m_{i}}$ la restriction de $st_{m_{i}}$ aux \'el\'ements elliptiques et compacts de $GL(m_{i};F)$. Il r\'esulte de ce qui pr\'ec\`ede que, pour $\pi\in {\mathbb C}[Irr_{unip}]$, il existe une unique $\sigma_{0}\in {\mathbb C}[Ell_{unip,n_{0}}]$ de sorte que l'on ait l'\'egalit\'e
$$proj_{ell,comp}(\pi_{M})=\underline{st}_{m_{1}}\otimes...\otimes\underline{st}_{m_{t}}\otimes proj_{ell}(\sigma_{0}).$$
Consid\'erons $(\lambda^+,\lambda^-)\in {\cal P}_{2}^{symp}(2n)$. A $\pi=\Pi^{st}(\lambda^+,\lambda^-)$ correspond ainsi une repr\'esentation $\sigma_{0}\in {\mathbb C}[Ell_{unip,n_{0}}]$. Les propri\'et\'es de stabilit\'e et de transfert se conservent par passage au module de Jacquet et par l'application $proj_{ell,comp}$. On en d\'eduit que $\sigma_{0,iso}$ est stable et que $-\sigma_{0,an}$ en est son transfert au groupe $G_{an,n_{0}}(F)$. D'apr\`es ce que l'on a dit au paragraphe pr\'ec\'edent, $\sigma_{0}$ est combinaison lin\'eaire des $\Pi^{st}(\nu^+,\nu^-)$ quand $(\nu^+,\nu^-)$ d\'ecrit ${\cal P}_{2}^{symp,disc}(2n_{0})$. Il y a donc une unique d\'ecomposition
$$(3) \qquad proj_{ell,comp}(\Pi^{st}(\lambda^+,\lambda^-)_{M})=\sum_{(\nu^+,\nu^-)\in {\cal P}_{2}^{symp,disc}(2n_{0})}c_{{\bf m}}(\lambda^+,\lambda^-;\nu^+,\nu^-)$$
$$\underline{st}_{m_{1}}\otimes...\otimes\underline{st}_{m_{t}}\otimes proj_{ell}(\Pi^{st}(\nu^+,\nu^-)),$$
o\`u les $c_{{\bf m}}(\lambda^+,\lambda^-;\nu^+,\nu^-)$ sont des coefficients complexes. 

 Montrons que
 
 (4) si $\lambda^-=\emptyset$, alors $c_{{\bf m}}(\lambda^+,\emptyset;\nu^+,\nu^-)=0$ pour tout $(\nu^+,\nu^-)\in {\cal P}_{2}^{symp,disc}(2n_{0})$ tel que $\nu^-\not=\emptyset$.
 
Consid\'erons d'abord un couple  $(\lambda^+,\lambda^-)$ quelconque et posons
$$\sigma_{0}=\sum_{(\nu^+,\nu^-)\in {\cal P}_{2}^{symp,disc}(2n_{0})}c_{{\bf m}}(\lambda^+,\lambda^-;\nu^+,\nu^-) \Pi^{st}(\nu^+,\nu^-).$$
 On a d\'ecrit en 2.1 la repr\'esentation de $\tilde{GL}(2n;F)$ correspondant par endoscopie tordue \`a $\Pi^{st}_{iso}(\lambda^+,\lambda^-)$, notons-la  $\tilde{\Pi}^{GL}(\lambda^+,\lambda^-)$.    Au sous-groupe parabolique $P_{iso}$ correspond un sous-groupe parabolique $P^{GL}$ de $GL(2n)$ dont une composante de Levi $M^{GL}$ est de la forme
   $$GL(m_{1})\times..\times GL(m_{t})\times GL(2n_{0})\times GL(m_{t})\times...\times GL(m_{1}).$$
 A $M$ est associ\'e un sous-espace de Levi $\tilde{M}^{GL}\subset \tilde{GL}(2n)$.   La correspondance d\'efinie par l'endoscopie tordue est compatible au passage au module de Jacquet et aux projections sur les \'el\'ements elliptiques et compacts (ces notions \'etant relatives \`a l'espace $\tilde{GL}(2n)$ et non pas au groupe $GL(2n)$). Il en r\'esulte que $proj_{ell,comp}(\Pi^{st}(\lambda^+,\lambda^-)_{M})$ correspond \`a $proj_{ell,comp}(\tilde{\Pi}^{GL}(\lambda^+,\lambda^-)_{\tilde{M}})$, avec des notations compr\'ehensibles. Le module de Jacquet $\tilde{\Pi}^{GL}(\lambda^+,\lambda^-)_{\tilde{M}}$ est somme de prolongements \`a $\tilde{M}^{GL}(F)$ de repr\'esentations  irr\'eductibles de $M(F)$. Celles-ci  interviennent dans le module de Jacquet ordinaire $\Pi^{GL}(\lambda^+,\lambda^-)_{M}$ et  sont de la forme
$$(5) \qquad \rho_{1}\otimes...\otimes \rho_{t}\otimes \rho_{n_{0}}\otimes \check{\rho}_{t}\otimes... \otimes \check{\rho}_{1},$$
 o\`u les $\rho_{j}$ sont des repr\'esentations irr\'eductibles  de $GL(m_{j};F)$, les $\check{\rho}_{j}$ en sont les contragr\'edientes et $\rho_{n_{0}}$ est une repr\'esentation irr\'eductible de $GL(n_{0};F)$ qui  se prolonge en une repr\'esentation $\tilde{\rho}_{n_{0}}$ de l'espace $\tilde{GL}(2n_{0};F)$. Toutes ces repr\'esentations sont de r\'eduction unipotente. Comme pr\'ec\'edemment, il en r\'esulte que, pour $i=1,...,t$, $proj_{ell,comp}(\rho_{i})$ est proportionnelle \`a $ \underline{st}_{m_{i}}$. La correspondance entre   $proj_{ell,comp}(\Pi^{st}(\lambda^+,\lambda^-)_{M})$ et 
 
 \noindent $proj_{ell,comp}(\tilde{\Pi}^{GL}(\lambda^+,\lambda^-)_{\tilde{M}})$ entra\^{\i}ne alors que $proj_{ell}(\sigma_{0})$ correspond \`a une certaine combinaison lin\'eaire des $proj_{ell}(\tilde{\rho}_{n_{0}})$, pour les $\rho_{0}$ intervenant ci-dessus.  De nouveau, pour une telle $\rho_{n_{0}}$, le prolongement  $\tilde{\rho}_{n_{0}}$   est somme d'une repr\'esentation tordue elliptique $\tilde{\rho}_{n_{0},ell}$ et d'une combinaison lin\'eaire d'induites propres. Puisqu'on  projette $\tilde{\rho}_{n_{0}}$ sur les \'el\'ements elliptiques, on peut aussi bien remplacer $\tilde{\rho}_{n_{0}}$ par $\tilde{\rho}_{n_{0},ell}$. D'autre part, on sait d'apr\`es 2.1 que $proj_{ell}(\sigma_{0})$ correspond \`a l'image par $proj_{ell}$ de
 $$\sum_{(\nu^+,\nu^-)\in {\cal P}_{2}^{symp,disc}(2n_{0})}c_{{\bf m}}(\lambda^+,\lambda^-;\nu^+,\nu^-)\tilde{ \Pi}^{GL}(\nu^+,\nu^-).$$
 Les repr\'esentations tordues  intervenant ici sont elliptiques. La projection sur les \'el\'ements elliptiques est injective sur les repr\'esentations tordues elliptiques. Donc la repr\'esentation ci-dessus est combinaison lin\'eaire des $\tilde{\rho}_{n_{0},ell}$ pr\'ec\'edents.   Cela entra\^{\i}ne que, pour tout couple $(\nu^+,\nu^-)$ intervenant avec un coefficient $c_{{\bf m}}(\lambda^+,\lambda^-;\nu^+,\nu^-)$ non nul, la repr\'esentation irr\'eductible sous-jacente $\Pi^{GL}(\nu^+,\nu^-)$ v\'erifie la condition suivante: il existe $\rho_{n_{0}}$ comme ci-dessus tel que   $\Pi^{GL}(\nu^+,\nu^-)$ soit une composante irr\'eductible de la repr\'esentation de $GL(2n_{0};F)$ sous-jacente \`a $\tilde{\rho}_{n_{0},ell}$.
 
  Faisons maintenant l'hypoth\`ese que $\lambda^-=\emptyset$. Alors le support cuspidal de $\Pi^{GL}(\lambda^+,\emptyset)$ est form\'e de caract\`eres du tore d\'eploy\'e maximal de $GL(2n)$ de la forme
 $$\vert .\vert _{F}^{h_{1}/2}\otimes...\otimes \vert .\vert _{F}^{h_{2n}/2},$$
 o\`u  $h_{i}\in {\mathbb Z}$ pour tout $i$. Les repr\'esentations (5) interviennent dans le module de Jacquet de $\Pi^{GL}(\lambda^+,\emptyset)$ donc le support cuspidal de la repr\'esentation $\rho_{n_{0}}$ est de la m\^eme forme, l'entier $n$ \'etant remplac\'e par $n_{0}$.  La  repr\'esentation $\tilde{\rho}_{n_{0},ell}$ se d\'eduit explicitement de  $\tilde{\rho}_{n_{0}}$ par la th\'eorie du quotient de Langlands. Il en r\'esulte que la repr\'esentation sous-jacente de $\tilde{\rho}_{n_{0},ell}$ a m\^eme support cuspidal que $\rho_{n_{0}}$. Cela entra\^{\i}ne que,  pour tout couple $(\nu^+,\nu^-)$ intervenant avec un coefficient $c_{{\bf m}}(\lambda^+,\emptyset;\nu^+,\nu^-)$ non nul, la repr\'esentation irr\'eductible sous-jacente $\Pi^{GL}(\nu^+,\nu^-)$ a pour support cuspidal des caract\`eres du tore d\'eploy\'e maximal de $GL(2n_{0})$ de la forme ci-dessus. Mais, si $\nu^-$ \'etait non vide,   les caract\`eres intervenant dans ce  support cuspidal contiendraient au moins une composante $\vert .\vert _{F}^{h/2}\chi^-$ o\`u $\chi^-$ est le caract\`ere non ramifi\'e de $F^{\times}$ d'ordre $2$. C'est impossible donc $\nu^-$ est vide. Cela prouve (4).

 Pour $\lambda\in {\cal P}^{symp}(2n)$ et $\nu\in {\cal P}^{symp,disc}(2n_{0})$, on pose simplement 
 $$c_{{\bf m}}(\lambda;\nu)=c_{{\bf m}}(\lambda,\emptyset; \nu,\emptyset).$$

 Notons ${\cal I}$ l'ensemble des paires $(I^+,I^-)$ d'ensembles tels que $I^+\sqcup I^-=\{1,...,t\}$. A toute telle paire, on associe les partitions ${\bf m} (I^+)$ et ${\bf m}(I^-)$: ${\bf m}(I^+)$ est form\'ee des $m_{i}$ pour $i\in I^+$ et ${\bf m}(I^-)$ est form\'ee des $m_{i}$ pour $i\in I^-$.
 
 \ass{Lemme}{Soient $(\lambda^+,\lambda^-)\in {\cal P}_{2}^{symp}(2n)$ et $(\nu^+,\nu^-)\in {\cal P}_{2}^{symp,disc}(2n_{0})$.   Alors on a l'\'egalit\'e
 $$c_{{\bf m}}(\lambda^+,\lambda^-;\nu^+,\nu^-)=\sum_{I^+,I^-}c_{{\bf m}(I^+)}(\lambda^+;\nu^+)c_{{\bf m}(I^-)}(\lambda^-;\nu^-),$$
 o\`u l'on somme sur les $(I^+,I^-)\in {\cal I}$ tels que $S(\lambda^+)=S(\nu^+)+2S({\bf m}(I^+))$ et 
$S(\lambda^-)=S(\nu^-)+2S({\bf m}(I^-))$.} 
 
 Preuve. En vertu de (2), l'\'egalit\'e (3) \'equivaut \`a
 $$(6) \qquad proj_{cusp}\circ Res^M(\Pi^{st}(\lambda^+,\lambda^-)_{M})=  \sum_{(\nu^+,\nu^-)\in {\cal P}_{2}^{symp,disc}(2n_{0})}c_{{\bf m}}(\lambda^+,\lambda^-;\nu^+,\nu^-)$$
$$proj_{cusp}\circ Res^M(st_{m_{1}}\otimes...\otimes st_{m_{t}}\otimes \Pi^{st}(\nu^+,\nu^-)).$$
On peut composer cette \'egalit\'e avec la dualit\'e $D^{par,M}$. Cette application  commute avec $proj_{cusp}$. 
On sait d'apr\`es  1.5(1) que $Res^M(\Pi^{st}(\lambda^+,\lambda^-)_{M})=res_{{\bf m}}\circ Res(\Pi^{st}(\lambda^+,\lambda^-))$. On a  aussi $D^{par,M}\circ res_{{\bf m}}=res_{{\bf m}}\circ D^{par}$ et $D^{par}\circ Res=Res\circ D$ d'apr\`es le lemme 1.7.    D'o\`u
$$D^{par,M}\circ proj_{cusp}\circ Res^M(\Pi^{st}(\lambda^+,\lambda^-)_{M})=proj_{cusp}\circ res_{{\bf m}}\circ  Res\circ D(\Pi^{st}(\lambda^+,\lambda^-)).$$
Par d\'efinition, on a
$$\Pi^{st}(\lambda^+,\lambda^-)=\sum_{\epsilon^+,\epsilon^-}\pi(\lambda^+,\epsilon^+,\lambda^-,\epsilon^-),$$
o\`u $(\epsilon^+,\epsilon^-)$ parcourt $\{\pm 1\}^{Jord_{bp}(\lambda^+)}\times \{\pm 1\}^{Jord_{bp}(\lambda^-)}$. Pour un tel couple, la proposition 1.11 dit que 
$$Res\circ D(\pi(\lambda^+,\epsilon^+,\lambda^-,\epsilon^-))=Rep\circ \rho\iota\circ j(\boldsymbol{\rho}_{\lambda^+,\epsilon^+}\otimes \boldsymbol{\rho}_{\lambda^-,\epsilon^-}).$$
Rappelons que $\boldsymbol{\rho}_{\lambda^+,\epsilon^+}$, resp. $\boldsymbol{\rho}_{\lambda^-,\epsilon^-}$, ont \'et\'e identifi\'es \`a des \'el\'ements de ${\cal S}_{n^+}$, resp. ${\cal S}_{n^-}$, o\`u $S(\lambda^+)=2n^+$, resp. $S(\lambda^-)=2n^-$.  D\'efinissons l'\'el\'ement
$$\boldsymbol{\rho}_{\lambda^+}=\sum_{\epsilon^+}\boldsymbol{\rho}_{\lambda^+,\epsilon^+}$$
de ${\cal S}_{n^+}$ et l'\'el\'ement similaire $\boldsymbol{\rho}_{\lambda^-}$ de ${\cal S}_{n^-}$. 
   En utilisant les propri\'et\'es de compatibilit\'e des applications $Rep$ et $Rep^M$, on obtient
$$D^{par,M}\circ proj_{cusp}\circ Res^M(\Pi^{st}(\lambda^+,\lambda^-)_{M})=Rep^M(\psi),$$
o\`u
$$\psi=   proj_{cusp}\circ res_{{\bf m}}\circ \rho\iota\circ j(\boldsymbol{\rho}_{\lambda^+}\otimes \boldsymbol{\rho}_{\lambda^-}).$$

On calcule aussi l'image par $D^{par,M}$ du membre de droite de (6). Pour $i=1,...,t$, l'image par l'involution de Zelevinsky de $st_{m_{i}}$ est la repr\'esentation triviale de $GL(m_{i};F)$. Son image par l'application analogue \`a $Res$ est l'image par l'application analogue \`a $Rep$ de la fonction constante \'egale \`a $1$ sur le groupe $\mathfrak{S}_{m_{i}}$. Rappelons qu'en 1.12, on a identifi\'e ${\mathbb C}[\hat{\mathfrak{S}}_{m_{i}}]_{cusp}$ \`a ${\mathbb C}$ pr\'ecis\'ement en choisissant pour base de cet espace la projection cuspidale de cette fonction constante. En faisant ainsi disparaitre ces espaces ${\mathbb C}[\hat{\mathfrak{S}}_{m_{i}}]_{cusp}$, on obtient que l'image par $D^{par,M}$ du membre de droite de (6) est \'egale \`a $Rep^M(\psi_{0})$, o\`u
$$\psi_{0}= \sum_{(\nu^+,\nu^-)\in {\cal P}_{2}^{symp,disc}(2n_{0})}c_{{\bf m}}(\lambda^+,\lambda^-;\nu^+,\nu^-) proj_{cusp}\circ\rho\iota\circ j_{n_{0}}(\boldsymbol{\rho}_{\nu^+}\otimes \boldsymbol{\rho}_{\nu^-}).$$
 Puisque $Rep^M$ est bijectif, l'\'egalit\'e (6) \'equivaut \`a 
 
 (7) $\psi=\psi_{0}$. 
 
 On a vu que l'\'egalit\'e (6) d\'eterminait les coefficients $c_{{\bf m}}(\lambda^+,\lambda^-;\nu^+,\nu^-)$. Il en est donc de m\^eme de l'\'egalit\'e (7). 
 
 On peut remplacer le couple $(\lambda^+,\lambda^-)$ par $(\lambda^+,\emptyset)$, en rempla\c{c}ant $n$ par $n^+$ et ${\bf m}$ par une partition ${\bf m}^+\in {\cal P}(n^+)$. L'assertion (4) nous dit que les $\nu^-$ intervenant dans la d\'efinition de $\psi_{0}$ disparaissent. En posant $\phi^+=\boldsymbol{\rho}_{\lambda^+}$, $n_{0}^+({\bf m}^+)=n^+-S({\bf m}^+)$ et  
 $$\phi^+[{\bf m}^+]= \sum_{\nu^+\in {\cal P}^{symp,disc}(2n^+_{0}({\bf m}^+))}c_{{\bf m}^+}(\lambda^+,\nu^+)\boldsymbol{\rho}_{\nu^+},$$
 l'analogue de l'\'egalit\'e (7) devient
 $$proj_{cusp}\circ res_{{\bf m}^+}\circ \rho\iota\circ j_{n^+,0}(\phi^+)=proj_{cusp}\circ\rho\iota\circ j_{n_{0}^+({\bf m}^+)}(\phi^+[{\bf m}^+]).$$
 La m\^eme chose vaut pour les donn\'ees affect\'ees d'un exposant $-$. Cela v\'erifie les hypoth\`eses du lemme 1.12. Appliquons ce lemme. Il affirme une \'egalit\'e dont le membre de gauche n'est autre que $\psi$. On voit que le membre de droite est \'egal \`a
  $$\psi'_{0}=\sum_{(\nu^+,\nu^-)\in {\cal P}_{2}^{symp,disc}(2n_{0})}c'_{{\bf m}}(\lambda^+,\lambda^-;\nu^+,\nu^-) proj_{cusp}\circ\rho\iota\circ j_{n_{0}}(\boldsymbol{\rho}_{\nu^+}\otimes \boldsymbol{\rho}_{\nu^-}),$$
  o\`u $c'_{{\bf m}}(\lambda^+,\lambda^-;\nu^+,\nu^-)$ est le membre de droite de l'\'egalit\'e du pr\'esent \'enonc\'e. On a donc $\psi=\psi'_{0}$. Puisqu'on a aussi $\psi=\psi_{0}$ et que cette \'egalit\'e d\'etermine les coefficients, on conclut $c_{{\bf m}}(\lambda^+,\lambda^-;\nu^+,\nu^-) =c'_{{\bf m}}(\lambda^+,\lambda^-;\nu^+,\nu^-) $ pour tout couple $(\nu^+,\nu^-)$, ce qui d\'emontre le lemme. $\square$
  
  \bigskip
  
  \subsection{Modules de Jacquet des repr\'esentations endoscopiques}
  On a vu en 2.4 que $\mathfrak{St}_{unip-quad}$ s'identifiait \`a ${\cal P}^{symp}_{2}(2n)$. L'ensemble $\mathfrak{Endo}_{unip-quad}$ s'identifie de m\^eme \`a ${\cal P}^{symp}_{4}(2n)$: \`a $(\lambda,s,h)\in \mathfrak{Endo}_{unip-quad}$, on associe le quadruplet de partitions $(\lambda^{++},\lambda^{-+},\lambda^{+-},\lambda^{--})$ o\`u, pour $\zeta,\xi=\pm$ et tout entier $i\geq1$, la multiplicit\'e de $i$ dans $\lambda^{\zeta\xi}$ est $m_{s,h}^{\zeta\xi}(i)$, cf. 2.2. Remarquons que le sous-ensemble $\mathfrak{Endo}_{unip,disc}$ de $\mathfrak{Endo}_{unip-quad}$ s'identifie \`a ${\cal P}^{symp,disc}_{4}(2n)$. 
  
  Soit ${\bf m}=(m_{1},...,m_{t}>0)\in {\cal P}(\leq n)$. On reprend les constructions  et notations du paragraphe pr\'ec\'edent associ\'ees \`a cette partition. Soient $(\lambda,s,h)\in \mathfrak{Endo}_{unip-quad}$ et $(\nu,z,k)\in \mathfrak{Endo}_{unip,disc,n_{0}}$. On d\'eduit de ces donn\'ees des quadruplets $(\lambda^{++},\lambda^{-+},\lambda^{+-},\lambda^{--})$ et $(\nu^{++},\nu^{-+},\nu^{+-},\nu^{--})$. Consid\'erons un quadruplet ${\bf I}=(I^{++},I^{-+},I^{+-},I^{--})$ d'ensembles tels que
   $$I^{++}\sqcup I^{-+}\sqcup I^{+-}\sqcup I^{--}=\{1,...,t\}.$$
   On en d\'eduit des partitions ${\bf m}^{\zeta\xi}({\bf I})$ pour $\zeta,\xi=\pm$, form\'ees des $m_{i}$ pour $i\in I^{\zeta\xi}$. Notons ${\cal I}$ l'ensemble des quadruplets comme ci-dessus qui v\'erifient
   $$S(\lambda^{\zeta\xi})=S(\nu^{\zeta\xi})+2S({\bf m}^{\zeta\xi}({\bf I}))$$
   pour tous $\zeta,\xi$. Pour ${\bf I}=(I^{++},I^{-+},I^{+-},I^{--})\in {\cal I}$ et pour $\zeta,\xi=\pm$, on a d\'efini en 2.5 le coefficient $c_{{\bf m}^{\zeta\xi}({\bf I})}(\lambda^{\zeta\xi};\nu^{\zeta\xi})$. Posons
   $$c_{{\bf I}}(\lambda,s,h;\nu,z,k)=\prod_{\zeta,\xi=\pm}c_{{\bf m}^{\zeta\xi}({\bf I})}(\lambda^{\zeta\xi};\nu^{\zeta\xi}),$$
   puis
   $$c(\lambda,s,h;\nu,z,k)=\sum_{{\bf I}\in {\cal I}}c_{{\bf I}}(\lambda,s,h;\nu,z,k).$$
   
   Remarquons la propri\'et\'e suivante
   
   (1) on a l'\'egalit\'e $c(\lambda,h,s;\nu,k,z)=c(\lambda,s,h;\nu,z,k)$.
   
   En effet, \'echanger $s$ et $h$, resp. $z$ et $k$, revient \`a \'echanger $\lambda^{-+}$ et $\lambda^{+-}$, resp.  $\nu^{-+}$ et $\nu^{+-}$. L'\'echange de $I^{-+}$ et $I^{+-}$ dans la d\'efinition ci-dessus conduit \`a l'\'egalit\'e (1).
   
   \ass{Lemme}{Supposons v\'erifi\'e le th\'eor\`eme 2.1. Soit $(\lambda,s,h)\in \mathfrak{Endo}_{unip-quad}$. Alors on a l'\'egalit\'e
   $$proj_{ell,comp}(\Pi(\lambda,s,h)_{M})=\sum_{(\nu,z,k)\in \mathfrak{Endo}_{unip,disc,n_{0}}}c(\lambda,s,h;\nu,z,k) \underline{st}_{m_{1}}\otimes...\otimes\underline{st}_{m_{t}}\otimes proj_{ell}(\Pi(\nu,z,k)).$$}
   
   Preuve. On introduit le quadruplet de partitions $(\lambda^{++},\lambda^{-+},\lambda^{+-},\lambda^{--})$  associ\'e \`a $(\lambda,s,h)$. Soit $\sharp=iso$ ou $an$. Le th\'eor\`eme 2.1 dit que
   $$\Pi(\lambda,s,h)_{\sharp}=sgn_{\sharp}transfert(\Pi_{iso}^{st}(\lambda^{++},\lambda^{-+})\otimes \Pi_{iso}^{st}(\lambda^{+-},\lambda^{--})),$$
   o\`u  on  a pos\'e $sgn_{iso}=1$ et $sgn_{an}=-1$.
   On d\'efinit les entiers $m^+$ et $m^-$ de sorte que $\Pi_{iso}^{st}(\lambda^{++},\lambda^{-+})$ soit une repr\'esentation de $G_{m^+,iso}(F)$ et $\Pi_{iso}^{st}(\lambda^{+-},\lambda^{--})$ soit une repr\'esentation de $G_{m^-,iso}(F)$. Le transfert est compatible au passage au module de Jacquet. On doit toutefois consid\'erer tous les sous-groupes de Levi de $G_{m^+,iso}\times G_{m^-,iso}$ qui se transf\`erent en notre Levi $M$. Pr\'ecis\'ement,  consid\'erons une paire d'ensembles $(J^+,J^-)$ telle que $J^+\cup J^-=\{1,...,t\}$. On lui associe comme toujours  deux partitions ${\bf m}(J^+)$ et ${\bf m}(J^-)$.  Notons ${\cal J}$ l'ensemble de ces paires telles que $S({\bf m}(J^+))\leq m^+$, $S({\bf m}(J^-))\leq m^-$. Une telle paire d\'efinit des groupes de Levi $M(J^+)\subset G_{m^+,iso}$ et $M(J^-)\subset G_{m^-,iso}$  associ\'es aux partitions ${\bf m}(J^+)$ et ${\bf m}(J^-)$. On a un isomorphisme
$$M(J^+)\times M(J^-)=GL(m_{1})\times...\times GL(m_{t})\times G_{m_{0}(J^+),iso}\times G_{m_{0}(J^-),iso},$$
o\`u $m_{0}(J^+)=m^+-S({\bf m}(J^+))$ et $m_{0}(J^-)=m^--S({\bf m}(J^-))$. Le groupe $G_{m_{0}(J^+),iso}\times G_{m_{0}(J^-),iso}$ est celui d'une donn\'ee endoscopique \'evidente de $G_{n_{0},\sharp}$. On dispose pour ces objets d'un transfert. En le tensorisant par les identit\'es sur les groupes $GL(m_{j})$, on obtient un transfert entre  distributions invariantes stables sur  $M(J^+)(F)\times M(J^-)(F)$ et  distributions invariantes sur $M(F)$. On a pr\'ecis\'ement
$$ \Pi(\lambda,s,h)_{\sharp,M_{\sharp}}=sgn_{\sharp}\sum_{(J^+,J^-)\in {\cal J}}transfert(\Pi^{st}_{iso}(\lambda^{++},\lambda^{-+})_{M(J^+)}\otimes \Pi^{st}_{iso}(\lambda^{+-},\lambda^{--})_{M(J^-)}).$$
Le transfert commute aussi aux projections sur les elliptiques compacts. Donc
$$(2) \qquad proj_{ell,comp}( \Pi(\lambda,s,h)_{\sharp,M_{\sharp}})=\sum_{(J^+,J^-)\in {\cal J}}sgn_{\sharp}$$
$$transfert(proj_{ell,comp}(\Pi^{st}_{iso}(\lambda^{++},\lambda^{-+})_{M(J^+)})\otimes proj_{ell,comp}(\Pi^{st}_{iso}(\lambda^{+-},\lambda^{--})_{M(J^-)})).$$
Fixons $(J^+,J^-)\in {\cal J}$. On applique la relation 2.5(3):
$$proj_{ell,comp}(\Pi^{st}_{iso}(\lambda^{++},\lambda^{-+})_{M(J^+)})\otimes proj_{ell,comp}(\Pi^{st}_{iso}(\lambda^{+-},\lambda^{--})_{M(J^-)})=$$
$$\sum_{(\nu^{++},\nu^{-+})\in {\cal P}_{2}^{symp,disc}(2m_{0}(J^+))}\sum_{(\nu^{+-},\nu^{--})\in {\cal P}_{2}^{symp,disc}(2m_{0}(J^-))}c_{{\bf m}(J^+)}(\lambda^{++},\lambda^{-+};\nu^{++},\nu^{-+})$$
$$c_{{\bf m}(J^-)}(\lambda^{+-},\lambda^{--};\nu^{+-},\nu^{--})\underline{st}_{m_{1}}\otimes ...\otimes \underline{st}_{m_{t}}\otimes proj_{ell}(\Pi^{st}_{iso}(\nu^{++},\nu^{-+}))\otimes proj_{ell}(\Pi^{st}_{iso}(\nu^{+-},\nu^{--})).$$
Les quadruplets $(\nu^{++},\nu^{-+},\nu^{+-},\nu^{--})$ qui interviennent dans cette formule appartiennent tous \`a ${\cal P}^{symp,disc}_{4}(2n_{0})$. La relation (2) se r\'ecrit
$$(3) \qquad proj_{ell,comp}( \Pi(\lambda,s,h)_{\sharp,M_{\sharp}})=\sum_{(\nu^{++},\nu^{-+},\nu^{+-},\nu^{--})\in {\cal P}^{symp,disc}_{4}(2n_{0})} C(\nu^{++},\nu^{-+},\nu^{+-},\nu^{--}) $$
$$sgn_{\sharp}\underline{st}_{m_{1}}\otimes ...\otimes \underline{st}_{m_{t}}\otimes \, sgn_{\sharp} transfert(proj_{ell}(\Pi^{st}_{iso}(\nu^{++},\nu^{-+}))\otimes proj_{ell}(\Pi^{st}_{iso}(\nu^{+-},\nu^{--}))),$$
o\`u  la constante $C(\nu^{++},\nu^{-+},\nu^{+-},\nu^{--})$ est d\'efinie de la fa\c{c}on suivante. Le quadruplet $(\nu^{++},\nu^{-+},\nu^{+-},\nu^{--})$ \'etant fix\'e, on note ${\cal J}^*$ le sous-ensemble des $(J^+,J^-)\in {\cal J}$ tels que $S(\nu^{++})+S(\nu^{-+})=2m_{0}(J^+)$ et $S(\nu^{+-})+S(\nu^{--})=2m_{0}(J^-)$. Alors
$$(4) \qquad C(\nu^{++},\nu^{-+},\nu^{+-},\nu^{--})=\sum_{(J^+,J^-)\in {\cal J}^* }c_{{\bf m}(J^+)}(\lambda^{++},\lambda^{-+};\nu^{++},\nu^{-+})$$
$$c_{{\bf m}(J^-)}(\lambda^{+-},\lambda^{--};\nu^{+-},\nu^{--}).$$
 Pour $(\nu^{++},\nu^{-+},\nu^{+-},\nu^{--})\in {\cal P}^{symp,disc}_{4}(2n_{0})$, les m\^emes arguments que ci-dessus (en particulier, on utilise le th\'eor\`eme 2.1) montrent que
 $$ sgn_{\sharp} transfert(proj_{ell}(\Pi^{st}_{iso}(\nu^{++},\nu^{-+}))\otimes proj_{ell}(\Pi^{st}_{iso}(\nu^{+-},\nu^{--})))=proj_{ell}(\Pi(\nu,z,k)_{\sharp}),$$
  o\`u $(\nu,z,k)\in \mathfrak{Endo}_{unip,disc,n_{0}}$ correspond \`a  $(\nu^{++},\nu^{-+},\nu^{+-},\nu^{--})$. La formule (3) devient celle de l'\'enonc\'e \`a condition de d\'emontrer l'\'egalit\'e
 $$(5) \qquad c(\lambda,s,h;\nu,z,k)=C(\nu^{++},\nu^{-+},\nu^{+-},\nu^{--}).$$
 Fixons donc $(\nu,z,k)$ et le quadruplet $(\nu^{++},\nu^{-+},\nu^{+-},\nu^{--})$ qui lui correspond. On applique le lemme 2.5 aux termes intervenant dans la d\'efinition (4). On obtient
 $$(6) \qquad C(\nu^{++},\nu^{-+},\nu^{+-},\nu^{--})=\sum_{(J^+,J^-)\in {\cal J}^* }\sum_{I^{++},I^{-+}}c_{{\bf m}(J^+)(I^{++})}(\lambda^{++};\nu^{++})$$
 $$c_{{\bf m}(J^+)(I^{-+})}(\lambda^{-+};\nu^{-+})\sum_{I^{+-},I^{--}}c_{{\bf m}(J^-)(I^{+-})}(\lambda^{+-};\nu^{+-})c_{{\bf m}(J^-)(I^{--})}(\lambda^{--};\nu^{--}).$$
 On v\'erifie que les quadruplets ${\bf I}=(I^{++},I^{-+},I^{+-},I^{--})$ intervenant dans cette formule appartiennent \`a l'ensemble ${\cal I}$ d\'efini avant l'\'enonc\'e. Un \'el\'ement ${\bf I}=(I^{++},I^{-+},I^{+-},I^{--})\in {\cal I}$ intervient dans la somme ci-dessus pour un unique $(J^+,J^-)$: on a $J^+=I^{++}\cup I^{-+}$ et $J^-=I^{+-}\cup I^{--}$. On a aussi ${\bf m}(J^+)(I^{++})={\bf m}(I^{++})$ etc...
Alors le membre de droite de (6) devient exactement  $ c(\lambda,s,h;\nu,z,k)$. Cela prouve (5) et ach\`eve la d\'emonstration. $\square$. 

\bigskip

\subsection{Egalit\'e des deux involutions}
On a introduit deux involutions ${\cal F}^{par}$ et $\mathfrak{F}^{par}$ de l'espace ${\cal R}^{par,glob}$, cf. 1.9 et 2.3.

\ass{Th\'eor\`eme}{Supposons v\'erifi\'e le th\'eor\`eme 2.1. Alors on a l'\'egalit\'e $\mathfrak{F}^{par}={\cal F}^{par}$.}

Preuve. En vertu de l'isomorphisme 1.5(2), il suffit de d\'emontrer que, pour toute partition ${\bf m}\in {\cal P}(\leq n)$, on a l'\'egalit\'e

(1) $proj_{cusp}\circ res_{{\bf m}}\circ \mathfrak{F}^{par}=proj_{cusp}\circ res_{{\bf m}}\circ {\cal F}^{par}$.

Fixons donc ${\bf m}=(m_{1},...,m_{t})$ et posons $n_{0}=n-S({\bf m})$.  Composons \`a droite les deux membres de l'\'egalit\'e ci-dessus par l'application $Res\circ D\circ \Pi$. On obtient des applications lin\'eaires d\'efinies sur ${\mathbb C}[\mathfrak{Endo}_{unip-quad}]$. 
En vertu du lemme 2.2 et de la remarque qui le suit,  il suffit de d\'emontrer que ces applications sont \'egales. Soit $(\lambda,s,h)\in \mathfrak{Endo}_{unip-quad}$. Par d\'efinition de $\mathfrak{F}^{par}$, on a l'\'egalit\'e
$$proj_{cusp}\circ res_{{\bf m}}\circ \mathfrak{F}^{par}\circ Res\circ D(\Pi(\lambda,s,h))=proj_{cusp}\circ res_{{\bf m}}\circ  Res\circ D(\Pi(\lambda,h,s)).$$
On a 
$$proj_{cusp}\circ res_{{\bf m}}\circ  Res\circ D=proj_{cusp}\circ res_{{\bf m}}\circ D^{par}\circ Res$$
d'apr\`es le lemme 1.7 puis
$$proj_{cusp}\circ res_{{\bf m}}\circ D^{par}\circ Res=D^{par,M}\circ proj_{cusp}\circ res_{{\bf m}}\circ  Res$$
d'apr\`es les propri\'et\'es de compatibilit\'e de l'application $D^{par}$. Enfin, 
$$res_{{\bf m}}\circ  Res(\Pi(\lambda,h,s))=Res^M(\Pi(\lambda,h,s)_{M})$$
d'apr\`es 1.5(1).
D'o\`u
 $$ proj_{cusp}\circ res_{{\bf m}}\circ \mathfrak{F}^{par}\circ Res\circ D(\Pi(\lambda,s,h))=D^{par,M}\circ proj_{cusp}\circ    Res^M(\Pi(\lambda,h,s)_{M}).$$
En appliquant 2.5(2) et  le lemme 2.6, on obtient
$$(2) \qquad  proj_{cusp}\circ res_{{\bf m}}\circ \mathfrak{F}^{par}\circ Res\circ D(\Pi(\lambda,s,h))=\sum_{(\nu,z,k)\in \mathfrak{Endo}_{unip,disc,n_{0}}}c(\lambda,h,s;\nu,z,k)X(\nu,z,k),$$
o\`u
$$X(\nu,z,k)=D^{par,M}\circ proj_{cusp}\circ    Res^M(st_{m_{1}}\otimes...\otimes st_{m_{t}}\otimes \Pi(\nu,z,k)).$$

Calculons maintenant $proj_{cusp}\circ res_{{\bf m}}\circ {\cal F}^{par}\circ Res\circ D(\Pi(\lambda,s,h))$. Dans l'espace ${\cal R}_{{\bf m}}^{par,glob}$, cf. 1.5, on d\'efinit les involutions ${\cal F}^{par}_{{\bf m}}$, resp. $\mathfrak{F}^{par}_{{\bf m}}$, qui sont les produits tensoriels des identit\'es des composantes $C^{GL(m_{i})}$ et de ${\cal F}^{par}_{n_{0}}$, resp. $\mathfrak{F}^{par}_{n_{0}}$, sur la composante ${\cal R}_{n_{0}}^{par,glob}$. Les propri\'et\'es de compatibilit\'e des involutions de Lusztig entra\^{\i}nent
 $$proj_{cusp}\circ res_{{\bf m}}\circ {\cal F}^{par}\circ Res\circ D(\Pi(\lambda,s,h))={\cal F}_{{\bf m}}^{par}\circ proj_{cusp}\circ res_{{\bf m}}\circ   Res\circ D(\Pi(\lambda,s,h)),$$
 puis, par les m\^emes arguments que ci-dessus,
 $$proj_{cusp}\circ res_{{\bf m}}\circ {\cal F}^{par}\circ Res\circ D(\Pi(\lambda,s,h))={\cal F}_{{\bf m}}^{par}\circ   D^{par,M}\circ proj_{cusp}\circ Res^M(\Pi(\lambda,s,h)_{M}).$$
 En appliquant 2.5(2) et  le lemme 2.6, on obtient
 $$(3) \qquad proj_{cusp}\circ res_{{\bf m}}\circ {\cal F}^{par}\circ Res\circ D(\Pi(\lambda,s,h))=\sum_{(\nu,z,k)\in \mathfrak{Endo}_{unip,disc,n_{0}}}c(\lambda,s,h;\nu,z,k)Y(\nu,z,k),$$
o\`u
$$Y(\nu,z,k)={\cal F}_{{\bf m}}^{par}\circ D^{par,M}\circ proj_{cusp}\circ    Res^M(st_{m_{1}}\otimes...\otimes st_{m_{t}}\otimes \Pi(\nu,z,k)).$$

Fixons $(\nu,z,k)\in \mathfrak{Endo}_{unip,disc,n_{0}}$. On a $D^{par,M}\circ proj_{cusp}\circ Res^M=proj_{cusp}\circ Res^M\circ D^M$. On a aussi ${\cal F}^{par}_{{\bf m}}\circ proj_{cusp}=proj_{cusp}\circ \mathfrak{F}_{{\bf m}}^{par}$ d'apr\`es le lemme 2.3. D'o\`u
$$Y(\nu,z,k)=proj_{cusp}\circ \mathfrak{F}_{{\bf m}}^{par}\circ Res^M\circ D^M(st_{m_{1}}\otimes...\otimes st_{m_{t}}\otimes \Pi(\nu,z,k)).$$
En appliquant la d\'efinition de l'involution $\mathfrak{F}_{{\bf m}}^{par}$,  on obtient
 $$Y(\nu,z,k)=proj_{cusp}\circ  Res^M\circ D^M(st_{m_{1}}\otimes...\otimes st_{m_{t}}\otimes \Pi(\nu,k,z)),$$
 puis, par les m\^emes arguments de compatibilit\'e,
 $$Y(\nu,z,k)=D^{par,M}\circ proj_{cusp}\circ  Res^M (st_{m_{1}}\otimes...\otimes st_{m_{t}}\otimes \Pi(\nu,k,z))$$
 Autrement dit, $Y(\nu,z,k)=X(\nu,k,z)$. En permutant $z$ et $k$ dans le membre de droite de (3), on obtient  
 $$(4) \qquad  proj_{cusp}\circ   res_{{\bf m}}\circ {\cal F}^{par}\circ Res\circ D(\Pi(\lambda,s,h))=\sum_{(\nu,z,k)\in \mathfrak{Endo}_{unip,disc,n_{0}}}c(\lambda,s,h;\nu,k,z)X(\nu,z,k).$$
 La relation 2.6(1) entra\^{\i}ne que les membres de droite de (2) et (4) sont \'egaux. Donc aussi les membres de gauche. Comme on l'a dit, cela d\'emontre (1). 
  $\square$
 
 \bigskip

 {\bf Index des notations}
 
${\mathbb C}[X]$ 1.4; $C'_{n'}$ 1.5; $C^{''\pm}_{n'',\sharp}$ 1.5; $C''_{n''}$ 1.5; $C^{GL(m)}$ 1.5; ${\mathbb C}[\hat{W}_{N}]_{cusp}$ 1.8; $D(n)$ 1.2; $D_{iso}(n)$ 1.2;  $D_{an}(n)$ 1.2; $D$ 1.7; $D^{par}$ 1.7; $\eta(Q)$ 1.1; $\eta^+(Q)$ 1.1; $\eta^-(Q)$ 1.1; $Ell_{unip}$ 1.4; $\mathfrak{Ell}_{unip}$ 1.4; $\mathfrak{Endo}_{tunip}$ 2.1; $\mathfrak{Endo}_{unip-quad}$ 2.2; $\mathfrak{Endo}_{unip-quad}^{red}$ 2.2; $\mathfrak{Endo}_{unip,disc}$ 2.4;  ${\cal F}^L$ 1.9; ${\cal F}^{par}$ 1.9; ${\cal F}$ 2.3: $\mathfrak{F}^{par}$ 2.3; $G_{iso}$ 1.1; $G_{an}$ 1.1; $\Gamma$ 1.8; $\boldsymbol{\Gamma}$ 1.8; $\tilde{GL}(2n)$ 2.1; $Irr_{tunip}$ 1.3; $\mathfrak{Irr}_{tunip}$ 1.3; $Irr_{unip-quad}$1.3; $\mathfrak{Irr}_{unip-quad}$ 1.3; $Jord(\lambda)$ 1.3; $Jord_{bp}(\lambda)$ 1.3; $Jord_{bp}^{k}(\lambda)$ 1.4; $K_{n',n''}^{\pm}$ 1.2; $k$ 1.9; $L^*$ 1.1; $L_{n',n''}$ 1.2; $l(\lambda)$ 1.3; $mult_{\lambda}$ 1.3; $\mathfrak{o}$ 1.1; $O^+(Q)$ 1.1; $O^-(Q)$ 1.1; $\varpi$ 1.1; $\pi_{n',n''}$ 1.3; ${\cal P}(N)$ 1.3; ${\cal P}^{symp}(2N)$ 1.3; $\boldsymbol{{\cal P}^{symp}}(2N)$ 1.3; $\pi(\lambda,s,\epsilon)$ 1.3; $\pi(\lambda^+,\epsilon^+,\lambda^-,\epsilon^-)$ 1.3; $\pi_{ell}(\lambda^+,\epsilon^+,\lambda^-,\epsilon^-)$ 1.4; $proj_{cusp}$ 1.5; ${\cal P}(\leq n)$ 1.5; ${\cal P}_{k}(N)$ 1.8; $\Pi(\lambda,s,h)$ 2.1; $\Pi^{st}(\lambda^+,\lambda^-)$ 2.4; ${\cal P}^{symp,disc}(2n)$ 2.4; $Q_{iso}$ 1.1; $Q_{an}$ 1.1; $\rho_{\lambda}$ 1.3;  ${\cal R}^{par}$ 1.5; ${\cal R}^{par,glob}$ 1.5; ${\cal R}^{par}_{cusp}$ 1.5; ${\cal R}^{par,glob}_{{\bf m}}$ 1.5; ${\cal R}^{par}_{{\bf m},cusp}$ 1.5; $res'_{m}$ 1.5; $res''_{m}$ 1.5; $res_{m}$ 1.5 et 1.8; $res_{{\bf m}}$ 1.5; ${\cal R}$ 1.8; ${\cal R}(\gamma)$ 1.8; ${\cal R}(\boldsymbol{\gamma})$ 1.8; ${\cal R}^{glob}$ 1.8; ${\cal R}_{cusp}$ 1.8; $Rep$ 1.9; $\rho\iota$ 1.10; $S(\lambda)$ 1.3; $\mathfrak{S}_{N}$ 1.8; $\hat{\mathfrak{S}}_{N}$ 1.8; $sgn$ 1.8; $sgn_{CD}$ 1.8; ${\cal S}_{n}$ 1.11; $\mathfrak{St}_{tunip}$ 2.1; $\mathfrak{St}_{unip-quad}$ 2.4; $\mathfrak{St}_{unip,disc}$ 2.4;  $sgn_{iso}$ 2.6; $sgn_{an}$ 2.6; $val_{F}$ 1.1; $V_{iso}$ 1.1; $V_{an}$ 1.1; $W_{N}$ 1.8; $\hat{W}_{N}$ 1.8; $w_{\alpha}$ 1.8; $w_{\alpha,\beta}$ 1.8; $w_{\alpha,\beta',\beta''}$ 1.8; $Z(\lambda)$ 1.3; $Z(\lambda,s)$ 1.3; ${\bf Z}(\lambda,s)$ 1.3; ${\bf Z}(\lambda,s)^{\wedge}$ 1.3; $\vert .\vert _{F}$ 1.1.

jean-loup.waldspurger@imj-prg.fr

\end{document}